\documentclass[preprint,11pt]{elsarticle} 
\usepackage[bookmarks=true,colorlinks=true,linkcolor=blue]{hyperref}

\pdfoutput=1

\biboptions{sort&compress}

\newcommand{\wdh}[1]{{\color{blue}#1}}

\usepackage{soul}
 
\usepackage[usenames,dvipsnames,svgnames,table]{xcolor}

\usepackage{amssymb,amsmath}
\usepackage{amsthm}

\usepackage{algorithm} 
\usepackage{algpseudocode} 

\usepackage{calc}
\usepackage[margin=1.in]{geometry}

\renewcommand{\url}[1]{}
\newcommand{\citeCount}[1]{}
\newlength{\ycbTop}
\newlength{\ycbMid}%

\newcommand{\f}[2]{\frac{#1}{#2}}
\newtheorem*{AMPFluidImpedence}{{AMP Fluid Impedence}}

\newcommand{\eqdef}{\overset{{\rm def}}{=}}

\newcommand{\dxs}{{\Delta \bar{x}}}

\newcommand{\dt}{\Delta t}

\def\ba#1\ea{\begin{align}#1\end{align}}

\def\bas#1\eas{\begin{align*}#1\end{align*}}

\def\bat#1\eat{\begin{alignat}{3}#1\end{alignat}}

\def\bats#1\eats{\begin{alignat*}{3}#1\end{alignat*}}

\newcommand{\bse}{\begin{subequations}}
\newcommand{\ese}{\end{subequations}}

\newcommand{\Lvh}{\Lv_h}
\newcommand{\Nvh}{\Nv_h}
    
\newcommand{\bogus}[1]{{}}

\newcommand{\dv}{\mathbf{ d}}
\newcommand{\ev}{\mathbf{ e}}

\newcommand{\iv}{\mathbf{ i}}

\newcommand{\nv}{\mathbf{ n}}

\newcommand{\qv}{\mathbf{ q}}
\newcommand{\rv}{\mathbf{ r}}

\newcommand{\tv}{\mathbf{ t}}
\newcommand{\uv}{\mathbf{ u}}
\newcommand{\vv}{\mathbf{ v}}

\newcommand{\xv}{\mathbf{ x}}

\newcommand{\Fv}{\mathbf{ F}}
\newcommand{\Gv}{\mathbf{ G}}
\newcommand{\Hv}{\mathbf{ H}}
\newcommand{\Iv}{\mathbf{ I}}

\newcommand{\Lv}{\mathbf{ L}}

\newcommand{\Nv}{\mathbf{ N}}

\newcommand{\Vv}{\mathbf{ V}}

\newcommand{\Real}{{\mathbb R}}

\newcommand{\Ac}{{\mathcal A}}
\newcommand{\Bc}{{\mathcal B}}
\newcommand{\Cc}{{\mathcal C}}

\newcommand{\Gc}{{\mathcal G}}

\newcommand{\Nc}{{\mathcal N}}

\newcommand{\tauv}{\boldsymbol{\tau}}

\newcommand{\sigmav}{\boldsymbol{\sigma}}

\newcommand{\grad}{\nabla}
\newcommand{\curl}{\times}
\newcommand{\grads}{\nabla_{\xsv}}

\newcommand{\tableFont}{\scriptsize}

\newcommand{\num}[2]{#1e#2} 

\newcommand{\rateLabel}{rate}

\clearpage

\newcommand{\tnv}{\tv}

\newcommand{\rhos}{\bar{\rho}}
\newcommand{\cp}{\bar{c}_p}
\newcommand{\cs}{\bar{c}_s}

\newcommand{\zs}{{\bar {z}}}  
\newcommand{\zp}{z_p}  
\newcommand{\us}{\bar{u}}
\newcommand{\rs}{\bar{r}}

\newcommand{\vs}{\bar{v}}
\newcommand{\sigmas}{\bar{\sigma}}

\newcommand{\xs}{\bar{x}}
\newcommand{\xsv}{\bar{\xv}}
\newcommand{\ys}{{\bar{y}}}

\newcommand{\zf}{z_f}

\newcommand{\Usv}{\bar{\uv}}
\newcommand{\Vsv}{\bar{\vv}}
\newcommand{\Fsv}{\bar{\Fv}}

\newcommand{\Qsv}{\bar{\qv}}

\newcommand{\Sigmasv}{\bar{\sigmav}}

\newcommand{\Grad}{\grad_h}
\newcommand{\Grads}{\bar{\grad}_h}

\newcommand{\isv}{{\bar{\iv}}}
\newcommand{\Sigmav}{{\sigmav}}

\newcommand{\usv}{\bar{\uv}}
\newcommand{\vsv}{\bar{\vv}}

\newcommand{\nsv}{\bar{\nv}}

\newcommand{\uvs}{\bar{\uv}}
\newcommand{\xvs}{\bar{\xv}}

\newcommand{\sigmasv}{\bar{\sigmav}}

\newcommand{\normalss}{\sffamily}

\newcommand{\Hs}{\bar{H}}

\newcommand{\lambdas}{\bar{\lambda}}
\newcommand{\mus}{\bar{\mu}}

\newcommand{\OmegaF}{\Omega}
\newcommand{\OmegaS}{{\bar\Omega}_0}

\newcommand{\OmegaFh}{\Omega_h}
\newcommand{\OmegaSh}{{\bar\Omega}_h}

\newcommand{\GammaIs}{{\bar\Gamma}_0}
\newcommand{\GammaIf}{\Gamma}
\newcommand{\GammaIfh}{\Gamma_h}
\newcommand{\GammaIsh}{{\bar\Gamma}_h}
\newcommand{\GammaIFh}{\Gamma_h}

\renewcommand{\zp}{\bar{z}_p}
\renewcommand{\zs}{\bar{z}_s}

\newcommand{\strutt}{\rule{0pt}{10pt}}

\newcommand{\nd}{{n_d}} 

\newcommand{\scf}{\delta}

\newcommand{\Gcd}{\Gc_d}

\newcommand{\Omegas}{\bar\Omega}

\usepackage{tikz}
\usetikzlibrary{calc}

\newlength{\tfwidth}
\newlength{\tfheight}
\newlength{\tfxa}
\newlength{\tfxb}
\newlength{\tfya}
\newlength{\tfyb}
\newcommand{\trimFigWithBox}[6]{%
\setlength\fboxsep{0pt}%
\setlength\fboxrule{1.0pt}
\fbox{\includegraphics[width=#2, clip, trim=#3 #4 #5 #6]{#1}}%
}
\newcommand{\trimFigNoBox}[6]{%
\setlength\fboxsep{1pt}
\setlength\fboxrule{0.0pt}
\fbox{\includegraphics[width=#2, clip, trim=#3 #4 #5 #6]{#1}}%
}
\newcommand{\trimFigHeightWithBox}[6]{%
\setlength\fboxsep{0pt}%
\setlength\fboxrule{1.0pt}
\fbox{\includegraphics[height=#2, clip, trim=#3 #4 #5 #6]{#1}}%
}
\newcommand{\trimFigHeightNoBox}[6]{%
\setlength\fboxsep{1pt}
\setlength\fboxrule{0.0pt}
\fbox{\includegraphics[height=#2, clip, trim=#3 #4 #5 #6]{#1}}%
}

\newsavebox\figBox

\newcommand{\trimw}[6]{%
\sbox\figBox{\includegraphics{#1}}
\setlength{\tfwidth}{\the\wd\figBox}
\setlength{\tfheight}{\the\ht\figBox}
\setlength{\tfxa}{\tfwidth*\real{#3}}%
\setlength{\tfxb}{\tfwidth*\real{#4}}%
\setlength{\tfya}{\tfheight*\real{#5}}%
\setlength{\tfyb}{\tfheight*\real{#6}}%
\trimFigNoBox{#1}{#2}{\tfxa}{\tfya}{\tfxb}{\tfyb}%
}

\newcommand{\trimwb}[6]{%

\sbox\figBox{\includegraphics{#1}}
\setlength{\tfwidth}{\the\wd\figBox}
\setlength{\tfheight}{\the\ht\figBox}
\setlength{\tfxa}{\tfwidth*\real{#3}}%
\setlength{\tfxb}{\tfwidth*\real{#4}}%
\setlength{\tfya}{\tfheight*\real{#5}}%
\setlength{\tfyb}{\tfheight*\real{#6}}%
\trimFigWithBox{#1}{#2}{\tfxa}{\tfya}{\tfxb}{\tfyb}%
}

\newcommand{\trimh}[6]{%
\sbox\figBox{\includegraphics{#1}}
\setlength{\tfwidth}{\the\wd\figBox}
\setlength{\tfheight}{\the\ht\figBox}
\setlength{\tfxa}{\tfwidth*\real{#3}}%
\setlength{\tfxb}{\tfwidth*\real{#4}}%
\setlength{\tfya}{\tfheight*\real{#5}}%
\setlength{\tfyb}{\tfheight*\real{#6}}%
\trimFigHeightNoBox{#1}{#2}{\tfxa}{\tfya}{\tfxb}{\tfyb}%
}

\newcommand{\trimhb}[6]{%

\sbox\figBox{\includegraphics{#1}}
\setlength{\tfwidth}{\the\wd\figBox}
\setlength{\tfheight}{\the\ht\figBox}
\setlength{\tfxa}{\tfwidth*\real{#3}}%
\setlength{\tfxb}{\tfwidth*\real{#4}}%
\setlength{\tfya}{\tfheight*\real{#5}}%
\setlength{\tfyb}{\tfheight*\real{#6}}%
\trimFigHeightWithBox{#1}{#2}{\tfxa}{\tfya}{\tfxb}{\tfyb}%
}
\usepackage{algorithm} 
\usepackage{algpseudocode} 
\usepackage{listings} 

\begin{document}

\begin{frontmatter}

\title{A stable added-mass partitioned (AMP) algorithm \\
       for elastic solids and incompressible flow}

\author[rpi]{D.~A.~Serino\fnref{NSFgrfp,NSFgrants}}
\ead{serind@rpi.edu}

\author[rpi]{J.~W.~Banks\fnref{DOEThanks,PECASEThanks}}
\ead{banksj3@rpi.edu}

\author[rpi]{W.~D.~Henshaw\fnref{DOEThanks,NSFgrants}}
\ead{henshw@rpi.edu}

\author[rpi]{D.~W.~Schwendeman\corref{cor1}\fnref{DOEThanks,NSFgrants}}
\ead{schwed@rpi.edu}

\address[rpi]{Department of Mathematical Sciences, Rensselaer Polytechnic Institute, Troy, NY 12180, USA}

\cortext[cor1]{Department of Mathematical Sciences, Rensselaer Polytechnic Institute, 110 8th Street, Troy, NY 12180, USA.}

\fntext[DOEThanks]{This work was supported by contracts from the U.S. Department of Energy ASCR Applied Math Program.}

\fntext[PECASEThanks]{Research supported by a U.S. Presidential Early Career Award for Scientists and Engineers.}

\fntext[NSFgrants]{Research supported by the National Science Foundation under grants DMS-1519934 and DMS-1818926.}

\fntext[NSFgrfp]{Research supported by the National Science Foundation Graduate Research Fellowship under Grant No. DGE-1744655.}

\begin{abstract}
A stable added-mass partitioned (AMP) algorithm is developed for fluid-structure
interaction (FSI) problems involving viscous incompressible flow and compressible elastic solids.
%
The AMP scheme is stable and second-order accurate even when added-mass, and added-damping, effects are large.
%
Deforming composite grids are used to effectively handle the evolving geometry and large deformations.
The fluid is updated with an implicit-explicit (IMEX) fractional-step scheme whereby the velocity
is advanced in one step, treating the viscous terms implicitly, and the pressure is computed
in a second step. The AMP interface conditions for the fluid arise from the outgoing characteristic variables
in the solid and are partitioned into a Robin (mixed) interface condition
for the pressure, and
interface conditions for the velocity. The latter conditions include an impedance-weighted average
between fluid and solid velocities using a fluid impedance of a special form.
A similar impedance-weighted average is used to define interface values for the solid.
The new algorithm is verified for accuracy and stability on a number of useful benchmark problems
including a \textit{radial-piston} problem where exact solutions for radial and azimuthal motions
are found and tested. Traveling wave exact solutions are also derived and numerically verified for a solid disk surrounded by
an annulus of fluid. Fluid flow in a channel past a deformable solid annulus is computed and errors are estimated from
a self-convergence grid refinement study.
The AMP scheme is found to be stable and second-order accurate even for very difficult cases
of very light solids.

\end{abstract}


\end{frontmatter}

\clearpage
\tableofcontents

\clearpage
\section{Introduction} \label{sec:intro}

Fluid-structure interaction (FSI) problems involving incompressible fluids and elastic solids occur in many areas of engineering and applied science.  Examples include modeling flow-induced vibrations of structures (i.e.,~aircraft wings, undersea cables, wind turbines, and bridges) and simulating blood flow in arteries and veins.  These problems are typically modeled using partial differential equations for both the fluid and solid in their respective domains together with coupling conditions involving velocity and stress at the fluid-solid interface.  FSI problems are challenging mathematically due, in part, to the evolving fluid and solid domains, and the coupling at the moving interface.  Numerical algorithms have been developed to compute solutions of the governing equations, and these schemes can be classified broadly as monolithic or partitioned.  Monolithic schemes integrate the equations implicitly as a single large system of coupled equations, while partitioned schemes use separate solvers for the fluid and solid together with some approach to couple the solutions at each time step.  Partitioned schemes are often preferred for modularity and computational performance, but they can suffer from instabilities associated with the method used to couple the solutions at the interface.  For example, a standard traditional partitioned (TP) scheme employs a coupling approach in which the solid provides a Dirichlet (no-slip) boundary condition for the fluid, and then the fluid supplies a Neumann (traction) boundary condition for the solid.  This approach is appealing from a physical point of view in the case when the density of the solid is much greater than that of the fluid, but it is well known that TP schemes suffer from {\em added-mass} instabilities in the case when the solid is {\em light} relative to the fluid.  A method of iteration of the coupling conditions can be included in a TP algorithm in effort to suppress added-mass instabilities, but the effectiveness of this sub-time-step iteration is limited and it increases the cost of the algorithm.

In recent work~\cite{fib2014,fis2014}, we developed a new class of Added-Mass Partitioned (AMP)
algorithms for FSI problems coupling incompressible flow and elastic solids.  The work
in~\cite{fib2014} for bulk solids and in~\cite{fis2014} for thin structures both considered Stokes
flow and linear elastic solids, with the fluid-solid interface linearized about a flat surface.  The
algorithms used a fractional-step approach for the fluid in which the velocity is advanced in one
stage followed by the solution of a Poisson problem for the pressure.  A key ingredient in the
algorithms is the AMP interface conditions, which includes a mixed Robin condition for the pressure.
For the case of a bulk solid, this condition includes the {\em acceleration} of the solid through a
consideration of the outgoing characteristic variables of the hyperbolic system for the elastic solid and the interface
conditions.  The velocity and stress at the interface are set by impedence-weighted averages so that
the AMP conditions behave correctly in the limits of heavy and light solids.  Thus, the resulting
AMP algorithm is stable for {\em any} ratio of the density of the solid to that of the fluid,
without sub-time-step iterations, effectively suppressing added-mass instabilities.
In~\cite{beamins2016}, the approach for thin structures was extended to incompressible flows with
nonlinear convection and Euler-Bernoulli beams with finite deflection, and it was found that
the AMP scheme remains stable for both heavy and light solids.

The aim of the present work is to develop a stable and second-order accurate AMP algorithm for FSI
problems coupling incompressible flow and linear elastic bulk solids.  The scheme follows the work
in~\cite{fib2014}, but with significant modifications to handle finite deformation of the
fluid-solid interface, nonlinear advection in the fluid, an IMEX-type scheme for the fluid, and viscous added-damping effects.
The Navier-Stokes equations for the
fluid are solved in either an Eulerian frame or an Arbitrary-Eulerian-Lagrangian (ALE) frame which is coupled to the motion of the interface, while the
equations governing the solid are solved in a static Lagrangian reference frame.  Both the fluid and
solid domains are covered by composite overlapping grids~\cite{fsi2012}, with the deformation of the
fluid grid performed using the approach described originally in~\cite{mog2006}.  For example,
Figure~\ref{fig:rpiFig} shows a composite grid for an FSI problem involving a fluid in a channel
flowing past three elastic solids in the shape of the letters ``R'', ``P'' and ``I''.  The viscous
fluid flows from left to right, and it interacts with the three solids whose computed displacements
determine the deformed shapes.  The image on the far-right of the figure shows an enlarged view of
the composite grid near the fluid-solid interface of the right-most deformed solid.

{
\newcommand{\trimfigRPI}[2]{\trimhb{#1}{#2}{.2}{.25}{.375}{.375}}
\newcommand{\trimfigRPIg}[2]{\trimhb{#1}{#2}{.15}{.5}{.35}{.35}}
\newcommand{\trimfigRPIz}[2]{\trimhb{#1}{#2}{.59}{.28}{.4}{.44}}

\begin{figure}[h]
  \centering
  \resizebox{14cm}{!}{
    \begin{tikzpicture}
      \useasboundingbox(0,0.25) rectangle (16,8);
      \def\cxl{-.1}
      \def\cxr{11.5}
      \def\cyu{0}
      \def\cyl{-5}
      \def\wx{3.75}

      \draw(0,4) node[anchor=north west,yshift=4cm] {\trimfigRPI{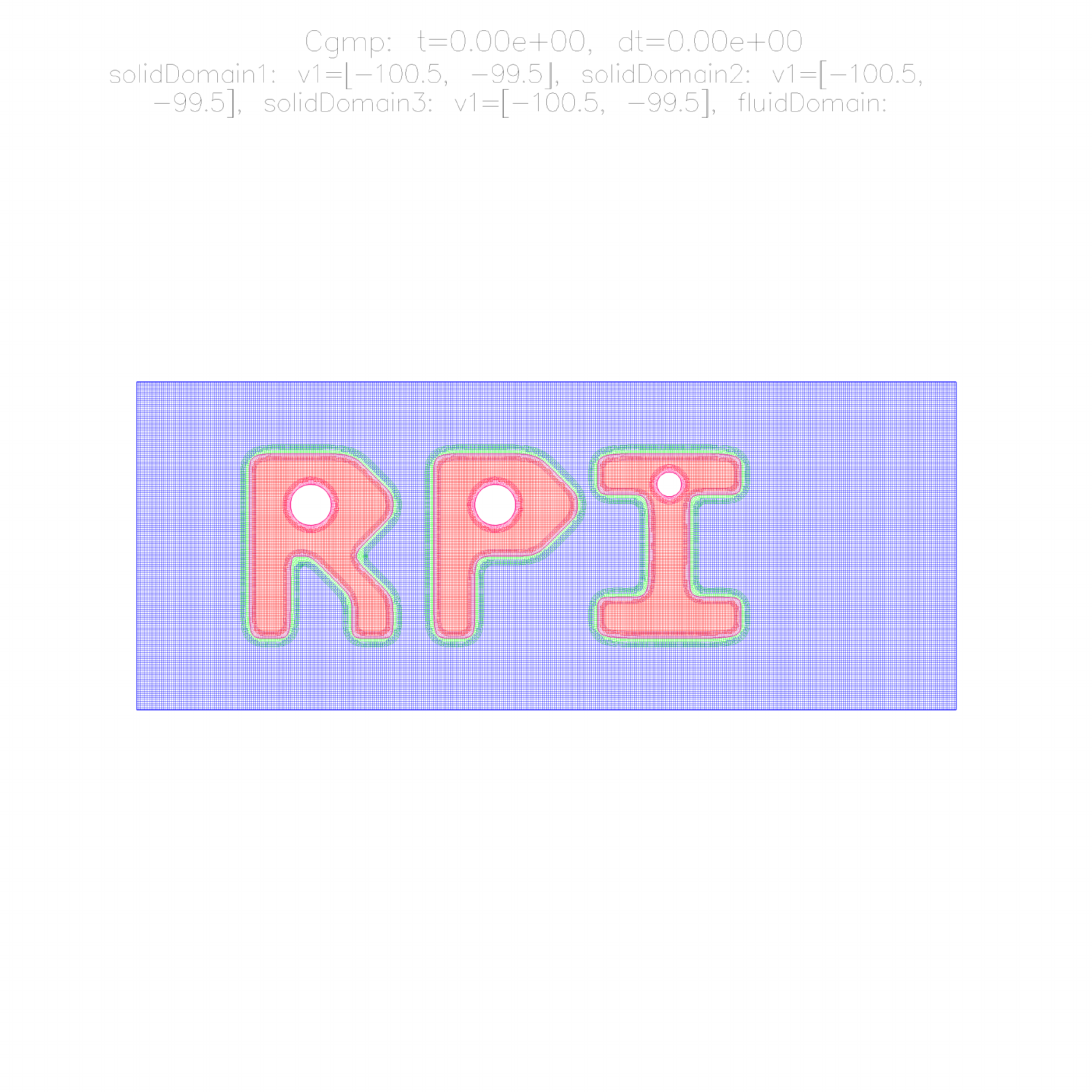}{\wx cm}};
      \draw(0,0) node[anchor=north west,yshift=4cm] {\trimfigRPI{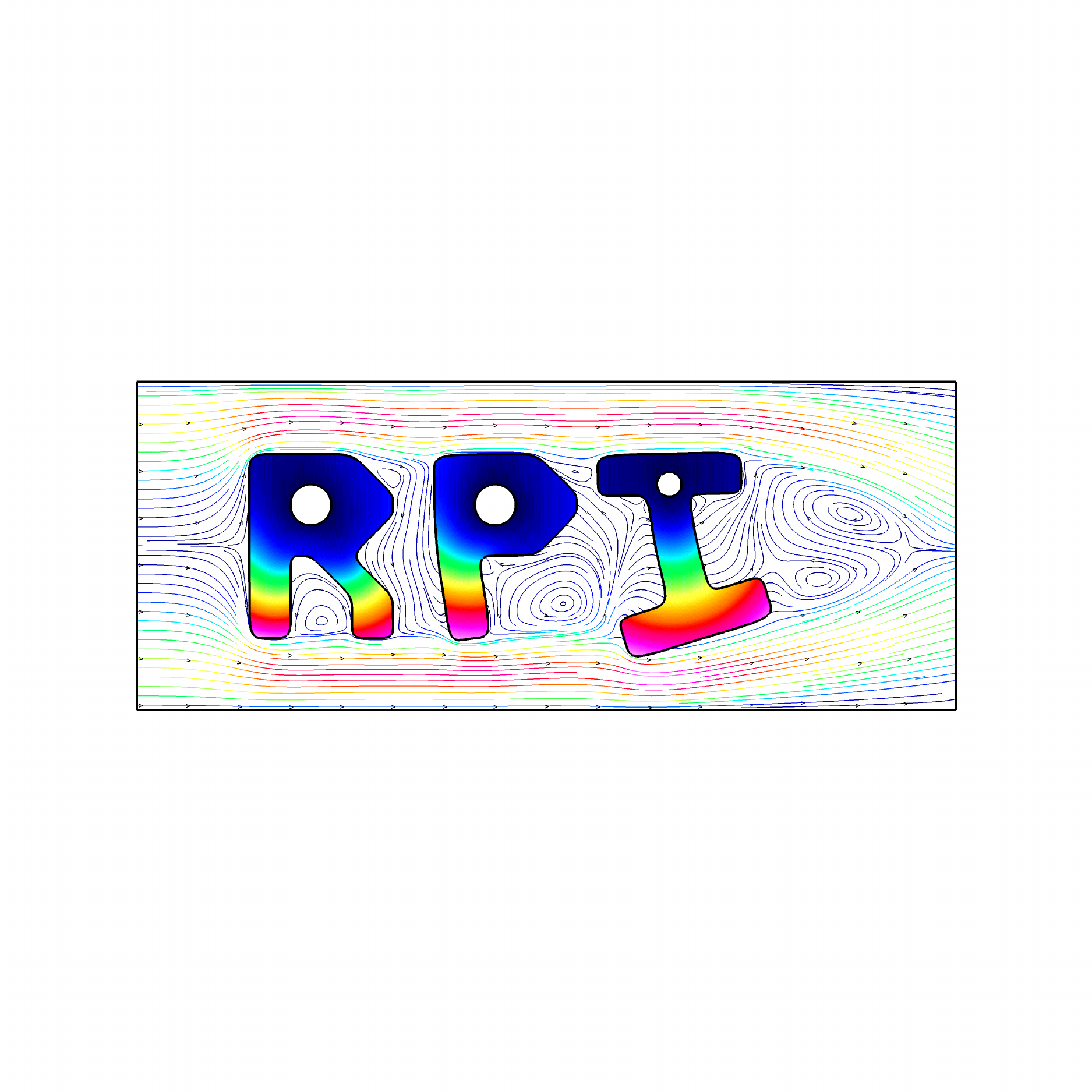}{\wx cm}};

      \begin{scope}[xshift=8.7cm,yshift=4cm]
        \draw(0,0) node[anchor=north west,yshift=4cm] {\trimfigRPIz{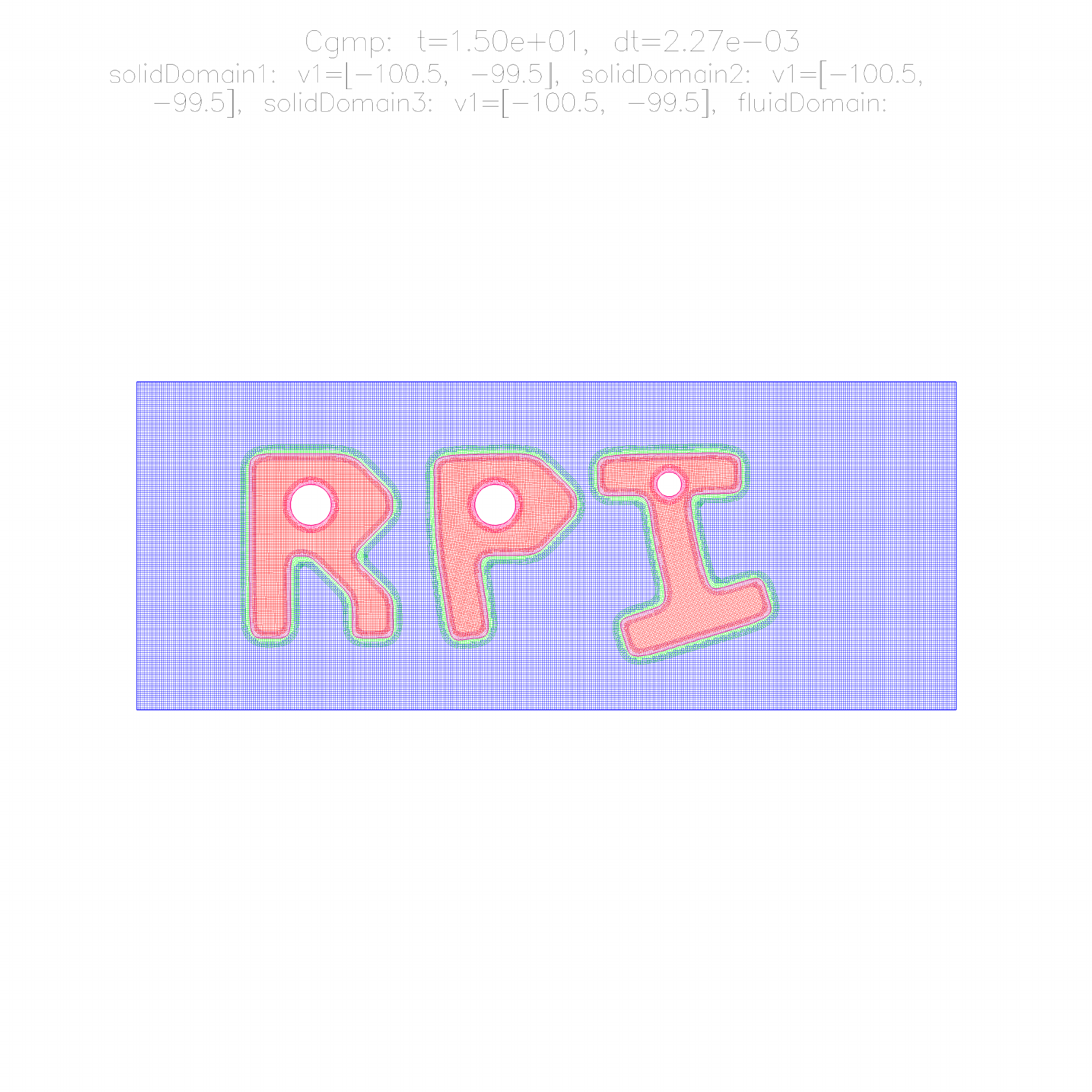}{7.75cm}};
        \draw(2,-2) node[draw,fill=white,anchor=west,xshift=2pt,yshift=1pt] {\footnotesize solid domain};
        \draw(4,2) node[draw,fill=white,anchor=west,xshift=2pt,yshift=1pt] {\footnotesize fluid domain};
        \draw(4,.5) node[draw,fill=white,anchor=west,xshift=2pt,yshift=1pt] {\footnotesize interface};
        \draw[->,thick,xshift=2pt] (4,.5) -- (3.1,0); 
      \end{scope}
    \end{tikzpicture}
  } 
  \caption{Top left: Initial grid is composed of a (blue) background grid and (green)
    deforming interface-fitted grids representing the fluid domain. The `RPI' shaped solid domain is 
    composed of (red) background grids, (purple) interface-fitted grids,
    and (pink) annular grids around interior cores where the displacement is set to zero.
    Bottom left: Solution at $t=15$ described by 
    streamlines in the fluid domain and shaded contours of the magnitude 
    of the solid displacement. 
    Right: Enlarged view of the grid after solid deformation.
    \label{fig:rpiFig}}
\end{figure}
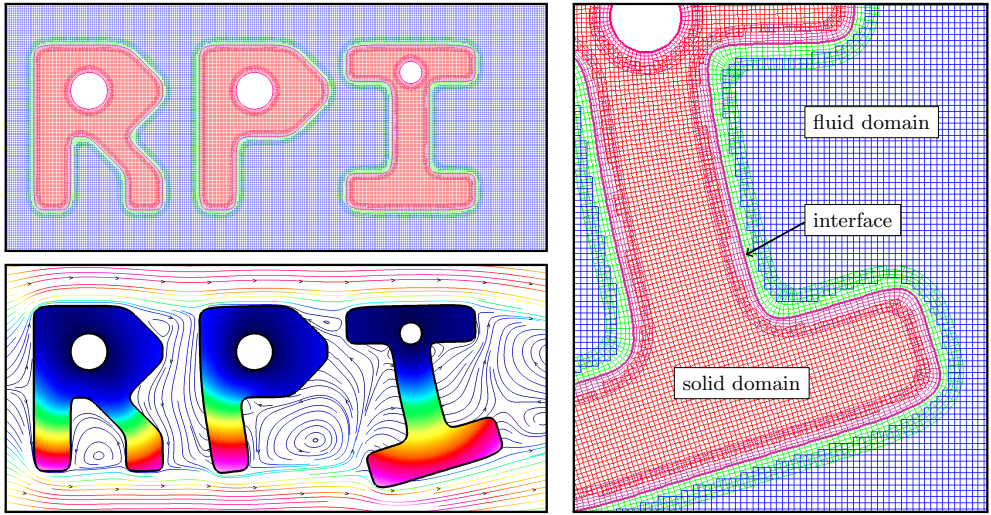

}

The equations of linear elasticity in the solid domain are advanced explicitly using the
time-stepping scheme described in~\cite{smog2012} for the hyperbolic system of equations in
first-order form.  The Navier-Stokes equations in the fluid domain are solved using a
fractional-step scheme.  Unlike the earlier work~\cite{fib2014}, the present solver uses an
implicit-explicit (IMEX) scheme based on~\cite{splitStep2003}.  In this scheme, the fluid velocity
is advanced first with the viscous terms treated implicitly while the pressure gradient and
nonlinear convection terms handled explicitly.  The fluid pressure is updated in a subsequent step
by solving a Poisson problem.  The incorporation of the IMEX scheme has necesitated important
changes to the original AMP interface conditions introduced in~\cite{fib2014}
to handle viscous added-damping effects that arise when the viscous CFL number becomes large.
Modifications to the AMP interface conditions are guided by the analysis of carefully chosen FSI model problems, as discussed in the companion manuscript~\cite{fibrmparXiv}.  
In particular, the analysis reveals a stable definition of the fluid impedance, which is a critical ingredient 
in the present AMP algorithm. The focus of the present work is a description of this algorithm for 
the full FSI problem, and its implementation in two space dimensions using deforming composite grids (DCG).
Performance of the algorithm is examined for a set of FSI
benchmark problems, including ones in which new exact solutions are derived. These results are used
to demonstrate the stability and accuracy of the algorithm, and a comparison to results from a traditional
partitioned scheme is also provided.


The overarching technique used to derive Added-Mass Partitioned schemes relies on the use of so-called compatibility boundary conditions, which are derived from the governing PDEs in concert with the interface conditions. This approach has also been used in a variety of other FSI regimes to yield stable and accurate partitioned solvers without the need for sub-iteration. For example in the case of FSI problems involving viscous incompressible flow, AMP schemes have been developed for thin elastic structures~\cite{beamins2016} (as noted above) and for rigid solids~\cite{rbinsmp2017,rbins2017,rbins3d2018}.  The first AMP schemes were developed for compressible flow~\cite{banks11a,sjogreen12}, and were subsequently extended into a general 2D framework using the DCG approach for FSI problems involving inviscid compressible flows and linear elastic solids~\cite{fsi2012}.  For the case of inviscid compressible flow, added-mass effects are localized due to finite propagation speeds of disturbances in the fluid~\cite{vanBrummelen2009}, and this leads to a somewhat simpler treatment of the coupling conditions in the corresponding AMP scheme relative to the incompressible case considered here.  Extensions of the scheme in~\cite{fsi2012} to FSI problems coupling inviscid compressible flow and rigid solids is described in~\cite{lrb2013}, with the case of nonlinear elastic solids detailed in~\cite{flunsi2016}.
In addition to the work cited above, there are other approaches to addressing added-mass-type instabilities in both partitioned and monolithic FSI solvers. For partitioned schemes, a typical strategy is to used under-relaxed subiterations, and previous research has shown that the number of necessary sub-iterations can be reduced by the use of Aitken acceleration or quasi-Newton methods, see \cite{Küttler2008,MEHL2016869}.  Additionally, sub-iteration schemes based on Robin-Neumann or Robin-Robin coupling, as opposed to the traditional Dirichlet-Neumann coupling, have also been shown to yield performance gains~\cite{Wang2018,BASTING2017312,BadiaNobileVergara2008,MokWallRamm2001,FernandezMullaertVidrascu2014,FernandezMullaertVidrascu2013,FernandezLandajuela2014,Gerardo_GiordaNobileVergara2010,BadiaNobileVergara2009,NobilePozzoliVergara2014}.
In related research, added-mass effects can also be treated by adding fictitious mass terms in the structure equations~\cite{BaekKarniadakis2012,YuBaekKarniadakis2013} or artificial compressibility in the fluid equations~\cite{RabackRuokolainenLyly2001,DegrooteSwillensBruggemanHaeltermanSegersVierendeels2010,Degroote2011}.  For difficult problems with large added-mass effects, 
monolithic schemes have been used to eliminate the need for sub-iterations~\cite{Forti2017,FigueroaVignonClementelJansenHughesTaylor2006}, and the performance of monolithic schemes has improved significantly through the application of multigrid~\cite{AULISA2018213,Heil2004,GeeKuttlerWall2011,HronTurekMonolithic2006}.

The remaining sections of the paper are organized as follows.  The equations governing the FSI problem are described in Section~\ref{sec:governingEquations}, and the AMP algorithm is described in Section~\ref{sec:algorithm}.  In the latter section, we begin with a discussion of the AMP interface conditions at a continuous level, which is followed by a detailed discussion of the AMP scheme.  Section~\ref{sec:numericalApproach} provides a description of the spatial approximations using deforming composite grids (DCG).  Numerical results confirming the stability and accuracy of the scheme and discussed in Section~\ref{sec:numericalResults}.  Some of the results use new exact solutions of benchmark FSI problems, which are described in the appendices.  Conclusions are given in Section~\ref{sec:conclusions}.

\section{Governing equations and interface conditions} \label{sec:governingEquations}

We consider the coupled evolution of an incompressible fluid and a linear elastic solid. 
The fluid occupies the domain $\xv \in \OmegaF(t)$, where
$\xv=(x_1,x_2,x_3)$ is a vector of physical coordinates and $t$ is time. 
The equations for the solid are written in terms of the Lagrangian coordinate $\xsv=(\bar x_1,\bar x_2,\bar x_3)$ for
a reference configuration $\xsv \in \OmegaS$ at $t=0$.  (An overbar is used here and elsewhere to denote
quantities associated with the solid.) The fluid and solid
are coupled at an interface described by
$\xv\in\GammaIf(t)$ in physical space and $\xsv \in \GammaIs$
in the corresponding reference space. Figure~\ref{fig:rpiFig} illustrates these geometric features using a model fluid-structure interaction problem.

The fluid is described by the incompressible Navier-Stokes equations
\bse
\begin{alignat}{2}
\grad \cdot \vv &= 0, \qquad &&\xv \in \OmegaF(t),
\label{eq:fluidContinuity} \\
\rho \vv_t + \rho (\vv \cdot \grad) \vv &=
  \grad \cdot \sigmav, \qquad &&\xv \in \OmegaF(t),
\label{eq:fluidMomentum}
\end{alignat}
\ese
where $\vv(\xv,t)$ is the velocity, $\sigmav(\xv,t)$ is the stress tensor and $\rho$ is the (constant) density.
The fluid stress tensor is given by
\begin{align}
\sigmav &= - p \Iv + \tauv, \label{eq:fluidStressTensor} 
\end{align}
where $p(\xv,t)$ is the pressure, $\Iv$ is the identity matrix and $\tauv(\xv,t)$ is
the viscous stress tensor.  The latter quantity is given in terms of the gradient of the fluid velocity by
\begin{align}
\tauv &= \mu \left(\grad \vv + \left(\grad \vv\right)^T \right), \label{eq:fluidViscousStressTensor}
\end{align}
where $\mu$ is the dynamic viscosity of the fluid (assumed to be constant).

We consider the governing equations for the fluid in velocity-pressure form.  In this form, the continuity equation in~\eqref{eq:fluidContinuity} is replaced by a Poisson equation for the pressure given by
\begin{alignat}{2}
\Delta p &= -\rho \grad \vv : \left(\grad \vv\right)^T, \qquad &&\xv \in \OmegaF(t),
\label{eq:pressurePoissonEquation}
\end{alignat}
where
\begin{align}
\grad \vv : \left(\grad \vv\right)^T \equiv 
\sum_{i=1}^3 \sum_{j=1}^3 \frac{\partial v_i}{\partial x_j} \frac{\partial v_j}{\partial x_i}.
\end{align}
Here, components of the fluid velocity $\vv$ are denoted by $(v_1,v_2,v_3)$, and a similar notation is used for the components of tensors, e.g.~the components of $\sigmav$ are given by $\sigma_{mn}$ with $m,n=1$, 2, 3.  In the velocity-pressure form of the equations, boundary conditions are required for the the pressure-Poisson equation. A suitable choice that ensures consistency with solutions of the velocity-divergence form is $\grad \cdot \vv = 0$ for $\xv\in\partial\OmegaF(t)$.

The evolution of the solid is described by the equations of linear elasticity.  These equations can be written in the form
\bse
\label{eq:solid}
\begin{alignat}{2}
\usv_t &= \vsv, \qquad &&\xvs \in \OmegaS ,
\label{eq:solidDisplacement} \\
\rhos \vsv_t &= \grads \cdot \sigmasv, \qquad &&\xvs \in \OmegaS,
\label{eq:solidMomentum}
\end{alignat}
\ese
where $\usv(\xsv,t)$ is the displacement of the solid, $\vsv(\xsv,t)$ its velocity, and $\rhos$ its density (assumed constant).
The Cauchy stress tensor $\sigmasv(\xsv,t)$ in~\eqref{eq:solidMomentum} is defined by
\begin{align}
\sigmasv &= \lambdas (\grads \cdot \usv) \Iv 
+ \mus \left(\grads \usv + \left(\grads \usv\right)^T\right),
\label{eq:solidStressTensor}
\end{align}
where $\lambdas$ and $\mus$ are Lam\'e parameters (taken to be constants).
The position of the solid in physical space is determined by the mapping
\begin{align}
\xv = \xsv + \usv(\xsv,t),
\label{eq:solidPhysicalCoordinates}
\end{align}
for $\xvs \in \OmegaS$.  The corresponding deformation gradient given by
$\bar F=\grads \xv=\Iv+\grads \usv$ is assumed to be a small perturbation of the
identity for a linear elastic solid so that a one-to-one mapping from the
reference space to the physical space of the solid exists. 
%
Following the approach described in~\cite{smog2012}, the solid equations are treated as a first order system of PDEs in time and space. This formulation
is achieved by taking the time derivative of \eqref{eq:solidStressTensor} to obtain 
\begin{align}
\sigmasv_t &= \lambdas (\grads \cdot \vsv) \Iv 
+ \mus \left(\grads \vsv + \left(\grads \vsv \right)^T\right),\qquad \xvs \in \OmegaS,
\label{eq:solidStressTensorDeriv}
\end{align}
which is then appended to previous system of equations in~\eqref{eq:solid}.
In this form, upwind solvers can be used to advance displacement, velocity and stress of the solid.  We note that~\eqref{eq:solidStressTensor} is enforced at $t=0$, and then the stress-strain relation is satisfied approximately for $t>0$ through compatibility conditions applied at the boundary and a stress-strain relaxation term in the numerical scheme (see~\cite{flunsi2016}).

The interface between the fluid and solid, in physical space given by $\xv\in\GammaIf(t)$, is determined by the mapping
in~\eqref{eq:solidPhysicalCoordinates} for the boundary of the solid given by $\xsv \in \GammaIs$.  The interface is
assumed to be smooth so that a well-defined normal to the interface exists.  Along the interface, the following matching
conditions hold:
\bse
\label{eq:matching}
\begin{alignat}{2}
\vv &= \vsv, \qquad &&\xv \in \GammaIf(t),
\label{eq:velocityMatching} \\
\sigmav \nv &= 
\sigmasv \nv, \qquad &&\xv \in \GammaIf(t),
\label{eq:tractionMatching}
\end{alignat}
\ese
%
where $\nv(\xv,t)$ is the outward unit normal to the fluid domain, i.e.~$\nv$ points from the fluid domain to the solid domain.  Suitable boundary
conditions are applied on the boundaries of the fluid and solid domians not included in $\GammaIf(t)$, and initial
conditions on $\vv$, $\usv$ and $\vsv$ are set to close the problem.

Since the equations governing the solid form a hyperbolic system, it 
is informative to consider the propagation of characteristics. The 
primary and secondary wave speeds are defined as
\begin{align}
\cp = \sqrt{\frac{\lambdas+2 \mus}{\rhos}}, \qquad
\cs = \sqrt{\frac{ \mus}{\rhos}},
\label{eq:solidWaveSpeeds}
\end{align}
which correspond to compression and sheer waves, respectively. At the interface, the Riemann variables corresponding to incoming and outgoing characteristics can be 
obtained from the first-order equations in~\eqref{eq:solidMomentum} and~\eqref{eq:solidStressTensorDeriv} projected onto the normal to the interface.  These variables for the compression and shear waves are given by
\bse
\label{eq:characteristicFunctions}
\begin{alignat}{2}
\Ac(\sigmasv, \vsv) &= \nv^T \sigmasv \nv - \zp \nv^T \vsv ,
\quad \phantom{m = 1,2,}
\qquad &&\xv \in \GammaIf(t)
\quad \text{(incoming),}
\label{eq:ACharacteristic}
\\
\Ac_m(\sigmasv, \vsv) &= \tnv_m^T \sigmasv \nv - \zs \tnv_m^T \vsv, 
\quad m = 1,2,
\qquad &&\xv \in \GammaIf(t)
\quad \text{(incoming),}
\label{eq:AmCharacteristic}
\end{alignat}
corresponding to information propagating into the solid domain, and
\begin{alignat}{2}
\Bc(\sigmasv, \vsv) &= \nv^T \sigmasv \nv + \zp \nv^T \vsv ,
\quad \phantom{m = 1,2,}
\qquad &&\xv \in \GammaIf(t)
\quad \text{(outgoing),}
\label{eq:BCharacteristic}
\\
\Bc_m(\sigmasv, \vsv) &= \tnv_m^T \sigmasv \nv + \zs \tnv_m^T \vsv ,
\quad m = 1,2,
\qquad &&\xv \in \GammaIf(t)
\quad \text{(outgoing),}
\label{eq:BmCharacteristic}
\end{alignat}
\ese
corresponding to information propagating out of the solid domain.
Here, $\tnv_1(\xv,t)$ and $\tnv_2(\xv,t)$ are mutually orthogonal unit vectors tangent to $\GammaIf(t)$, and thus orthogonal to $\nv$, and
$\zp = \rhos \cp$ and $\zs = \rhos \cs$ are the solid impedances 
for compression and shear waves, respectively. 
The functions given in~\eqref{eq:characteristicFunctions}, corresponding to the Riemann variables for the solid, can be used to derive an equivalent set of characteristic-based interface conditions which take into 
account the direction of propagation.  These conditions are
\bse
\label{eq:matching-char}
\begin{align}
\Ac(\sigmasv, \vsv)    &= \Ac(\sigmav, \vv),\qquad \Ac_m(\sigmasv, \vsv)  = \Ac_m(\sigmav, \vv),
 \qquad \xv \in \GammaIf(t),
\label{eq:OutgoingMatching} \\
\Bc(\sigmav, \vv)    &= \Bc(\sigmasv, \vsv), \qquad\; \Bc_m(\sigmav, \vv) =  \Bc_m(\sigmasv, \vsv),
 \qquad \xv \in \GammaIf(t),
 \label{eq:IncomingMatching}
\end{align}
\ese
which are linear combinations of the matching conditions in~\eqref{eq:matching}.
%

This paper is concerned with partitioned algorithms for solving FSI problems involving incompressible fluids and
linear elastic solids.  Such algorithms employ PDE solvers for the equations governing the fluid and solid separately
with the matching conditions at the fluid-solid interface providing boundary conditions for each domain solver (along
with conditions on the boundaries not included in the shared interface).  For example, in traditional
partitioned algorithms, the fluid velocity is determined by the solid, while the solid traction is determined by
the fluid.  This ordering can often be used for the case of heavy solids, but it suffers from added-mass
instabilities when the solid density is light or the solid grid is very fine.  In the reverse scheme, referred to as
the anti-traditional scheme, the solid imposes a traction condition on the fluid, while the fluid velocity provides
a condition for the solid velocity.  This reverse ordering works well for light solids or fine grids, but 
may exhibit instability for heavy solids or coarse grids.

The added-mass partitioned (AMP) scheme described in the next section follows the work in~\cite{fib2014}, and uses
the characteristic-based matching conditions in~\eqref{eq:matching-char}.
At a high level, these are implemented as boundary conditions for the solid solver as
\bse
\label{eq:BCchar}
\begin{equation}
\Ac(\sigmasv, \vsv)  = \Ac(\sigmav^{(p)}, \vv^{(p)}),\qquad \Ac_m(\sigmasv, \vsv)  = \Ac_m(\sigmav^{(p)}, \vv^{(p)}),
 \qquad \xv \in \GammaIf(t),
\label{eq:ACharacteristicBC}
\end{equation}
and as boundary conditions for the fluid solver as
\begin{equation}
\Bc(\sigmav, \vv)   = \Bc(\sigmasv^{(p)}, \vsv^{(p)}), \qquad \Bc_m(\sigmav, \vv) =  \Bc_m(\sigmasv^{(p)}, \vsv^{(p)}),
 \qquad \xsv \in \GammaIf(t),
 \label{eq:BCharacteristicBC}
\end{equation}
\ese
where the $(p)$ superscript denotes suitable {\em predicted} values (as described in more detail below) so that the right-hand sides of the boundary conditions in~\eqref{eq:BCchar} are treated as known.  This characteristic-based decomposition of the interface matching conditions is a key step in the derivation of the AMP
algorithm, and it leads to a stable scheme for a wide range of fluid and solid densities.  In particular, in the limits of heavy solids and light solids, 
the AMP scheme approaches the traditional and anti-traditional partitioned schemes, respectively.



\section{FSI algorithms} \label{sec:algorithm}

In this section, we describe the AMP algorithm, and we also discuss the traditional partitioned~(TP) algorithm and the anti-traditional partitioned~(ATP) algorithm for comparative purposes.  All three partitioned schemes use the same fluid and solid domain solvers.  For the fluid domain, we integrate the governing equations in velocity-pressure form using a fractional-step scheme in which the fluid velocity is advanaced in one step and the pressure-Poisson equation is solved to update the fluid pressure in the other step (see~\cite{splitStep2003}).  An explicit upwind Godunov method is used to advance the velocity and stress of the solid in its reference domain following the approach described in~\cite{smog2012}, while a Lax-Wendroff scheme is used to advance the solid displacement.  These two domain solvers, both second-order accurate in space and time, are embedded in the FSI algorithms using a predictor-corrector approach.
The first step of the FSI time-stepping loop involves advancing the solid variables to the next time level.  The solid displacement predicts the new position of the fluid-solid interface and this information is used to deform the fluid grid and compute its {\em grid velocity}.  The fluid-domain solver is then called to advance the fluid velocity to the next time level and update the fluid pressure.  In the fluid solver, the diffusion terms of the momentum equations are treated implicitly, while the convection and pressure-gradient terms are approximated explicitly.  Following this fluid-domain step, the solid interface data is updated using fluid data at the new time level.  The steps outlined above are essentially repeated in the corrector step.  Depending on the densities of the fluid and solid (and on the grid spacing), the TP and ATP schemes may require addition sub-iterations for stability, in which case the corrector step is repeated.

The AMP algorithm differs from that of the TP and ATP schemes in the implementation of the fluid-solid interface conditions.  To this end, we begin in Section~\ref{sec:interfaceConditions} by describing the various treatments of the interface conditions, including compatability conditions, at a continuous level.  This discussion is then followed in Section~\ref{sec:partitionedAlgorithms} by a detailed description of the various FSI algorithms at a discrete level.

\subsection{Continuous interface conditions} \label{sec:interfaceConditions}

Using the formulas for the fluid stress tensor in~\eqref{eq:fluidStressTensor} and~\eqref{eq:fluidViscousStressTensor}, the AMP interface conditions in~\eqref{eq:BCharacteristicBC} for the fluid become
\bse
\begin{align}
-p + \nv^T \tauv \nv + \zp \nv^T \vv = \nv^T \sigmasv \nv + \zp \nv^T \vsv, \quad \phantom{m = 1,2,}\qquad \xv \in \GammaIf(t),
\label{eq:AMPNormalCharacteristic} \\
\tnv_m^T \tauv \nv + \zs \tnv_m^T \vv = \tnv_m^T \sigmasv \nv + \zs \tnv_m^T \vsv,  \quad m = 1,2, \qquad \xv \in \GammaIf(t),\label{eq:AMPTangentialCharacteristic}
\end{align}
\ese
with the $(p)$ superscipt on the solid variables suppressed.  These conditions, along with $\grad \cdot \vv = 0$ applied on the boundary, are a sufficient set of conditions for the fluid equations in velocity-pressure form.  The left-hand side of the interface conditions in~\eqref{eq:AMPTangentialCharacteristic} only involve the fluid velocity (and its derivatives), and these conditions with the divergence-free constraint are used as boundary conditions for the momentum equation in~\eqref{eq:fluidMomentum} to advance the fluid velocity.  The stress and velocity of the solid appearing on the right-hand side of~\eqref{eq:AMPTangentialCharacteristic} are obtained from an explicit time step of the hyperbolic system of equations for the solid.  Once the fluid velocity is advanced to the next time step, a Poisson problem is solved to update the fluid pressure.  Following the analysis in~\cite{fib2014,fibrmparXiv}, the interface condition in~\eqref{eq:AMPNormalCharacteristic} is used with the momentum equation in~\eqref{eq:fluidMomentum} to derive a Robin condition for the pressure-Poisson problem that takes the form
\begin{align}
-p - {\zp \dt\over\rho} \partial_n p =  \nv^T(\sigmasv \nv-\tauv \nv) + \zp \dt\nv^T\bigl( \vsv_t + \nu\grad\times\grad\times\vv \bigr), \qquad \xv \in \GammaIf(t),
\label{eq:AMPpressureBC}
\end{align}
where $\partial_n=\nv\cdot\grad$ is the normal derivative and $\nu=\mu/\rho$ is the kinematic viscosity of the fluid.  Following~\cite{splitStep2003}, we have used the identity, $\Delta \vv = - \grad \curl \grad \curl \vv$, noting that $\grad \cdot \vv = 0$, to replace $\Delta \vv$ on the right-hand side of~\eqref{eq:AMPpressureBC} in favor of the curl-curl operator.  This is done for improved stability of the fractional-step scheme.  The condition in~\eqref{eq:AMPpressureBC}, along with suitable conditions for $\xv\in\partial\OmegaF(t)\backslash\GammaIf(t)$, is used for the Poisson equation in~\eqref{eq:pressurePoissonEquation} for the pressure.


The AMP algorithm requires an additional \textit{interface projection} to ensure the fluid and solid
interface values match at the end of the time-step. This projection defines common interface values
for the velocity and traction using an impedance-weighted average.  In particular, the normal and
tangential components of the velocity on the interface
$\xv \in \GammaIf(t)$ are given by 
\bse
\label{eq:AMPinterfaceVelocity}
\begin{alignat}{2}
\nv^T\vv^I=\;&{\zf\over\zf+\zp}\nv^T\vv+{\zp\over\zf+\zp}\nv^T\vsv+{1\over\zf+\zp}\nv^T\bigl(\sigmasv\nv-\sigmav\nv\bigr), \quad &&\phantom{m = 1,2,}
\label{eq:AMPinterfaceNormalVelocity} \\
\tnv_m^T\vv^I=\;&{\zf\over\zf+\zs}\tnv_m^T\vv+{\zs\over\zf+\zs}\tnv_m^T\vsv+{1\over\zf+\zs}\tnv_m^T\bigl(\sigmasv\nv-\sigmav\nv\bigr) ,  \quad &&m = 1,2,
\label{eq:AMPinterfaceTangentialVelocity}
\end{alignat}
\ese
while the traction is defined from an inverse impedance average, 
\bse
\label{eq:AMPinterfaceTraction}
\begin{alignat}{2}
\nv^T\sigmav^I\nv=\;&{\zf^{-1}\over\zf^{-1}+\zp^{-1}}\nv^T\sigmav\nv+{\zp^{-1}\over\zf^{-1}+\zp^{-1}}\nv^T\sigmasv\nv
+{1\over\zf^{-1}+\zp^{-1}}\nv^T\bigl(\vsv-\vv\bigr), \quad &&\phantom{m = 1,2,}
\label{eq:AMPinterfaceNormalTraction}\\
\tnv_m^T\sigmav^I \nv=\;&{\zf^{-1}\over\zf^{-1}+\zs^{-1}}\tnv_m^T\sigmav\nv+{\zs^{-1}\over\zf^{-1}+\zs^{-1}}\tnv_m^T\sigmasv\nv
+{1\over\zf^{-1}+\zs^{-1}}\tnv_m^T\bigl(\vsv-\vv\bigr) ,  \quad &&m = 1,2.
\label{eq:AMPinterfaceTangentialTraction}
\end{alignat}
\ese
Since the fluid is incompressible there is apparently no natural definition for the fluid impedance~$z_f$.
In our previous paper~\cite{fib2014}, we used $z_f=\rho H/\Delta t$, where $H$
measured the height of the fluid domain, and the corresponding scheme was found to be very insensitive to the choice of $H$. 
This previous scaling was sufficient when using explicit time-stepping
in the fluid. In the current paper, however, we use an IMEX scheme in the fluid where viscous effects
are more important since the viscous CFL parameter $\nu\dt/h^2$ may become very large.
In this regime, the previous choice of $\zf$ does not result in a stable scheme.
The analysis of a model problem that accounts for viscous effects performed in the companion paper~\cite{fibrmparXiv},
provides the suitable form for $z_f$ which is a {\em key result} that leads to a stable scheme.
\newcommand{\Cam}{\Cc_{\rm AM}}
\newcommand{\Cad}{\Cc_{\rm AD}}
\begin{AMPFluidImpedence}
The fluid impedance $\zf$ accounting for both added-mass and viscous added-damping effects is chosen as
\ba
  z_f &\eqdef  \Cam \Bigl(\frac{\rho h}{\dt}\Bigr) + \Cad \Bigl(\frac{\mu}{h}\Bigr)  \label{eq:fluidImpedance} ,
\ea
where $\rho$ is the fluid density, $h$ measures the fluid grid spacing in the normal direction to the interface, $\dt$ is the time-step, $\mu$ is the dynamic viscosity of the fluid, and $\Cam=1$ and $\Cad=2$ are constants following the analysis in~\cite{fibrmparXiv}.
\label{eq:ZfDef}
\end{AMPFluidImpedence}
In agreement with~\cite{fibrmparXiv}, the choice~\eqref{eq:fluidImpedance} is found to yield a stable AMP algorithm for all FSI problems considered, including those covering a wide range of solid densities from very heavy to light. On the other hand, stable AMP algorithms also result when scaling the constants $\Cad$ and $\Cam$ by a factor of $1/2$ or $2$. This indicates that there is some robustness in the choice of the constants in~\eqref{eq:fluidImpedance}. However, it appears to be important to keep both terms in~\eqref{eq:fluidImpedance}, since failing to do so can lead to subtle instabilities at late times.

We note that with appropriate scales for length and time from the analysis given by~$h$ and~$\dt$, respectively, there is a simple scaling argument leading to the form for $\zf$ in~\eqref{eq:fluidImpedance}.  Consider a dimensional analysis of the fluid momentum equation in one dimension,
\ba
   \rho v_t + p_x = \mu v_{xx}.
\label{eq:simpleModel}
\ea
In addition, let $V$ and $P$ be reference scales for velocity and pressure (stress), respectively.  In the inertia-dominated limit, the two terms on the left-hand side of~\eqref{eq:simpleModel} set the scale for pressure leading to the balance
\[
P=\left({\rho h\over\dt}\right)V.
\]
The implied impedance, $(z_f)_{{\rm AM}}=\rho h/\dt$, is expected to be dominant for added-mass effects.  In the viscous-dominated limit when added-damping effects are important, the pressure gradient term in~\eqref{eq:simpleModel} balances the viscous stress term leading to
\[
P=\left({\mu\over h}\right)V.
\]
This balance suggests a second expression for the fluid impedance given by $(z_f)_{{\rm AD}}=\mu/h$.  A linear combination of the limiting forms agrees with the choice for $\zf$ in~\eqref{eq:fluidImpedance}.
The impedance-weighted averages~\eqref{eq:AMPinterfaceVelocity} and~\eqref{eq:AMPinterfaceTraction}
 are defined to ensure that the interface conditions behave correctly in the limits of a heavy and light solids.
In the limit of a heavy solid, i.e.~$\rhos\rightarrow\infty$ ($\zp$, $\zs\rightarrow\infty$), the interface conditions for the fluid should approach those of a (moving) no-slip wall.  For $\zs\rightarrow\infty$, \eqref{eq:AMPTangentialCharacteristic} reduces to
\begin{align}
\tnv_m^T\vv=\tnv_m^T\vsv,\quad m=1,2,\qquad \xv \in \GammaIf(t),
\label{eq:tangentialVelocityHeavySolid}
\end{align}
so that the tangential components of the fluid velocity match those of the solid wall.  We note that the tangential components of the interface velocity defined in~\eqref{eq:AMPinterfaceTangentialVelocity} agree with those given in~\eqref{eq:tangentialVelocityHeavySolid} in this limit.  The normal component of the fluid is given by the normal component of the interface velocity, defined in~\eqref{eq:AMPinterfaceNormalVelocity}, which becomes
\begin{align}
\nv^T\vv=\nv^T\vv^I=\nv^T\vsv,\qquad \xv \in \GammaIf(t),
\label{eq:normalVelocityHeavySolid}
\end{align}
so that the normal component of the velocity of the fluid also matches that of the solid wall.  In this heavy-solid limit, the divergence-free constraint applied on the boundary becomes a compatibility condition used to specify values of the normal component of the fluid velocity at ghost points.  The Robin condition for the pressure given in~\eqref{eq:AMPpressureBC} in the limit $\zp\rightarrow\infty$ reduces to
\begin{align}
- {1\over\rho} \partial_n p =  \nv^T\vsv_t + \nu\nv^T\bigl(\grad\times\grad\times\vv\bigr), \qquad \xv \in \GammaIf(t).
\label{eq:pressureBCHeavySolid}
\end{align}
This Neumann condition for the pressure is often used for the pressure-Poisson problem at a no-slip wall.

In the limit of a light solid, i.e.~$\rhos\rightarrow0$ ($\zp$, $\zs\rightarrow0$), the interface conditions for the fluid should approach those of a free surface.  With $\zs\rightarrow0$, \eqref{eq:AMPTangentialCharacteristic} becomes
\begin{align}
\tnv_m^T\sigmav\nv=\tnv_m^T\sigmasv\nv,\quad m=1,2,\qquad \xv \in \GammaIf(t),
\label{eq:tangentialStressLightSolid}
\end{align}
and thus the tangential components of the fluid stress match those of the light solid.  In this limit, the conditions in~\eqref{eq:tangentialStressLightSolid} and the divergence-free constraint on the fluid velocity provide the boundary conditions needed to advance the velocity using~\eqref{eq:fluidMomentum}, along with suitable boundary conditions for $\xv\in\partial\OmegaF(t)\backslash\GammaIf(t)$ which we do not discuss.  The interface condition in~\eqref{eq:AMPNormalCharacteristic} or the derived Robin condition in~\eqref{eq:AMPpressureBC} both reduce to
\begin{align}
-p = - \nv^T \tauv \nv + \nv^T \sigmasv \nv, \qquad \xv \in \GammaIf(t),
\label{eq:normalStressLightSolid}
\end{align}
in the limit $\zp\rightarrow0$, and this Dirichlet condition on pressure is used for the pressure-Poisson problem.  Essentially, the solid traction sets the fluid traction along the interface in the light-solid limit, and the definitions of the interface traction in~\eqref{eq:AMPinterfaceTraction} only confirms this assignment.  Finally, the impedance-weighted averages in~\eqref{eq:AMPinterfaceVelocity} reduce to
\begin{align}
\vv^I=\vv, \qquad \xv \in \GammaIf(t),
\label{eq:fluidVelocityLightSolid}
\end{align}
so that the velocity of the interface is specified by the fluid velocity, as expected, which then prescribes the velocity of the solid as a boundary condition.

\subsection{Discrete partitioned algorithms} \label{sec:partitionedAlgorithms}

In this section, the details of the AMP algorithm are presented, along with the special cases of the TP and ATP schemes for comparison.  The description of the AMP algorithm here extends the framework of the original scheme given in~\cite{fib2014} to moving grids and deforming interfaces.  A semi-implicit predictor-corrector scheme is used to solve the fluid equations following the work in~\cite{splitStep2003} and extending the explicit scheme in~\cite{fib2014}, while a second-order accurate upwind scheme is used to advance the solid equations.  The evolving fluid domain is discretized on a moving overlapping grid, which is discussed in more detail in Section~\ref{sec:numericalApproach}.  The solid domain is discretized using a static overlapping grid in its reference configuration.

Algorithm~\ref{alg:amp} outlines the steps for the AMP~scheme.  Let $\xv_\iv^n$ for $\iv \in \OmegaFh$ denote the time-dependent grid points belonging to the fluid domain (including grid points on the boundaries and the interface), and let $\dot{\xv}_\iv^n$ denote the velocity of the grid points at the $n\sp{{\rm th}}$ time step, $t^n$.  The indices of grid points on the interface of the fluid grid are denoted by $\iv \in \GammaIfh$.  Similarly, let $\xsv_\isv$ for $\isv \in \OmegaSh$ denote the grid points belonging to the static overlapping grid for the solid (including boundaries and the interface), while $\isv \in \GammaIsh$ are the indicies associated with grid points on the interface of the solid grid.  Discrete approximations for the fluid variables are given by $\vv_\iv\sp{n}\approx \vv(\xv_\iv\sp{n},t\sp{n})$ and $p_\iv\sp{n}\approx p(\xv_\iv\sp{n},t\sp{n})$, while the solid variables are given by $\usv_\isv\sp{n}\approx \usv(\xsv_\isv,t\sp{n})$, $\vsv_\isv\sp{n}\approx \vsv(\xsv_\isv,t\sp{n})$ and $\sigmasv_\isv\sp{n}\approx \sigmasv(\xsv_\isv,t\sp{n})$.

{
    \def\alignspace{\hspace{1.2em}}
\begin{algorithm}\caption{\rm Added-mass partitioned (AMP) scheme \label{alg:ampAlgorithm}}
\small
\[
\begin{array}{l}
\hbox{// \textsl{Predictor steps}}\smallskip\\
1.\text{ Predict solid:}\\ 
\alignspace 
\begin{cases}
\Usv_\isv^{(p)} = \Usv_\isv^n
+ \dt \Vsv_\isv^n
+ \frac{\dt^2}{2\rhos} \Grads \cdot \Sigmasv_\isv^n, 
& \qquad \isv \in \OmegaSh. \\
\Qsv_\isv^{(p)} = \Qsv_\isv^n 
- \dt \sum_{m=1}^3 \frac{1}{\dxs_m} 
\Bigl( \Fsv_{m, \; \isv}^+ - \Fsv_{m, \; \isv}^- \Bigr), & \qquad \isv \in \OmegaSh,
\end{cases} 
\medskip\\
2.\text{ Predict fluid grid: advance fluid grid to $t^{n+1}$ using $\Usv_\isv^{(p)}$ for $\isv \in \GammaIsh$, and compute grid velocity.}
\medskip\\
3.\text{ Predict fluid velocity:}\\ 
\alignspace 
\begin{cases}
\vv_\iv^{(p)} = \vv_\iv^n + \frac{\dt}{2}\bigl(3\Nvh(\vv_\iv\sp n,p_\iv\sp n) - \Nvh(\vv_\iv\sp{n-1},p_\iv\sp{n-1})\bigr)
+ \frac{\dt}{2}\bigl(\Lvh(\vv_\iv\sp{(p)}) + \Lvh(\vv_\iv\sp n)\bigr),
   & \qquad \iv \in \OmegaFh\backslash\GammaIfh, \\
 \tnv_m^T \tauv_{\iv}^{(p)} \nv + \zs \tnv_m^T \vv_\iv^{(p)} =
 \tnv_m^T \Sigmasv_\isv^{(p)} \nv
+ \zs \tnv_m^T \Vsv_{\isv}^{(p)},
 & \qquad\iv \in \GammaIfh, \; \isv \in \GammaIsh,  \\
\Grad\cdot\vv_\iv\sp{(p)}=0,
 & \qquad \iv \in \GammaIfh, \\
 \nv^T\vv_\iv^{(p)}=\f{\zf}{\zf+\zp}\nv^T\Vv^p_h(\vv_\iv\sp{(p)})+ \f{\zp}{\zf+\zp}\nv^T\vsv_\isv^{(p)}, \qquad \tnv_m^T\vv_\iv^{(p)}=\tnv_m^T\Vv\sp{p}_h(\vv_\iv\sp{(p)}), & \qquad \iv \in \GammaIfh, \; \isv \in \GammaIsh, \\
%
\text{Velocity boundary conditions on $\partial \OmegaFh \backslash \GammaIfh$.}
\end{cases}
\medskip\\
4.\text{ Predict fluid pressure:}\\ 
\alignspace 
\begin{cases}
\Delta_h p_\iv^{(p)} = -\rho \Grad \vv_\iv^{(p)} : \bigl(\Grad \vv_\iv^{(p)}  \bigr)^T + \alpha_\iv\Grad \cdot \vv_\iv^{(p)},
&\quad \iv \in \OmegaFh, \\
-p_\iv^{(p)}
-\frac{\zp \dt}{\rho} (\nv\cdot\Grad) p_\iv^{(p)}
= 
\nv^T \bigl(\Sigmasv_\isv^{(p)} \nv - \tauv_\iv^{(p)} \nv\bigr) 
+ \zp \dt
 \nv^T \bigl((\vsv_t)_\isv^{(p)} + \nu \Grad \curl \Grad \curl \vv_\iv^{(p)} \bigr), &\quad \iv \in \GammaIFh, \; \isv \in \GammaIsh, \\
\text{Pressure boundary conditions on $\partial \OmegaFh \backslash \GammaIfh$.}
\end{cases}
\medskip\\
5.\text{ Project solid interface:}\\ 
\alignspace 
\begin{cases}
\nv^T \Vsv_\isv^I = \frac{\zf}{\zf+\zp} \nv^T \vv_\iv^{(p)}
     + \frac{\zp}{\zf+\zp} \nv^T \Vsv_\isv^{(p)} 
     + { \f{1}{\zf+\zp}\bigl( \nv^T\Sigmasv^{(p)}_\iv\nv - \nv^T\Sigmav^{(p)}_\iv\nv \bigr)},
     &\qquad \isv \in \GammaIsh, \iv \in \GammaIfh,\\
\tv_m^T \Vsv_\isv^I = \frac{\zf}{\zf+\zs} \tv_m^T \vv_\iv^{(p)}
     + \frac{\zs}{\zf+\zs} \tv_m^T \Vsv_\isv^{(p)}
     + { \f{1}{\zf+\zs}\bigl( \tv_m^T\Sigmasv^{(p)}_\iv\nv - \tv_m^T\Sigmav^{(p)}_\iv\nv \bigr)},
     &\qquad \isv \in \GammaIsh, \iv \in \GammaIfh,\\
\nv^T \Sigmasv^I_{\isv} \nv = \frac{\zf^{-1}}{\zf^{-1}+\zp^{-1}}\nv^T\Sigmav^{(p)}_\iv \nv
+ \frac{\zp^{-1}}{\zf^{-1}+\zp^{-1}} \nv^T\Sigmasv_\isv^{(p)} \nv
+ { \f{1}{\zf^{-1}+\zp^{-1}}\bigl( \nv^T\vsv^{(p)}_\isv - \nv^T\vv^{(p)}_\iv \bigr)},
     &\qquad \isv \in \GammaIsh, \iv \in \GammaIfh. \\
\tv_m^T \Sigmasv^I_{\isv} \nv = \frac{\zf^{-1}}{\zf^{-1}+\zs^{-1}}\tv_m^T\Sigmav^{(p)}_\iv \nv
    + \frac{\zs^{-1}}{\zf^{-1}+\zs^{-1}} \tv_m^T\Sigmasv_\isv^{(p)} \nv
    + { \f{1}{\zf^{-1}+\zs^{-1}}\bigl( \tv_m^T\vsv^{(p)}_\isv - \tv_m^T\vv^{(p)}_\iv \bigr)},
        &\qquad \isv \in \GammaIsh, \iv \in \GammaIfh,\\
 { \vsv_\isv^{(p)} \leftarrow \vsv_\isv^I, \quad \sigmasv_\isv^{(p)}\nv \leftarrow \sigmasv_\isv^I\nv,}  &\qquad \isv \in \GammaIsh, \iv \in \GammaIfh, \\
\text{Apply solid boundary conditions and set all ghost points.}
\end{cases}
\medskip\\
\hbox{// \textsl{Corrector steps}}\smallskip\\
6.\text{ Correct fluid grid: recompute grid velocity using $\Vsv_\isv^I$ for $\isv \in \GammaIsh$.}
\medskip\\
7.\text{ Correct fluid velocity:} \\
\alignspace 
\begin{cases}
\vv_\iv^{n+1} = \vv_\iv^n + \frac{\dt}{2}\bigl(\Nvh(\vv_\iv^{(p)},p_\iv^{(p)}) + \Nvh(\vv_\iv^{n},p_\iv^{n})\bigr)
+ \frac{\dt}{2}\bigl(\Lvh(\vv_\iv^{n+1}) + \Lvh(\vv_\iv^ n)\bigr),
&\qquad    \iv \in \OmegaFh\backslash\GammaIfh, \\
 \tnv_m^T \tauv_{\iv}^{n+1} \nv + \zs \tnv_m^T \vv_\iv^{n+1} =
 \tnv_m^T \Sigmasv_\isv^{I} \nv
+ \zs \tnv_m^T \Vsv_{\isv}^{I},
 & \qquad\iv \in \GammaIfh, \; \isv \in \GammaIsh,  \\
\Grad\cdot\vv_\iv\sp{n+1}=0,
 & \qquad \iv \in \GammaIfh, \\
 \nv^T\vv_\iv^{n+1}=\f{\zf}{\zf+\zp}\nv^T\Vv_h(\vv_\iv\sp{n+1})+ \f{\zp}{\zf+\zp}\nv^T\vsv_\isv^{I}
, \qquad \tnv_m^T\vv_\iv^{n+1}=\tnv_m^T\vv_\iv^{(e)},  \qquad m = 1,2, & \qquad \iv \in \GammaIfh, \; \isv \in \GammaIsh, \\
\text{Velocity boundary conditions on $\partial \OmegaFh \backslash \GammaIfh$.}
\end{cases}
\medskip\\
8.\text{ Correct fluid pressure.} \\
\alignspace 
\begin{cases}
\Delta_h p_\iv^{n+1} = -\rho \Grad \vv_\iv^{n+1} : \bigl(\Grad \vv_\iv^{n+1}  \bigr)^T + \alpha_\iv\Grad \cdot \vv_\iv^{n+1},
&\quad \iv \in \OmegaFh, \\
-p_\iv^{n+1}-\frac{\zp \dt}{\rho} (\nv\cdot\Grad) p_\iv^{n+1}
= 
\nv^T \bigl(\Sigmasv_\isv^{I} \nv - \tauv_\iv^{n+1} \nv\bigr) 
+ \zp \dt
 \nv^T \bigl((\vsv_t)_\isv^{I} + \nu \Grad \curl \Grad \curl \vv_\iv^{n+1} \bigr), &\quad \iv \in \GammaIFh, \; \isv \in \GammaIsh, \\
\text{Pressure boundary conditions on $\partial \OmegaFh \backslash \GammaIfh$.}
\end{cases}
\medskip\\
9.\text{ Correct solid interface.} \\
\alignspace 
\begin{cases}
 { \Vsv_\isv^{n+1} =  \vv_\iv^{n+1}} , &\qquad \isv \in \GammaIsh, \iv \in \GammaIfh,\\ 
 { \Sigmasv^{n+1}_{\isv} \nv = \Sigmav^{n+1}_\iv \nv, } &\qquad \isv \in \GammaIsh, \iv \in \GammaIfh, \\
%
\text{Reset ghost points corresponding to $\isv \in \GammaIsh$.}
\end{cases}
\\
\end{array}
\]
\label{alg:amp}
\end{algorithm}
}

In the description of the algorithm, we assume that all quantities belonging to the fluid and solid are known at $t=t^n$, and the task is to advance the solution to $t^{n+1}=t^n+\dt$.  Algorithm~\ref{alg:amp} gives a summary of all of the stages in the AMP~scheme.  The first five stages of the algorithm form the {\em predictor} step.

\smallskip\noindent
\textbf{Stage 1 (predict solid)}: 
Predicted values for the solid displacement are obtained by integrating~\eqref{eq:solidDisplacement} using the explicit Lax-Wendroff-type scheme
\begin{align}
\Usv_\isv^{(p)} = \Usv_\isv^n
+ \dt \Vsv_\isv^n
+ \frac{\dt^2}{2\rhos} \Grads \cdot \Sigmasv_\isv^n, 
\qquad \isv \in \OmegaSh.
\label{eq:LaxWendroffSolidDisplacement}
\end{align}
Predicted values for the solid velocity and stress, contained in the vector $\Qsv_\isv=(\vsv_\isv,\sigmasv_\isv)$, are obtained using a second-order accurate Godunov upwind scheme of the form 
\begin{align}
\frac{\Qsv_\isv^{(p)} - \Qsv_\isv^n }{\dt} = 
- \sum_{m=1}^3 \frac{1}{\dxs_m} 
\left( \Fsv_{m, \; \isv}^+ - \Fsv_{m, \; \isv}^- \right), \qquad \isv \in \OmegaSh,
\label{eq:GodunovSolidVariables}
\end{align}
where $\Fsv_{m, \; \isv}^\pm$ represents the upwind fluxes corresponding to the $\xs_m$ coordinate direction.
In practice, we add a stress-relaxation term to the right-hand side of~\eqref{eq:GodunovSolidVariables} following the approach in~\cite{flunsi2016} so that the stress-strain relation in~\eqref{eq:solidStressTensor} is satisfied approximately throughout the solid domain.


\medskip\noindent
\textbf{Stage 2 (predict fluid grid motion)}: 
The computed solid displacement, $\Usv_\isv^{(p)}$, is used to predict the position, $\xv_\iv^{(p)}$ for $\iv \in \GammaIfh$, of the fluid-solid interface at $t^{n+1}$.  The overlapping-grid generator is called to move the fluid grid to match the predicted interface position, and to determine the grid velocity $\dot{\xv}_\iv^{(p)}$.

\medskip\noindent
\textbf{Stage 3 (predict fluid velocity)}: 
The fluid velocity at $t^{n+1}$ is predicted using an explicit Adams-Bashforth scheme for the convection and pressure-gradient terms and an implicit Crank-Nicholson scheme for the viscous terms. The scheme takes the form
\begin{align}
\frac{\vv_\iv^{(p)} - \vv_\iv^n}{\dt} = \frac{1}{2}\left(3\Nvh(\vv_\iv\sp n,p_\iv\sp n) - \Nvh(\vv_\iv\sp{n-1},p_\iv\sp{n-1})\right)
+ \frac{1}{2}\left(\Lvh(\vv_\iv\sp{(p)}) + \Lvh(\vv_\iv\sp n)\right),
\qquad    \iv \in \OmegaFh\backslash\GammaIfh,
\label{eq:VelocityPredictorUpdate}
\end{align}
%
where $\Nvh$ and $\Lvh$ represent grid operators associated with the explicit and implicit terms in~\eqref{eq:VelocityPredictorUpdate}, respectively, given by
\begin{equation}
\Nvh(\vv_\iv,p_\iv) = - \bigl( (\vv_\iv - \dot{\xv}_\iv) \cdot \Grad \bigr) \vv_\iv 
- \frac{1}{\rho} \Grad \,p_\iv ,\qquad
\Lvh(\vv_\iv) = \nu \Delta_h \vv_\iv.
\label{eq:operators}
\end{equation}
Since $\vv_\iv\sp{(p)}$ is determined implicitly in~\eqref{eq:VelocityPredictorUpdate}, boundary conditions are required for $\iv\in\partial\OmegaFh$.  Standard conditions are used for the portion of the boundary away from the interface, while care must be used in the application of the AMP interface conditions for the velocity as described in Section~\ref{sec:interfaceConditions} so that the correct limiting behaviors are obtained.  On the interface, the characteristic condition in~\eqref{eq:AMPTangentialCharacteristic} and the divergence-free constraint give
\bse
\label{eq:AMPVelocityDiscreteBCs}
\begin{alignat}{2}
 \tnv_m^T \tauv_{\iv}^{(p)} \nv + \zs \tnv_m^T \vv_\iv^{(p)} &=
 \tnv_m^T \Sigmasv_\isv^{(p)} \nv
+ \zs \tnv_m^T \Vsv_{\isv}^{(p)}
,  \qquad m = 1,2,
\qquad &&\iv \in \GammaIfh, \; \isv \in \GammaIsh, \label{eq:AMPVelocityDiscrete} \\
\Grad\cdot\vv_\iv\sp{(p)}&=0,
\qquad &&\iv \in \GammaIfh,
\label{eq:DivergenceFreeDiscrete}
\end{alignat}
\ese
where the interface normal, $\nv$, and tangent vectors, $\tnv_m$, are evaluated at grids points $\iv\in\GammaIfh$ along the predicted interface at $t^{n+1}$.
Since~\eqref{eq:VelocityPredictorUpdate} excludes the grids points on the interface, another set of equations are required to complete the system for the implicit time step.  We use the following impedance-weighted averages to set the velocity on the interface:
\bse
\label{eq:AMPinterfaceVelocityDiscrete}
\begin{alignat}{2}
 \nv^T\vv_\iv^{(p)}&=\f{\zf}{\zf+\zp}\nv^T\Vv\sp{p}_h(\vv_\iv\sp{(p)})+ \f{\zp}{\zf+\zp}\nv^T\vsv_\isv^{(p)}, \quad &&\phantom{m = 1,2,}\qquad \iv \in \GammaIfh, \; \isv \in \GammaIsh,
        \label{eq:AMPinterfaceNormalVelocityDiscrete} \\
\tnv_m^T\vv_\iv^{(p)}&=\tnv_m^T\Vv\sp{p}_h(\vv_\iv\sp{(p)}),  \quad &&m = 1,2, \qquad \iv \in \GammaIfh, \; \isv \in \GammaIsh,
             \label{eq:AMPinterfaceTangentialVelocityDiscrete}
\end{alignat}
\ese
where $\Vv\sp{p}_h$ is defined by
\ba
\Vv\sp{p}_h(\vv_\iv\sp{(p)}) \eqdef \vv_\iv^n + \frac{\dt}{2}\left(3\Nvh(\vv_\iv^ n,p_\iv^ n) - \Nvh(\vv_\iv^{n-1},p_\iv^{n-1})\right)
+ \frac{\dt}{2}\Big(\Lvh(\vv_\iv^{(p)}) + \Lvh(\vv_\iv^ n)\Big). \label{eq:veDef}
\ea
Note that~\eqref{eq:AMPinterfaceVelocityDiscrete}, with $\Vv\sp{p}_h$ in~\eqref{eq:veDef}, forms a set of equations for $\vv_\iv\sp{(p)}$, $\iv \in \GammaIfh$, which couples to the evolution equations in~\eqref{eq:VelocityPredictorUpdate} and the boundary conditions in~\eqref{eq:AMPVelocityDiscreteBCs}.  Together, this linear system is solved to obtain the predicted fluid velocity $\vv_\iv\sp{(p)}$.

%

The TP scheme sets the fluid velocity on the interface to equal that of the solid, so that the interface conditions in~\eqref{eq:AMPinterfaceVelocityDiscrete} are replaced by
\begin{align}
\vv_\iv^{(p)} = \Vsv_\isv^{(p)}, \quad \iv \in \GammaIfh, \; \isv \in \GammaIsh.
\label{eq:TPVelocityDiscrete}
\end{align}
The divergence-free constraint in~\eqref{eq:DivergenceFreeDiscrete} is used to specify the normal component of the fluid velocity in the ghost points, while the tangential components of the velocity in the ghost points are extrapolated.  The ATP scheme, on the other hand, applies a tangential stress boundary condition on the fluid along with the divergence-free constraint.  These conditions are equivalent to the conditions in~\eqref{eq:AMPVelocityDiscreteBCs}, but with the terms in~\eqref{eq:AMPVelocityDiscrete} involving the fluid and solid velocity dropped (effectively setting $\zs=0$), and then replacing the interface conditions in~\eqref{eq:AMPinterfaceVelocityDiscrete} with $\vv_\iv^{(p)}=\Vv\sp{p}_h(\vv_\iv\sp{(p)})$ for $\iv\in\GammaIfh$.  We note that these limiting cases for TP and ATP schemes correspond to the limits of the continuous conditions discussed above in Section~\ref{sec:interfaceConditions}.

\medskip\noindent
\textbf{Stage 4 (predict fluid pressure):}
Predicted values for the fluid pressure are obtained by solving the discrete Poisson problem
\begin{align}
\Delta_h p_\iv^{(p)} = -\rho \Grad \vv_\iv^{(p)} : \Bigl(\Grad \vv_\iv^{(p)}  \Bigr)^T + \alpha_\iv\Grad \cdot \vv_\iv^{(p)},
\qquad \iv \in \OmegaFh,
\label{eq:predictPressure}
\end{align}
with the predicted fluid velocity from the previous stage used to determine the right-hand side of~\eqref{eq:predictPressure}.  We have included a {\em divergence damping} term to the right-hand side of~\eqref{eq:predictPressure} with coefficient~$\alpha_\iv$ to suppress the growth of the divergence of the velocity in the numerical solution that may occur due to discretization errors.  The discrete approximation of the Robin condition in~\eqref{eq:AMPpressureBC} is
\begin{align}
-p_\iv^{(p)}
-\frac{\zp \dt}{\rho} (\nv\cdot\Grad) p_\iv^{(p)}
&= 
\nv^T \left(\Sigmasv_\isv^{(p)} \nv - \tauv_\iv^{(p)} \nv\right) 
\nonumber \\
&\qquad + \zp \dt
 \nv^T \left((\vsv_t)_\isv^{(p)} + \nu \Grad \curl \Grad \curl \vv_\iv^{(p)}\right), \qquad \iv \in \GammaIFh, \; \isv \in \GammaIsh,
\label{eq:discreteAMPpressureBC}
\end{align}
where the acceleration of the solid along the interface appearing on the right-hand side of~\eqref{eq:discreteAMPpressureBC} is obtained from the solid velocity using a second-order accurate finite difference approximation in time.

In the TP scheme, the Robin condition in~\eqref{eq:discreteAMPpressureBC} is replaced by a discrete approximation of~\eqref{eq:pressureBCHeavySolid} given by
\[
-{1\over\rho}(\nv\cdot\Grad) p_\iv^{(p)} = \nsv^T \left((\vsv_t)_\isv^{(p)} + \nu \Grad \curl \Grad \curl \vv_\iv^{(p)} \right), \qquad \iv \in \GammaIFh, \; \isv \in \GammaIsh,
\]
whereas the ATP scheme uses an approximation of~\eqref{eq:normalStressLightSolid} given by
\[
-p_\iv^{(p)} = \nv^T \left(\Sigmasv_\isv^{(p)} \nv - \tauv_\iv^{(p)} \nv\right), \qquad \iv \in \GammaIFh, \; \isv \in \GammaIsh.
\]
As before, the both cases are consistent with the limiting cases of a heavy (TP) and light (ATP) solid.


\medskip\noindent
\textbf{Stage 5 (predict solid interface):}  The projections given by~\eqref{eq:AMPinterfaceVelocity} and~\eqref{eq:AMPinterfaceTraction} are used in this stage to define values for the solid velocity and traction on the interface $\isv \in \GammaIsh$ ($\iv \in \GammaIfh$) as
\bse
\begin{alignat}{2}
\nv^T \Vsv_\isv^I =\;& \frac{\zf}{\zf+\zp} \nv^T \vv_\iv^{(p)}
     + \frac{\zp}{\zf+\zp} \nv^T \Vsv_\isv^{(p)}
     + { \f{1}{\zf+\zp}\Big( \nv^T\Sigmasv^{(p)}_\isv\nv - \nv^T\Sigmav^{(p)}_\iv\nv \Big)},\\
\tv_m^T \Vsv_\isv^I =\;& \frac{\zf}{\zf+\zs} \tv_m^T \vv_\iv^{(p)}
     + \frac{\zs}{\zf+\zs} \tv_m^T \Vsv_\isv^{(p)}
     + { \f{1}{\zf+\zs}\Big( \tv_m^T\Sigmasv^{(p)}_\isv\nv - \tv_m^T\Sigmav^{(p)}_\iv\nv \Big)},
\end{alignat}
\ese
and
\bse
\begin{alignat}{2}
\nv^T \Sigmasv^I_{\isv} \nv =\;& \frac{\zf^{-1}}{\zf^{-1}+\zp^{-1}}\nv^T\Sigmav^{(p)}_\iv \nv
+ \frac{\zp^{-1}}{\zf^{-1}+\zp^{-1}} \nv^T\Sigmasv^{(p)}_\isv\nv
+ { \f{1}{\zf^{-1}+\zp^{-1}}\Big( \nv^T\vsv^{(p)}_\isv - \nv^T\vv^{(p)}_\iv \Big)},\\
\tv_m^T \Sigmasv^I_{\isv} \nv =\;& \frac{\zf^{-1}}{\zf^{-1}+\zs^{-1}}\tv_m^T\Sigmav^{(p)}_\iv \nv
    + \frac{\zs^{-1}}{\zf^{-1}+\zs^{-1}} \tv_m^T\Sigmasv^{(p)}_\isv\nv
    + { \f{1}{\zf^{-1}+\zs^{-1}}\Big( \tv_m^T\vsv^{(p)}_\isv - \tv_m^T\vv^{(p)}_\iv \Big)},
\end{alignat}
\ese
respectively.  These projected values are used to specify $\Vsv_\isv^{(p)}$ and $\Sigmasv^{(p)}_{\isv} \nv$ on the interface, and then the corresponding values in the ghost points are set using extrapolation. 

In the TP scheme, we use $\Sigmasv^{(p)}_{\isv} \nv=\Sigmav^{(p)}_\iv \nv$ to set the predicted solid traction on the interface,
but then use the compatibility condition
\begin{align}
\lambdas \left(\Grads \cdot \Vsv_\isv^{(p)} \right) \nsv
+ \mus \left(\Grads \Vsv_\isv^{(p)} + \left(\Grads \Vsv_\isv^{(p)} \right)^T  \right) \nsv
 &= \dot{\Sigmasv}_{\isv} \nsv,
\qquad \isv \in \GammaIsh ,
\label{eq:predictedSolidInterfaceVelocityCompatibility}
\end{align}
which can be interpreted as a set of discrete Neumann conditions for the components of $\Vsv_\isv^{(p)}$.  The time-derivative of the solid stress on the right-hand side of~\eqref{eq:predictedSolidInterfaceVelocityCompatibility} is obtained using a second-order accurate approximation involving $\Sigmasv_{\isv}^{(p)}$ and the solid stress at $t^n$ and $t^{n-1}$.  For the ATP scheme, we set the solid velocity on the interface using $\Vsv_\isv^{(p)}=\vv_\iv^{(p)}$, and use the compability condition
\begin{align}
\Grads \cdot \Sigmasv_\isv^{(p)} = \rho \dot{\Vsv}_\isv\sp{(p)},
\qquad \isv \in \GammaIsh, \label{eq:predictedSolidStressGhostBC}
\end{align}
which can be interpreted as a set of discrete Neumann conditions for the components of $\Sigmasv_\isv^{(p)}\nv$.  The acceleration of the solid on the interface is computed using a finite-difference approximation similar to that described to compute the stress for the TP scheme.

\medskip
This completes the set of stages used for the predictor step of Algorithm~\ref{alg:amp}.
The remaining four stages form the {\em corrector} step.
Stages~7 and~8 in the corrector step are similar in form to the corresponding two stages in the predictor step.
In Stage~9 the solid velocity and traction are set equal to the current fluid values.
This ensures that the interface jump conditions are exactly matched at the end of the correction step.
The boundary conditions for the TP and ATP scheme in the corrector stages are handled in an analogous
fashion to the predictor stages.

\section{Spatial approximations} \label{sec:numericalApproach}

\begin{figure}[hbt]
\begin{center}
\begin{tikzpicture}[scale=.7]
\useasboundingbox (.75,1.25) rectangle (15.5,6);  
%
\begin{scope}[xshift=1cm,yshift=1cm]
\fill[black!10!white,xshift=.5cm,yshift=.5cm] (0,0) -- (2.583333,0) arc (180:90:1.416667) -- (4.,4.) -- (0,4.) -- (0,0);
\draw[-,thick,blue,yshift=.0 cm] 
   \foreach \x/\y in {1.5/0,1.5/.5,2/1,2/1.5,2.5/2,3/2.5,4/3,5/3.5,5/4,5/4.5,5/5}{ (0,\y) -- (\x,\y) }
   \foreach \x/\y in {0/0,.5/0,1/0,1.5/0,2/1,2.5/2,3/2.5,3.5/3,4/3,4.5/3.5,5/3.5}{ (\x,\y) -- (\x,5) };
  \begin{scope}[xshift=4.5cm,yshift=0.5cm]
    \draw[thick,green] \foreach \r in {1.000000,1.416667,1.833333,2.250000,2.666667,3.083333,3.500000}{ (0,\r) arc (90:190:\r)  (0,\r) arc (90:80:\r) };
    \draw[thick,green]
     (0.173648,0.984808)  -- (0.607769,3.446827)
     (0.000000,1.000000)  -- (0.000000,3.500000)
     (-0.173648,0.984808) -- (-0.607769,3.446827)
     (-0.342020,0.939693) -- (-1.197071,3.288924)
     (-0.500000,0.866025) -- (-1.750000,3.031089)
     (-0.642788,0.766044) -- (-2.249757,2.681156)
     (-0.766044,0.642788) -- (-2.681156,2.249757)
     (-0.866025,0.500000) -- (-3.031089,1.750000)
     (-0.939693,0.342020) -- (-3.288924,1.197071)
     (-0.984808,0.173648) -- (-3.446827,0.607769)
     (-1.000000,0.000000) -- (-3.500000,0.000000)
     (-0.984808,-0.173648) -- (-3.446827,-0.607769);
  \end{scope}
  \draw[very thick,red,xshift=.5cm,yshift=.5cm] (0,0) -- (2.583333,0) arc (180:90:1.416667) -- (4.,4.) -- (0,4.) -- (0,0);
%
   \filldraw[green] (1.5,.5)  circle (3pt)
                 (1.5,1 )  circle (3pt)
                 (2  ,1 )  circle (3pt)
                 (2 ,1.5)  circle (3pt)
                 (2 , 2 )  circle (3pt)
                 (2.5,2 )  circle (3pt)
                 (2.5,2.5) circle (3pt)
                 (3 , 2.5) circle (3pt)
                 (3  ,3 )  circle (3pt)
                 (3.5,3 )  circle (3pt)
                 (4  ,3. ) circle (3pt)
                 (4  ,3.5) circle (3pt)
                 (4.5,3.5) circle (3pt);
%
  \begin{scope}[xshift=4.5cm,yshift=0.5cm]
      \filldraw[blue]
       (0.000000,3.500000)    circle (3pt)
       (-0.607769,3.446827)   circle (3pt)
       (-1.197071,3.288924)  circle (3pt) 
       (-1.750000,3.031089)  circle (3pt) 
       (-2.249757,2.681156)  circle (3pt) 
       (-2.681156,2.249757)  circle (3pt) 
       (-3.031089,1.750000)  circle (3pt) 
       (-3.288924,1.197071)  circle (3pt) 
       (-3.446827,0.607769)  circle (3pt) 
       (-3.500000,0.000000)  circle (3pt);
  \end{scope}
   \draw (1.25,3.5) node[thick,draw=blue,fill=white] {\large$G_1$};
   \draw (3.25,1.95) node[thick,draw=green,fill=white] {\large$G_2$};
\end{scope}
%
\definecolor{ghostColour}{named}{DodgerBlue}
\newcommand{\mytrix}{(\x,-.15) -- ++(.3,0) -- ++(-.15,.26) -- (\x,-.15)}
\newcommand{\mytriy}{(-.15,\y) -- ++(.3,0) -- ++(-.15,.26) -- (-.15,\y)}
\begin{scope}[xshift=7cm,yshift=2.25cm,scale=.75]
\draw[-,thick,blue,yshift=.0 cm] 
   \foreach \x in {0,.5,...,5}{ (\x,0) -- (\x,5) }
   \foreach \y in {0,.5,...,5}{ (0,\y) -- (5,\y) };
  \draw[very thick,red,xshift=.5cm,yshift=.5cm] (1.,0) -- (.0,0) -- (.0,4.) -- (4.,4.) -- (4.,3.);
   \filldraw[green] (1.5,.5)  circle (3pt)
                 (1.5,1 )  circle (3pt)
                 (2  ,1 )  circle (3pt)
                 (2 ,1.5)  circle (3pt)
                 (2 , 2 )  circle (3pt)
                 (2.5,2 )  circle (3pt)
                 (2.5,2.5) circle (3pt)
                 (3 , 2.5) circle (3pt)
                 (3  ,3 )  circle (3pt)
                 (3.5,3 )  circle (3pt)
                 (4  ,3. ) circle (3pt)
                 (4  ,3.5) circle (3pt)
                 (4.5,3.5) circle (3pt);
  \filldraw[fill=white,draw=black]  \foreach \x in {2,2.5,...,5}{ (\x,.0) circle (3.5pt) };
  \filldraw[fill=white,draw=black]  \foreach \x in {2,2.5,...,5}{ (\x,.5) circle (3.5pt) };
  \filldraw[fill=white,draw=black]  \foreach \x in {2.5,3,...,5}{ (\x,1.) circle (3.5pt) };
  \filldraw[fill=white,draw=black]  \foreach \x in {2.5,3,...,5}{ (\x,1.5) circle (3.5pt) };
  \filldraw[fill=white,draw=black]  \foreach \x in {3,3.5,...,5}{ (\x,2.0) circle (3.5pt) };
  \filldraw[fill=white,draw=black]  \foreach \x in {3.5,4,...,5}{ (\x,2.5) circle (3.5pt) };
  \filldraw[fill=white,draw=black]  \foreach \x in {4.5,5}      { (\x,3.0) circle (3.5pt) };
  \draw[fill=ghostColour,xshift=-.15cm,yshift=0cm]  \foreach \x in {.5,1.,1.5}{ \mytrix };  
  \draw[fill=ghostColour,xshift=-.15cm,yshift=5cm]  \foreach \x in {.5,1.,...,5}{ \mytrix };  
  \draw[fill=ghostColour,xshift=0cm,yshift=-.15cm]  \foreach \y in {0,.5,...,5}{ \mytriy };
  \draw[fill=ghostColour,xshift=5cm,yshift=-.15cm]  \foreach \y in {3.5,4,4.5}{ \mytriy };
   \draw (1.25,3.5) node[thick,draw=blue,fill=white] {\large$G_1$};
\end{scope}
\begin{scope}[xshift=11.5cm,yshift=2.25cm,scale=.75]
\draw[-,thick,green,yshift=.0 cm] 
   \foreach \x in {0,.454545,...,5}{ (\x,0) -- (\x,5) }
   \foreach \y in {0,.833333,...,5}{ (0,\y) -- (5,\y) };
 \draw[very thick,red,xshift=.454545cm,yshift=.833333cm] (0.,4) -- (.0,0) -- (4.0909,0.) -- (4.0909,4);
 \filldraw[blue]  \foreach \x in {.454545,.909090,...,4.545454}{ (\x,5) circle (3.5pt) };
 \draw[fill=ghostColour,xshift=-.15cm]  \foreach \x in {.454545,.909090,...,4.545454}{ \mytrix };
 \draw[fill=ghostColour,yshift=-.15cm]  \foreach \y in {0,.833333,...,5}{ \mytriy };
 \draw[fill=ghostColour,xshift=5cm,yshift=-.15cm]  \foreach \y in {0,.833333,...,5}{ \mytriy };
\end{scope}
\begin{scope}[xshift=7cm,yshift=.7cm]
  \fill[black!10!white,xshift=-.1cm,yshift=-.25cm] (0,0) -- (4,0) -- (4.,1.3) -- (0,1.3) -- (0,0);
  \filldraw[green,xshift=.0cm,yshift=.8cm] (.25,.0)  circle (3pt);
  \filldraw[blue,xshift=.3cm,yshift=.8cm] (.25,.0)  circle (3pt);
  \draw[xshift=.0cm,yshift=.8cm] (.5,0) node[anchor=west,xshift=6] {\small interpolation};
  \draw[fill=ghostColour,xshift=.0cm,yshift=.4cm] (.35,0) \foreach \x in {.1}{ \mytrix } node[anchor=west,xshift=12,yshift=3] {\small ghost};
  \draw[fill=white,draw=black,xshift=.0cm,yshift=.0cm] (.25,0) circle (3.5pt) node[anchor=west,xshift=6] {\small unused};
\end{scope}
\begin{scope}[xshift=11.5cm,yshift=2.25cm,scale=.75]
   \draw (1.6,3.27) node[thick,draw=green,fill=white] {\large$G_2$};
\end{scope}
\end{tikzpicture}
\end{center}
\caption{Left: composite of a background grid ($G_1$, blue) 
  and a boundary-fitted grid ($G_2$, green) in physical space for the domain defined by
  the interior of the red boundary. The grid points on $G_1$ with green dots interpolate from
  $G_2$ and the grid points on $G_2$ with blue dots interpolate from $G_1.$
  Middle: Plot of $G_1$ showing interpolation points, ghost points (grid points which exist 
  outside the physical boundary), and unused points (grid points which do not affect the computation).
  Right: The green boundary fitted grid, $G_2,$ is mapped to a unit square. The plot shows
  interpolation points and ghost points.
\label{fig:overlappingGrid}}
\end{figure}
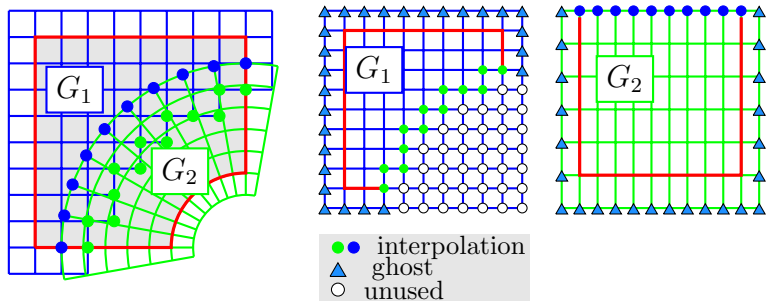

The domain for the FSI problem is discretized using an overlapping grid $\Gc$, which is defined by a set of component grids $\{G_g\}$, $g = 1,\, \ldots,\, \Nc$, that cover both the fluid and solid domains.  Static Cartesian background grids are used for the bulk of each domain, while boundary and interface-fitted grids are used to resolve curved boundaries and moving fluid-solid interfaces (see Figure~\ref{fig:overlappingGrid}).  In the fluid, each component grid is defined by a mapping from the physical space $\xv$ to a unit computational space $\rv$ (a square in two dimensions or a box in three dimensions).  The mapping is time dependent in general to account for the moving fluid-solid interface, and is given by
\begin{align}
\xv = \Gv_g(\rv,t), \qquad \rv \in [0,1]^{\nd}, \quad \xv \in \Real^{\nd},
\end{align}
where $\nd$ is the number of dimensions.  Each component grid for the solid is defined by a mapping from the (fixed) reference space $\xsv$ to a unit computational space, and is given by
\begin{align}
\xsv = \bar{\Gv}_g(\rv), \qquad \rv \in [0,1]^{\nd}, \quad \xsv \in \Real^{\nd}.
\end{align}
The corresponding physical space $\xv$ for the solid is then determined by the mapping in~\eqref{eq:solidPhysicalCoordinates} using the computed displacement $\usv(\xsv,t)$.  Within the fluid and solid domains, the component grids may overlap and interpolation formulas are used to match the discrete solutions across the component grids.  The position of the fluid-solid interface is determined by the interface-fitted component grids of the solid, and then the interface-fitted fluid grids deform to match the moving interface position in physical space.  A hyperbolic grid generator~\cite{HyperbolicGuide} is used to regenerate the fluid interface-fitted grids, and then the {\tt Ogen} grid generator~\cite{OGEN} is called to regenerate the overlapping grid connectivity between the (moving) interface-fitted grids of the fluid and the (static) background grids.

The Navier-Stokes equations for the fluid in velocity-pressure form are transformed to the unit-square coordinates, $\rv$, using standard chain-rule differentiation formulas.  On moving grids, the equations are transformed to a moving coordinate system which introduces the grid velocity, $\dot{\xv}$, into the advection terms as indicated in~\eqref{eq:VelocityPredictorUpdate} and~\eqref{eq:operators}.  The resulting equations are then discretized using standard finite-difference approximations for the derivatives with respect to $\rv$ following the approach in~\cite{ICNS,splitStep2003,mog2006}. The equations of linear elasticity for the solid are also transformed to the unit-square coordinates and then integrated in time using a (Godunov-based) upwind scheme similar to \cite{smog2012,fsi2012,flunsi2016}.  The approximations of the equations governing the fluid and solid are second-order accurate in the present implementation of the AMP time-stepping scheme.

The time step, $\Delta t$, is taken as the minimum value between the stable time steps for the fluid and solid domain solvers independently.  There is no restriction on the time step imposed by the AMP interface conditions.  For the fluid solver, the time step is generally determined by a CFL-type stability constraint based on the advection terms in the fluid momentum equations.  The viscous terms are integrated implicitly so that there is no stability constraint on $\Delta t$ arising from these terms.  The time step for the solid solver is determined based on a CFL-type stability constraint involving the compression wave speed, but there also may be an adjustment imposed by the stress-relaxation term (see~\cite{smog2012,flunsi2016}).

\section{Numerical results} \label{sec:numericalResults}

We now present numerical results to 
verify the stability and accuracy of the AMP scheme for a variety of FSI problems.
We start with fundamental problems for which exact solutions are known. 
The first problem involves an annular fluid region surrounding an elastic solid disk
with the motion of the solid constrained to the radial direction normal to the fluid-solid interface. 
The second problem uses the same geometry, but assumes circumferential motion of the solid.
Exact solutions exist for both problems, as discussed in~\ref{sec:radialElasticPistonExactSolution}
and~\ref{sec:rotatingElasticDiskExactSolution}, while the results are given in Sections~\ref{sec:radialElasticPiston}
and~\ref{sec:rotatingElasticDisk} below.
We next test the AMP scheme for a more general FSI problem involving a radial traveling wave.  Results for this
problem are given in Section~\ref{sec:travelingWave} with the corresponding exact solution for a linearized
traveling wave derived in~\ref{sec:radialTravelingWaveSolution}.
For each of these test cases, we examine the convergence of the numerical solution using the exact solution
for different fluid-solid density ratios to confirm that the AMP scheme is both stable and second-order accurate.
The remaining two tests involve more complex FSI problems for which exact solutions are unavailable.  The first
problem, discussed in Section~\ref{sec:elasticDiskInAChannel}, considers a pressure-driven flow in a channel about
an elastic solid.  A self-convergence study is performed for this clean benchmark problem which further verifies the
accuracy and stability of the AMP scheme.  We also examine the behavior of the TP scheme for this problem.  The final
test problem involves the flow past multiple solids with a wide range of densities.  This test is discussed in
Section~\ref{sec:multipleBodies} and is designed to illustrate the robustness of the AMP scheme.

\subsection{Radial elastic piston} \label{sec:radialElasticPiston}

Let us begin by considering two benchmark FSI problems for which exact solutions are available.  The geometry of the two problems, shown in Figure~\ref{fig:radialElasticPistonGeometry}, is similar.  A solid elastic disk with initial radius $\bar r=r_I(0)$ is surrounded by an incompressible fluid which occupies the annular region $r_I(t)<r<R$.  In the first problem, it is assumed that the motion of the solid and fluid is confined to the radial direction, while circumferential motion is assumed in the second problem.

{
\newcommand{\lbFont}{\small}
\def\rI{1.4}
\def\rIDef{1.7}
\def\rF{3.7}
\def\ep{.15}
\def\thetaA{210}
\def\thetaB{-30}
\def\xR{9}
\def\gap{.5}
\begin{figure}[hbt]
	\newcommand{\textFont}{\normalss}
	\begin{center}
            \resizebox{10cm}{!}{
		\begin{tikzpicture}[scale=.75]
                  \useasboundingbox (-5,-3) rectangle (12.5,4);  
                  
		%
                \path [draw=blue,fill=blue!20,even odd rule,line width=2pt] (0,0) circle (\rF) (0,0) circle (\rI);
                \path [draw=red,fill=red!20,line width=2pt] (0,0) circle (\rI-\ep);
                
                \draw [thick,black,anchor=south] (0,0) node[yshift=-4pt] {\lbFont solid: $\bar{\Omega}(0)$};
                \draw [thick,black,anchor=south] (0,\rI) node {\lbFont interface: $\Gamma(0)$};
                \draw [thick,black] (0,-.5*\rI-.5*\rF) node {\lbFont fluid: $\OmegaF(0)$};

                \draw [draw=black,line width=2pt,->,anchor=north west] (0,0) -- (\thetaB:\rI-\ep) node {\lbFont$\rs=R_I(0)$};
                \draw [draw=black,line width=2pt,->,anchor=north east] (0,0) -- (\thetaA:\rF) node {\lbFont$\rs=R_f$};

                \path [draw=blue,fill=blue!20,even odd rule,line width=2pt] (\xR,0) circle (\rF) (\xR,0) circle (\rIDef);
                \path [draw=red,fill=red!20,line width=2pt] (\xR,0) circle (\rIDef-\ep);
                
                \draw [thick,black,anchor=south] (\xR,0) node {\lbFont solid: $\bar{\Omega}(t)$};
                \draw [thick,black,anchor=south] (\xR,\rIDef) node {\lbFont interface: $\Gamma(t)$};
                \draw [thick,black] (\xR,-.5*\rIDef-.5*\rF) node {\lbFont fluid: $\OmegaF(t)$};

                \draw [draw=black,line width=2pt,->,anchor=north west] (\xR,0) -- +(\thetaB:\rIDef-\ep) node {\lbFont$r=R_I(t)$};
                \draw [draw=black,line width=2pt,->,anchor=north east] (\xR,0) -- +(\thetaA:\rF) node {\lbFont$r=R_f$};


		\end{tikzpicture}
          }
	\end{center}
	\caption{The geometry for the radial elastic piston solution.
          The left figure shows the initial domain without deformation
          and the right figure shows the deformed domain for $t > 0.$
        } \label{fig:radialElasticPistonGeometry}
\end{figure}
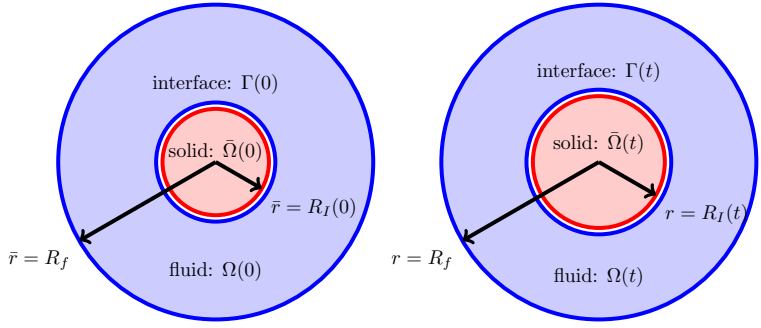
}

\subsubsection{Radial motion of a circularly symmetric elastic piston}\label{sec:radialElasticPiston}


For the case of radial motion, the solid disk behaves as an elastic piston that drives a radial flow in the fluid surrounding it (as illustrated in Figure~\ref{fig:radialElasticPistonGeometry}).
The exact solution for this problem is derived in~\ref{sec:radialElasticPistonExactSolution}.  It is convenient to define the
solution in terms of a specified radial component of the solid displacement, $\us_r(\rs,t)$, having the form
\begin{equation}
\us_r(\rs,t)=\beta J_1(\omega\rs/\cp)\sin(\omega t),
\label{eq:sDisplacement}
\end{equation}
where $\omega$ is the frequency of the oscillation, $\beta$ determines the amplitude, and $J_1(z)$ is the Bessel function of the first kind of order one.  The corresponding time-dependent interface radius, $r_I(t)$, is given by 
\begin{equation}
r_I(t)=r_0+\us_r(r_0,t)=r_0+b\sin(\omega t),\qquad b=\beta J_1(\omega r_0/\cp),
\label{eq:pistionDisplacement}
\end{equation}
where $r_0 = r_I(0)$ is the initial interface radius. The radial component of the fluid velocity and its pressure are given by
\bas
& v_r(r,t) = \frac{R}{r} V(t), \\
& p(r,t) = P(t)+{\rho\over2}\left(1-{R\sp2\over r\sp2}\right)V(t)\sp2+\rho R\log\left({R\over r}\right)\dot{V}(t), 
\eas
where
\bas
V(t)=\;& \frac{r_I(t)}{R} \omega\beta J_1(\omega r_0/\cp) \cos(\omega t) , \\
P(t)=\;& -{\rho\over2}\left(1-{R\sp2\over r_I(t)\sp2}\right)V(t)\sp2-\rho R\log\left({R\over r_I(t)}\right)\dot{V}(t) \\
&\qquad -\beta\left[(\lambdas+2\mus){\omega\over\cp}J_1\sp\prime(\omega r_0/\cp)+{\lambdas\over r_0}J_1(\omega r_0/\cp)\right]\sin(\omega t).
\eas
Here, $V(t)$ and $P(t)$ are the radial velocity and pressure of the fluid at $r=R$, respectively, and $\dot{V}(t)$ is the radial acceleration.  These functions are specified in the exact solution based on the choice for the interface motion given in~\eqref{eq:pistionDisplacement}.
We take $\beta=0.05$ to be the specified amplitude of the motion 
and $\omega=\pi$ as the specified frequency. 
Additionally, we choose $r_0=0.5$, $R=1$, $\rho=1$, and $\rhos=\mus=\lambdas=\delta$, where $\delta$ is a parameter that determines the ratio of the density of the solid to that of the fluid.

\newcommand{\Gcrep}{\Gc_{\rm rep}}

{
%
\newcommand{\drawContour}[9]{%
\begin{scope}[#1]
\draw(0.0,0) node[anchor=south west,xshift=-4pt,yshift=+0pt] {\trimfig{fig/#2}{\figWidth}};
  \draw(1.3,1.3) node[draw,fill=white,anchor=west,xshift=2pt,yshift=1pt] {\scriptsize #3};
  \draw(2.2,2.2) node[draw,fill=white,anchor=west,xshift=2pt,yshift=1pt] {\scriptsize #4};
  \draw(.5,5) node[draw,fill=white,anchor=west,xshift=2pt,yshift=-3pt] {\scriptsize #5}; 
\begin{scope}[xshift=5pt,yshift=-7pt]
  \draw (\xcb,\ycb) node[anchor=south west,xshift=0.25cm,yshift=.5cm,rotate=-90] {\trimfigcb{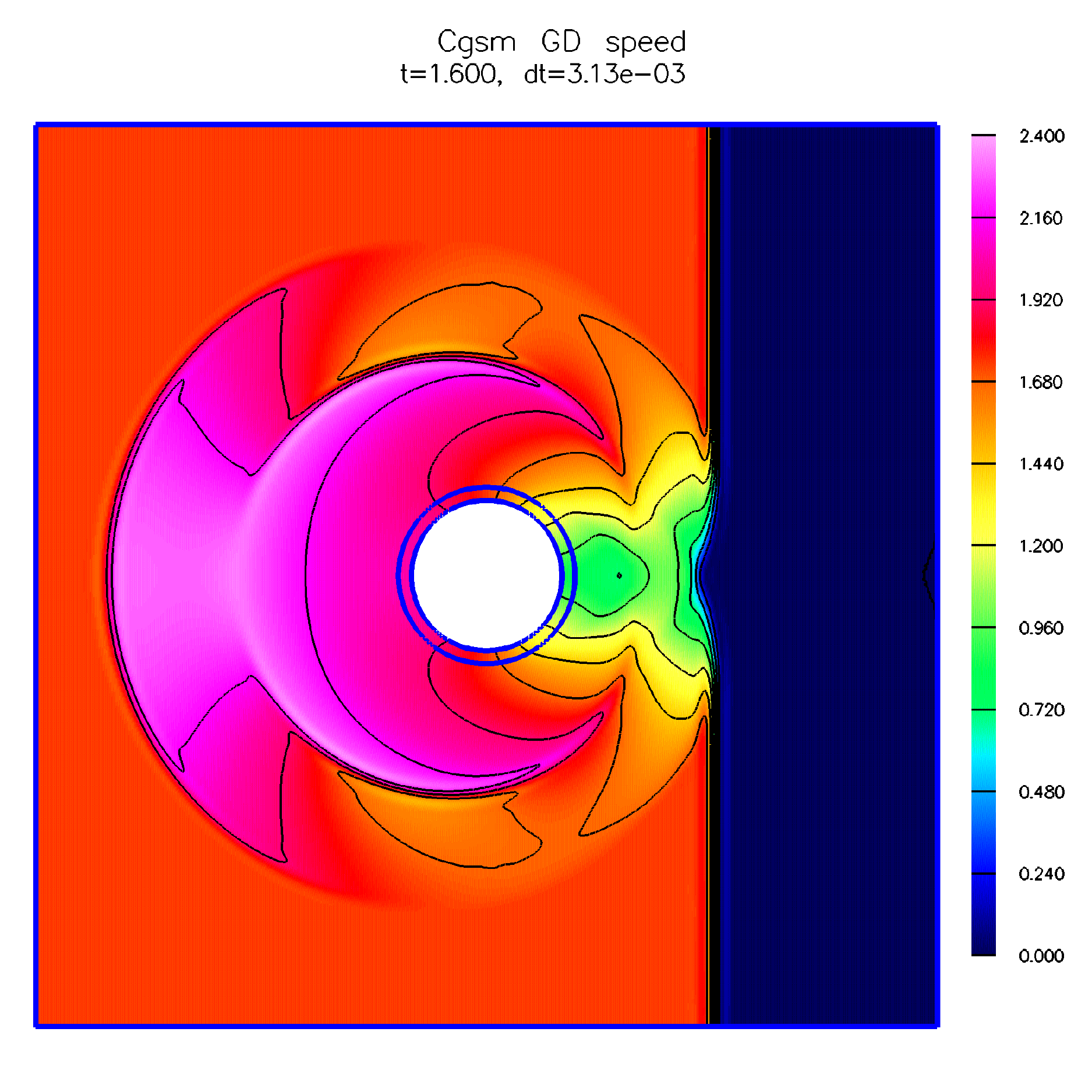}{\cbWidth}{\cbHeight}};
  \draw (.8,0) node[anchor=north,xshift=+3pt,yshift=+2pt] {\scriptsize $#6$};
  \draw (2.8,0) node[anchor=north,xshift=+0pt,yshift=+2pt] {\scriptsize fluid};
  \draw (4.8,0) node[anchor=north,xshift=+0pt,yshift=+2pt] {\scriptsize $#7$};
  \begin{scope}[yshift=14pt]
    \draw (.8,0) node[anchor=north,xshift=+3pt,yshift=+2pt] {\scriptsize $#8$};
    \draw (2.8,0) node[anchor=north,xshift=+0pt,yshift=+2pt] {\scriptsize solid};
    \draw (4.8,0) node[anchor=north,xshift=+0pt,yshift=+2pt] {\scriptsize $#9$};
 \end{scope}
\end{scope}
\end{scope}
}
\newcommand{\cbWidth}{.2cm}
\newcommand{\cbHeight}{4cm}
\newcommand{\xcb}{.5cm}
\newcommand{\ycb}{-.2cm}
\setlength{\ycbTop}{\ycb+\cbHeight}
\setlength{\ycbMid}{\ycb+\cbHeight*\real{.5}}
\newcommand{\trimfigcb}[3]{\includegraphics[width=#2, height=#3, clip, trim=17cm 2.35cm 1.65cm 2.35cm]{#1}}
\newcommand{\figWidth}{5cm}
\newcommand{\trimfig}[2]{\trimh{#1}{#2}{.11}{.11}{.11}{.11}}
\begin{figure}[htb]
\begin{center}
\begin{tikzpicture}[scale=1]
  \useasboundingbox (0,0.25) rectangle (15.5,5.25);  

   \draw(0.0,0.0) node[anchor=south west,xshift=-15pt,yshift=0pt] {\trimfig{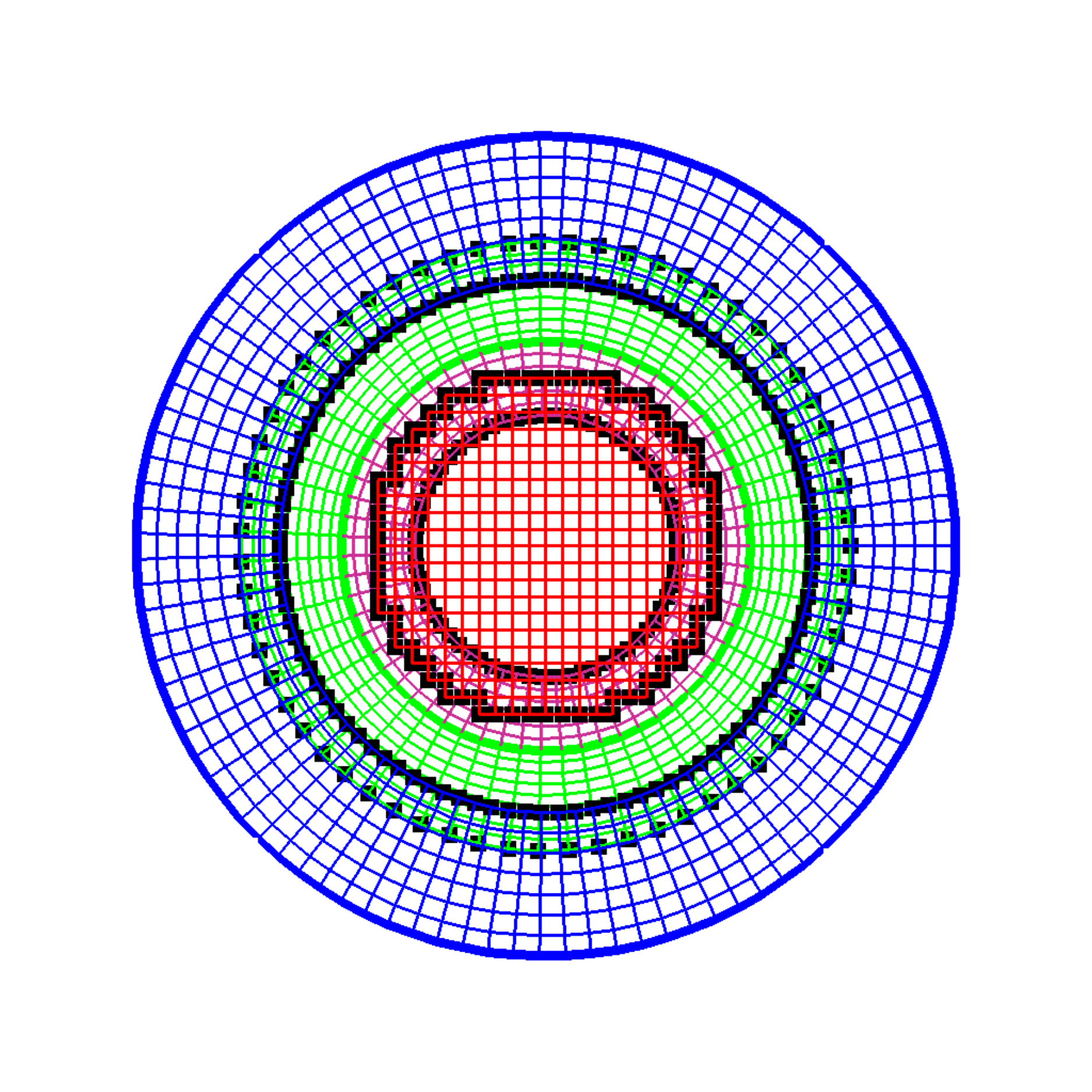}{\figWidth}};
   \drawContour{xshift= 5.cm,yshift=0cm}{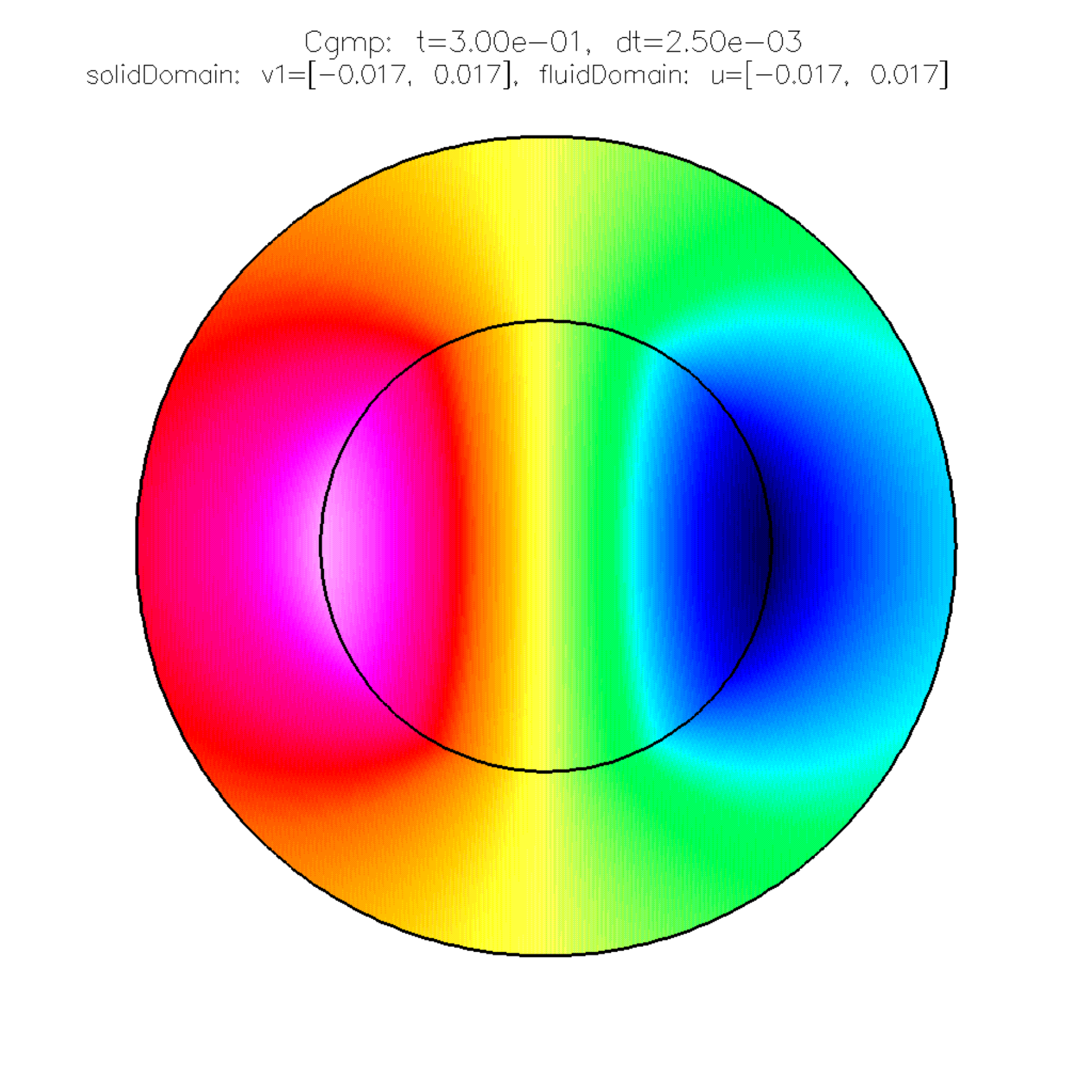}{$v_1$}{$\vs_1$}{$t=0.3$}{$-.015$}{$.017$}{$-.017$}{$.015$};
   \drawContour{xshift=10.cm,yshift=0cm}{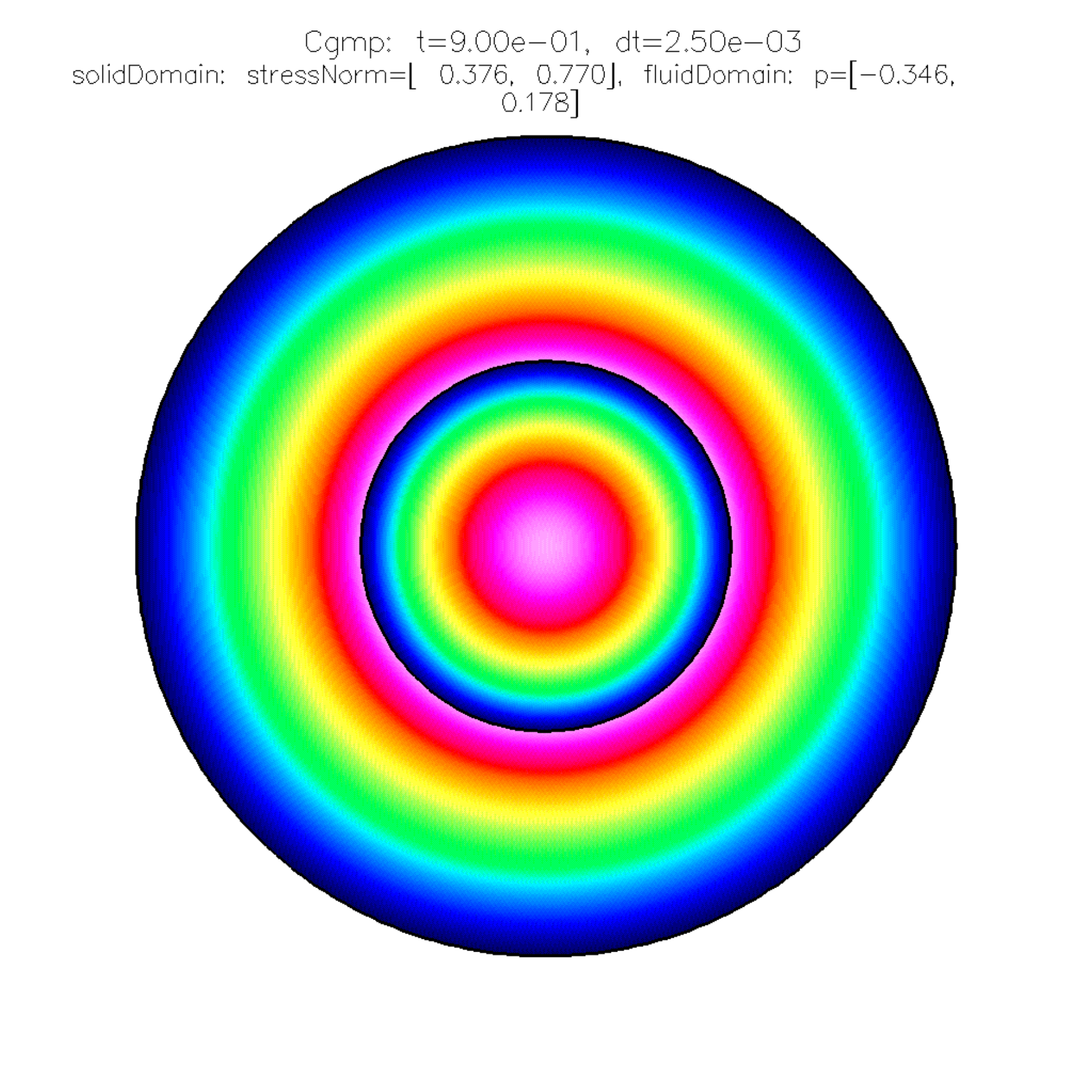}{$p$}{$|\sigmasv|$}{$t=0.9$}{$-.35$}{$.18$}{$.38$}{$.77$};
\end{tikzpicture}
\end{center}
  \caption{Left: composite grid $\Gcrep^{(2)}$ for the radial elastic piston problem.
    The green fluid grid adjacent to the interface  moves over time. The blue fluid background grid remains fixed.
    The red and pink grids for the solid reference domain also remain fixed over time.
    Middle: $v_1$ and $\vs_1$ at $t=0.3$. Right: $p$ and $|\sigmasv|$ at $t=0.9$. Solutions computed on $\Gcrep^{(8)}$
    with $\scf=1$. 
  \label{fig:radialElasticPiston}
}
\end{figure}
}

Numerical solutions are computed using a composite grid, denoted by
$\Gcrep^{(j)}$, and shown in~Figure~\ref{fig:radialElasticPiston}.
The grid $\Gcrep^{(j)}$ has a target grid spacing of $\Delta s = 1/(10 j)$ for a chosen resolution parameter~$j$.
The solid disk is represented by an inner Cartesian (red) grid and
an interface-fitted annular (pink) grid.  These two grids are fixed in the
reference domian, ${\bar\Omega}_0$, of the solid.
The fluid domain is covered by an interface-fitted annular (green)
grid that deforms over time with the radial motion of the interface, and an outer annular (blue) grid
that is fixed to the outer boundary of the fluid.

{
\def\width{8}
\def\hscale{.75}
\newcommand{\figWidth}{\width cm}
\newcommand{\trimfig}[2]{\trimw{#1}{#2}{.0}{.05}{.0}{.0}}
\begin{figure}[h]
\begin{center}
\resizebox{14cm}{!}{
\begin{tikzpicture}[scale=1]
\useasboundingbox (-.1,.1) rectangle (2*\width+.1,-3*\width*\hscale+.75);  
\draw(0.0,0.0) node[anchor=north west,xshift=-15pt,yshift=8pt] {\trimfig{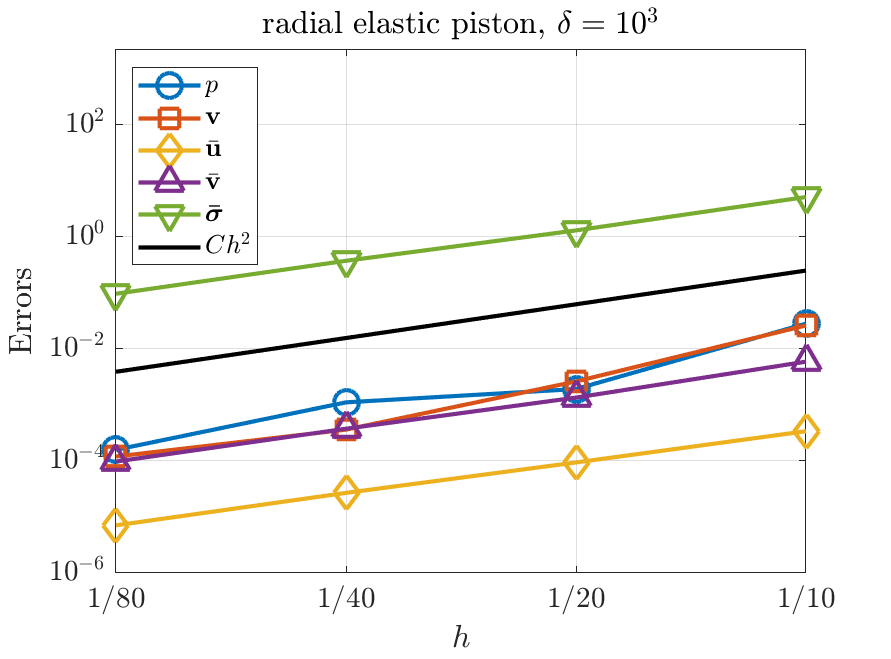}{\figWidth}};
\draw(\width,0.0) node[anchor=north west,xshift=-15pt,yshift=8pt] {\trimfig{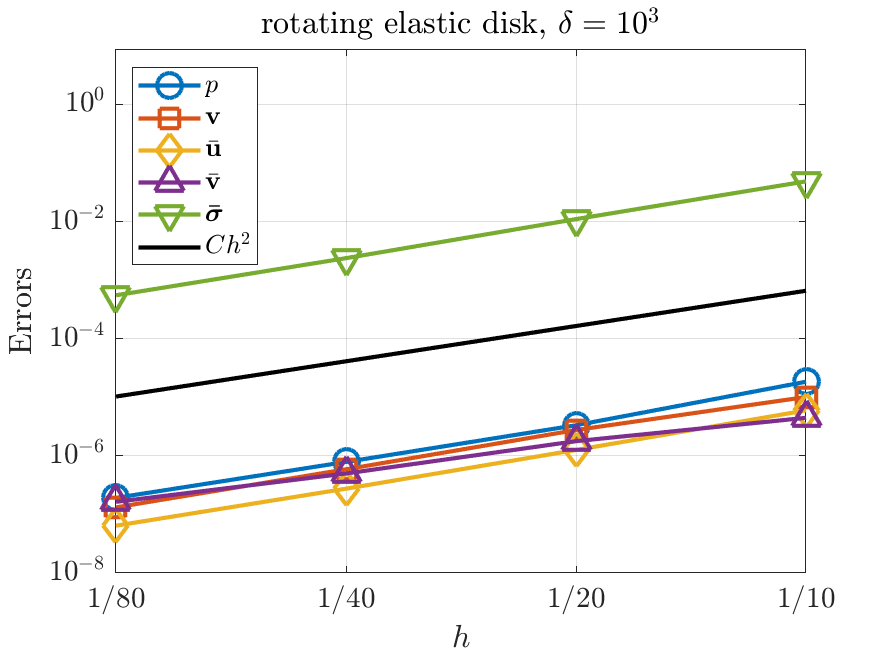}{\figWidth}};
\draw(0.0,-\hscale*\width) node[anchor=north west,xshift=-15pt,yshift=6pt] {\trimfig{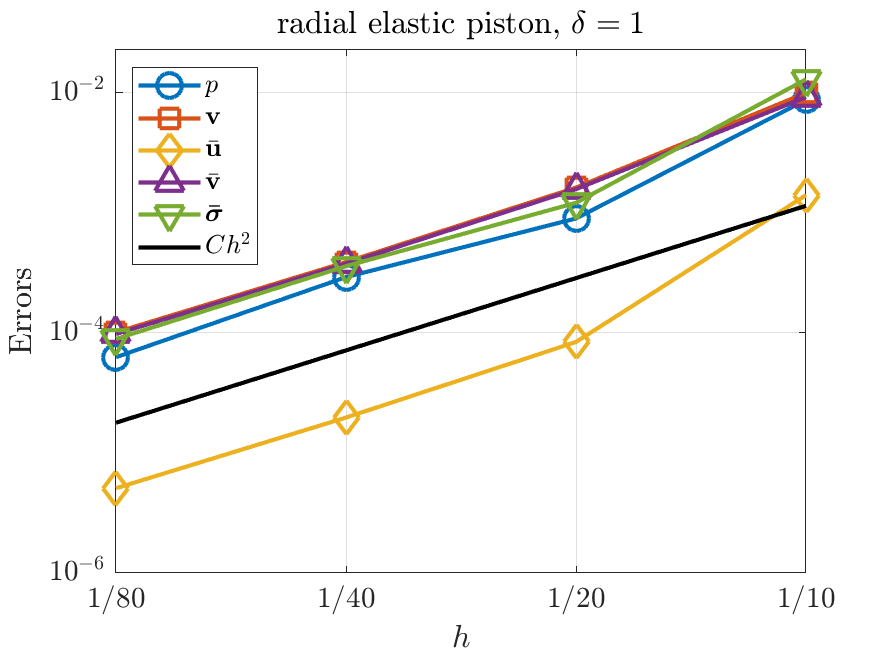}{\figWidth}};
\draw(\width,-\hscale*\width) node[anchor=north west,xshift=-15pt,yshift=6pt] {\trimfig{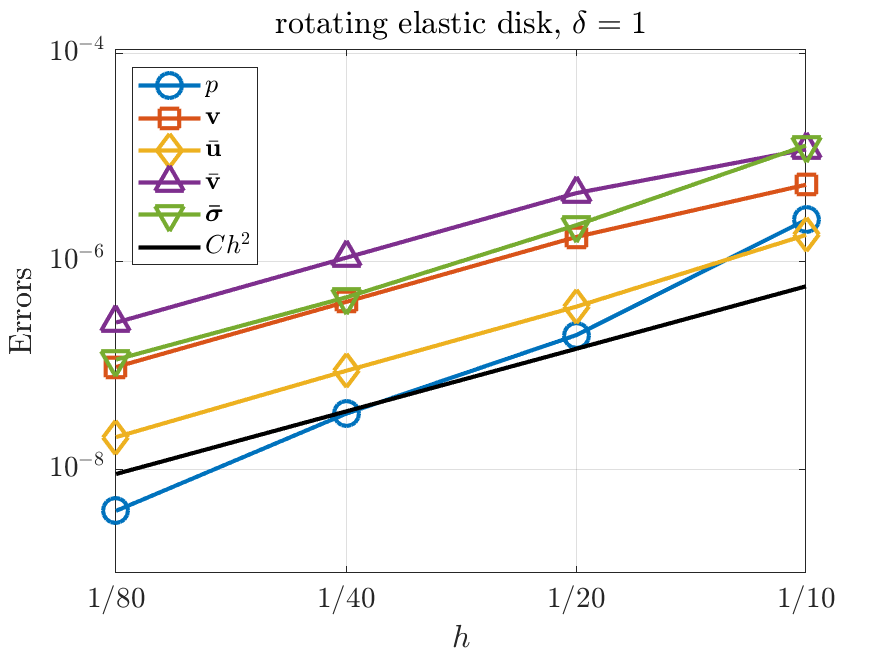}{\figWidth}};
\draw(0.0,-2*\hscale*\width) node[anchor=north west,xshift=-15pt,yshift=4pt] {\trimfig{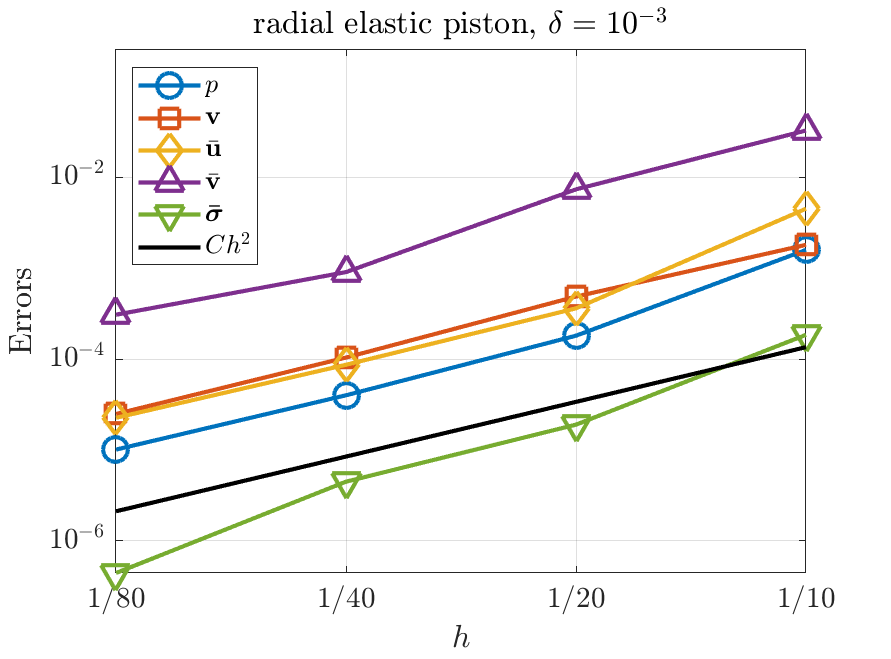}{\figWidth}};
\draw(\width,-2*\hscale*\width) node[anchor=north west,xshift=-15pt,yshift=4pt] {\trimfig{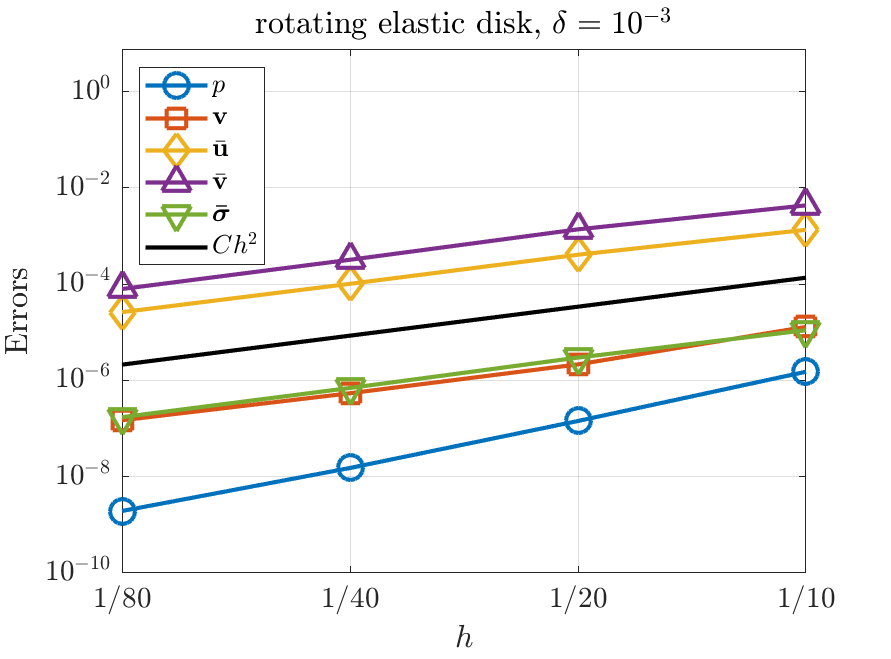}{\figWidth}};
\end{tikzpicture}
}
\end{center}
\caption{
Grid convergence studies of the max-norm errors at time $t_f=0.6$ 
for the circularly symmetric elastic piston with density ratios
$\delta=10^3$ (top left), $\delta=1$ (middle left), and $\delta=10^{-3}$ (bottom left) 
and the rotating elastic disk with density ratios
$\delta=10^3$ (top right), $\delta=1$ (middle right), and $\delta=10^{-3}$ (bottom right). 
\label{fig:radialElasticPistonConv}}
\end{figure}
}

Figure~\ref{fig:radialElasticPiston} also shows the behavior of the solution at times $t=0.3$ and $0.9$ computed using the AMP~algorithm and the grid $\Gcrep^{(8)}$.  The density ratio for this solution is $\scf=1$.  We note that while the viscous terms in the fluid momentum equation vanish identically in the derivation of the exact solution, the numerical scheme requires a choice for the visosity which is taken to be $\nu=0.05$.
The interface between the fluid and solid expands and contracts over time.
At $t=0.3$ the interface is near its largest radius, while at $t=0.9$ it is near its smallest.

Grid convergence studies of the max-norm errors are given in the left column of plots of Figure~\ref{fig:radialElasticPistonConv}
for density ratios of $\scf=10^3$, $\scf=1$ and $\scf=10^{-3}$.
The computations are stable for all values of $\delta$ considered, and we observe that the max-norm errors converge at close to second-order accuracy as indicated by the
(black) reference line in each plot.





\subsubsection{Rotating elastic disk in a fluid annulus} \label{sec:rotatingElasticDisk}

The second FSI benchmark problem has a geometry similar to that shown in Figure~\ref{fig:radialElasticPistonGeometry},
but instead of the radial motion shown we consider a circumferential motion of the disk. 
The exact solution, derived in~\ref{sec:rotatingElasticDiskExactSolution}, for the circumferential components of the fluid velocity and solid displacement takes the form 
\begin{equation}
v_\theta(r,t) = \hat{v}_\theta(r) e^{i \omega t}, \qquad \us_\theta(r,t) = \hat{\us}_\theta(r) e^{i \omega t},
\label{eq:radialTransSolns}
\end{equation}
where
\begin{equation}
\hat{v}_\theta(r)=b\bigl[J_1(\lambda r)Y_1(\lambda R)-J_1(\lambda R)Y_1(\lambda r)\bigr],\qquad \hat{\us}_\theta(r)=\bar bJ_1(k_sr).
\label{eq:radialTransCoefficients}
\end{equation}
Here,
\begin{align}
b = \frac{i \omega \us_0}{J_1(\lambda r_0) Y_1 (\lambda r_0) - J_1 (\lambda R) Y_1(\lambda r_0)}, \qquad
\lambda^2= {i \omega \over \nu}, \qquad
\bar{b} = \frac{\us_0}{J_1(k_s r_0)}, \qquad
k_s = {\omega \over \cs},
\end{align}
and $Y_1(z)$ is the Bessel function of the second kind of order one.  The parameter $\us_0$ determines the amplitude of the rotational
motion.  Solutions for the circumferential components of the fluid velocity and the 
solid displacement are taken from the real parts of $v_\theta(r,t)$ and $\us_\theta(r,t)$ defined in~\eqref{eq:radialTransSolns}.
The pressure is then determined by
\begin{align}
p(r,t) = \rho \int_{r_0}^r \frac{v_\theta(s,t)^2}{s} \, ds.
\end{align}
The frequency $\omega$ is obtained from the solution of the dispersion relation in~\eqref{eq:radialTransConstraint}.
For our computations, we choose $r_0= \frac{1}{2}$, $R=1$, $\rho=1$, $\nu=0.1$, $\rhos=\mus=\delta$ and $\us_0=10^{-5}$. 
Selected values of $\omega$ and the corresponding values for $b/\us_0$ and $\bar{b}/\us_0$ are given in Table~\ref{fig:ElasticShearFrequencies}
for different values of $\delta$. 
{
\begin{table}[h]\tableFont
\caption{Selected solutions for the frequency $\omega$ of the dispersion relation in~\eqref{eq:radialTransConstraint} for $\delta=10\sp{-3}$, $1$ and $10\sp{3}$.  The corresponding values for $b/\us_0$ and $\bar{b}/\us_0$ are also given.  Additional parameters chosen are $r_0= \frac{1}{2}$, $R=1$, $\rho=1$, $\nu=0.1$ and $\rhos=\mus=\delta$. 
\label{fig:ElasticShearFrequencies}}
\begin{center}
\begin{tabular}{cccc}
\hline
$\delta$ & $\omega$ & $b / \us_0$ & $ \bar{b} / \us_0$ \\
\hline
$10^{-3}$ & $7.664+0.001497i$ & $4.938-4.327i$ & $-1523+2316i$ \\
$1$ & $8.778+0.7854i$ & $4.355-6.155i$ & $-3.512+1.913i$ \\
$10^3$ & $10.27+0.002055i$ & $2.216-5.864i$ & $-2.944+0.0005914i$ \\
\hline
\end{tabular}
\end{center}
\end{table}
}

Numerical solutions are computed using the AMP~algorithm and the composite grid $\Gcrep^{(j)}$ as described for the previous problem.
Figure~\ref{fig:radialShearElasticPiston} shows the composite grid, and the computed solution at $t=0.1$ and $0.3$ for the case $\scf=1$.
The Lam\'e parameter $\lambdas$ does not appear in the exact solution for this problem, but a value is required for the numerical
scheme and we set $\lambdas=\delta$ in the
numerical calculations.  

Grid convergence studies of the max-norm errors in the right column of plots of 
Figure~\ref{fig:radialElasticPistonConv} for density ratios of $\scf=10^3$ , $\scf=1$ and 
$\scf=10^{-3}$.
The computations are stable, as before, and the max-norm errors all converge at close to second-order accuracy.

{
%
\newcommand{\drawContour}[9]{%
\begin{scope}[#1]
\draw(0.0,0) node[anchor=south west,xshift=-4pt,yshift=+0pt] {\trimfig{fig/#2}{\figWidth}};
  \draw(1.3,1.3) node[draw,fill=white,anchor=west,xshift=2pt,yshift=1pt] {\scriptsize #3};
  \draw(2.2,2.2) node[draw,fill=white,anchor=west,xshift=2pt,yshift=1pt] {\scriptsize #4};
  \draw(.5,5) node[draw,fill=white,anchor=west,xshift=2pt,yshift=-3pt] {\scriptsize #5}; 
\begin{scope}[xshift=5pt,yshift=-7pt]
  \draw (\xcb,\ycb) node[anchor=south west,xshift=0.25cm,yshift=.5cm,rotate=-90] {\trimfigcb{fig/colourBarLines}{\cbWidth}{\cbHeight}};
  \draw (.8,0) node[anchor=north,xshift=+5pt,yshift=+2pt] {\scriptsize $#6$};
  \draw (2.8,0) node[anchor=north,xshift=+0pt,yshift=+2pt] {\scriptsize fluid};
  \draw (4.8,0) node[anchor=north,xshift=-2pt,yshift=+2pt] {\scriptsize $#7$};
  \begin{scope}[yshift=14pt]
    \draw (.8,0) node[anchor=north,xshift=+5pt,yshift=+2pt] {\scriptsize $#8$};
    \draw (2.8,0) node[anchor=north,xshift=+0pt,yshift=+2pt] {\scriptsize solid};
    \draw (4.8,0) node[anchor=north,xshift=-2pt,yshift=+2pt] {\scriptsize $#9$};
 \end{scope}
\end{scope}
\end{scope}
}
\newcommand{\cbWidth}{.2cm}
\newcommand{\cbHeight}{4cm}
\newcommand{\xcb}{.5cm}
\newcommand{\ycb}{-.2cm}
\setlength{\ycbTop}{\ycb+\cbHeight}
\setlength{\ycbMid}{\ycb+\cbHeight*\real{.5}}
\newcommand{\trimfigcb}[3]{\includegraphics[width=#2, height=#3, clip, trim=17cm 2.35cm 1.65cm 2.35cm]{#1}}
\newcommand{\figWidth}{5cm}
\newcommand{\trimfig}[2]{\trimh{#1}{#2}{.11}{.11}{.11}{.11}}
\begin{figure}[htb]
\begin{center}
\begin{tikzpicture}[scale=1]
  \useasboundingbox (0,0.25) rectangle (15.5,5.25);  

   \draw(0.0,0.0) node[anchor=south west,xshift=-15pt,yshift=0pt] {\trimfig{fig/radialElasticPistonGrid}{\figWidth}};
   \drawContour{xshift= 5.cm,yshift=0cm}{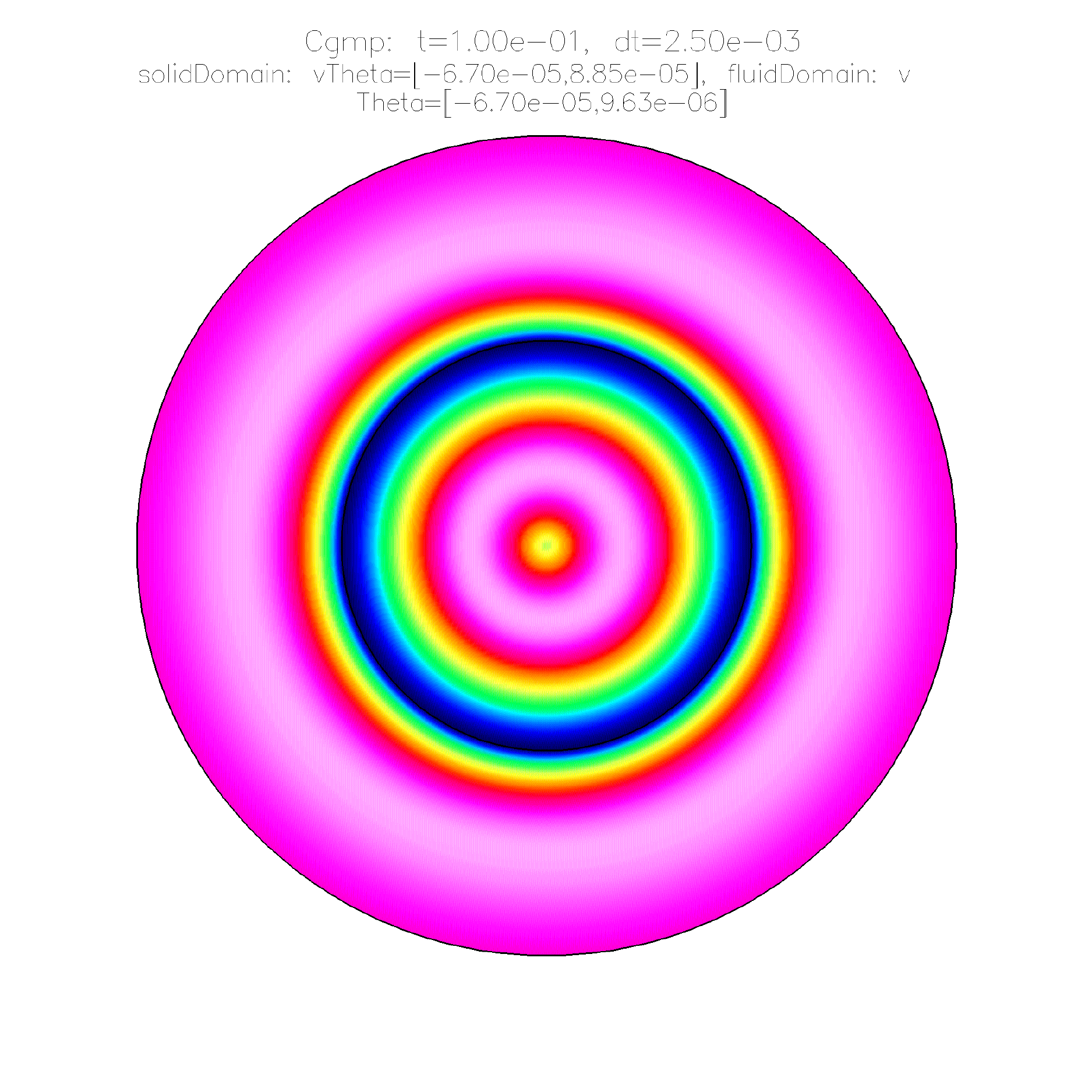}{$v_\theta$}{$\vs_\theta$}{$t=0.1$}{$-6.70e-05$}{$9.63e-05$}{$-6.70e-05$}{$8.85e-05$};
   \drawContour{xshift=10.cm,yshift=0cm}{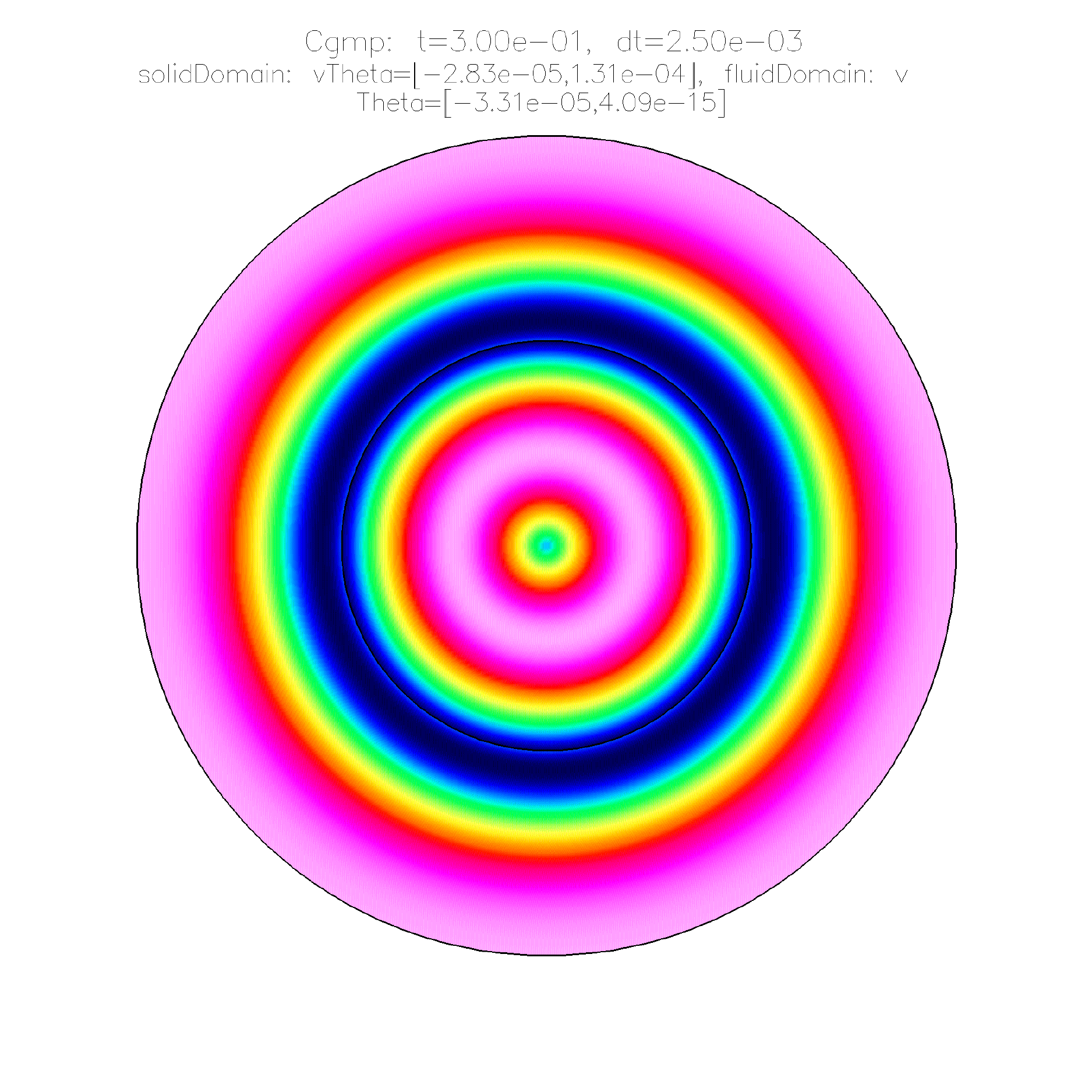}{$v_\theta$}{$\vs_\theta$}{$t=0.3$}{$-3.31e-05$}{$0$}{$-2.83e-05$}{$1.31e-04$};
\end{tikzpicture}
\end{center}
  \caption{Left: composite grid $\Gcrep^{(2)}$ for the radial elastic piston problem.
    The green fluid grid adjacent to the interface  moves over time. The blue fluid background grid remains fixed.
    The red and pink grids for the solid reference domain also remain fixed over time.
    Middle: $v_\theta$ and $\vs_\theta$ at $t=0.1$. Right: $v_\theta$ and $\vs_\theta$ at $t=0.3$. Solutions computed on $\Gcrep^{(8)}$
    with $\scf=1$.  
}
  \label{fig:radialShearElasticPiston}
\end{figure}
}


%
%








\subsection{Radial traveling wave solution} \label{sec:travelingWave}

We next consider traveling wave solutions involving both radial and circumferential motion of an incompressible fluid in a circular region $0 < r < r_0$ coupled to a surrounding linearly elastic solid in the annular region $r_0 < r < R$.  The geometry of the FSI problem is shown in Figure~\ref{fig:radialTravelingWaveGrid}.  An exact solution is derived for this problem in~\ref{sec:radialTravelingWaveSolution} assuming a Stokes fluid and an linearized interface about $r=r_0$.  In this solution, the fluid velocity and pressure take the form 
\begin{align*}
v_r(r,\theta,t) = \hat{v}_r(r) e^{i (n \theta - \omega t)}, \qquad
v_\theta(r,\theta,t) = \hat{v}_\theta(r) e^{i (n \theta - \omega t)} , \qquad
p(r,\theta,t) = \hat{p}(r) e^{i (n \theta - \omega t)},
\end{align*}
where $\omega$ is a frequency, $n$ is a positive integer, and $\hat{v}_r(r)$, $\hat{v}_\theta(r)$ and $\hat{p}(r)$
are coefficient functions.  Similarly, the displacement in the solid is given by
\begin{align*}
\us_r(r,\theta,t) = \hat{\us}_r(r) e^{i (n \theta - \omega t)}, \qquad
\us_\theta(r,\theta,t) = \hat{\us}_\theta(r) e^{i (n \theta - \omega t)},
\end{align*}
with coefficient functions $\hat{\us}_r(r)$ and $\hat{\us}_\theta(r)$.  Following the analysis in~\ref{sec:radialTravelingWaveSolution}, the coefficient functions are found to be
\begin{align*}
\hat{v}_r(r) = d \frac{J_n(\lambda r)}{r}
+ p_I \frac{n r^{n-1}}{\mu \lambda^2 r_0^n},
\qquad
\hat{v}_{\theta}(r) =
d \frac{i\lambda}{n} J_n'(\lambda r)
+ p_I \frac{i n r^{n-1}}{\mu \lambda^2 r_0^n},
\qquad
\hat{p}(r) = p_I \left(\frac{r}{r_0} \right)^n ,
\end{align*}
and
\newcommand{\db}{\bar{d}}
\begin{align*}
\hat{\us}_r(r) &= 
k_p \left(\db_1 J_n'(k_p r) + \db_2 Y_n'(k_p r) \right)
+ \frac{i n}{r} \left( \db_3 J_n(k_s r) + \db_4 Y_n(k_s r)\right), \\
\hat{\us}_\theta(r) &= 
-k_s \left(\db_3 J_n'(k_s r) + \db_4 Y_n'(k_s r) \right)
+ \frac{i n}{r} \left( \db_1 J_n(k_p r) + \db_2 Y_n(k_p r) \right),
\end{align*}
where $\lambda^2= i \omega / \nu$, $k_s = \omega / \cs$, $k_p = \omega / \cp$, $(p_I,d,\db_1,\db_2,\db_3,\db_4)$ are constants of integration, and $J_n$ and $Y_n$ are Bessel functions of order $n$.  An application of traction-free boundary conditions at $r=R$ and matching conditions involving velocity and stress at the interface $r=r_0$ leads to a homogeneous linear system of equations for the six constants.  Nontrivial solutions are found if the system is singular, which leads to the dispersion relation in~\eqref{eq:dispersion} involving the frequency $\omega$ and circumferential wave number~$n$.

\newcommand{\Gcrtw}{\Gc_{\rm rtw}}
Numerical solutions are computed for the case $r_0=1$ and $R=1.2$.  The fluid has $\rho=1$ and $\nu=0.1$, while the solid has $\mus= \lambdas= \rhos=\delta$.  The domain for the FSI problem is represented by the composite grid, $\Gcrtw^{(j)}$, with target grid spacing  $\Delta s = 1/(10j)$ as shown in Figure~\ref{fig:radialTravelingWaveGrid}.  The circular fluid domain is covered by a fixed background Cartesian grid (blue) 
and interface-fitted grid (green). The solid domain is represented by an annular grid (red).  While the exact solution assumes a linearization about a fixed interface position, $r=r_0$, the numerical solution does not make this assumption so that the interface-fitted grid of the fluid moves in time with the (very small) deformation of the interface.  Also, the numerical solution does not assume a Stokes fluid.  However, since the velocity of the fluid is very small for this problem, the nonlinear convective terms in the momentum equation are negligible.  Initial conditions for the AMP~time-stepping scheme are taken from the exact (linearized) solution at $t=0$ for an $(n,\omega)$-pair satisfying the dispersion relation in~\eqref{eq:dispersion} and the corresponding set of constants $(p_I,d,\db_1,\db_2,\db_3,\db_4)$ with a normalization taken to be ${\bar u}_0=10\sp{-7}$.  The exact solution is described in terms of complex-valued functions, but we use the real part to set the initial conditions and to compare with the numerical solution.

Figure~\ref{fig:radialTravelingWaveGrid} shows representative solutions at $t=0.1$ and~$0.3$ computed using the grid $\Gcrtw^{(8)}$.  The density ratio is chosen to be $\delta=1$, and a $(n,\omega)$-pair satisfying the dispersion relation is given by $n=3$ and $\omega=3.491-1.154\,i$.  The shaded contours of the fluid pressure and the norm of the solid displacement both show three-fold symmetry in the circumferential direction in agreement with the choice of $n=3$.
Grid convergence studies of the max-norm errors are given in 
Figure~\ref{fig:radialTravelingWaveConv} for density ratios of $\scf=10$, $\scf=1$ and $\scf=10^{-1}$.
The computations are stable and the max-norm errors all converge at close to second-order accuracy.

{
%
\newcommand{\drawContour}[9]{%
\begin{scope}[#1]
\draw(0.0,0) node[anchor=south west,xshift=-4pt,yshift=+0pt] {\trimfig{fig/#2}{\figWidth}};
  \draw(1.1,1.1) node[draw,fill=white,anchor=west,xshift=2pt,yshift=1pt] {\scriptsize #3};
  \draw(2.2,2.2) node[draw,fill=white,anchor=west,xshift=2pt,yshift=1pt] {\scriptsize #4};
  \draw(.5,5) node[draw,fill=white,anchor=west,xshift=2pt,yshift=-3pt] {\scriptsize #5}; 
\begin{scope}[xshift=5pt,yshift=-7pt]
  \draw (\xcb,\ycb) node[anchor=south west,xshift=0.25cm,yshift=.5cm,rotate=-90] {\trimfigcb{fig/colourBarLines}{\cbWidth}{\cbHeight}};
  \draw (.8,0) node[anchor=north,xshift=+5pt,yshift=+2pt] {\scriptsize $#6$};
  \draw (2.8,0) node[anchor=north,xshift=+0pt,yshift=+2pt] {\scriptsize fluid};
  \draw (4.8,0) node[anchor=north,xshift=-2pt,yshift=+2pt] {\scriptsize $#7$};
  \begin{scope}[yshift=14pt]
    \draw (.8,0) node[anchor=north,xshift=+5pt,yshift=+2pt] {\scriptsize $#8$};
    \draw (2.8,0) node[anchor=north,xshift=+0pt,yshift=+2pt] {\scriptsize solid};
    \draw (4.8,0) node[anchor=north,xshift=-2pt,yshift=+2pt] {\scriptsize $#9$};
 \end{scope}
\end{scope}
\end{scope}
}
\newcommand{\cbWidth}{.2cm}
\newcommand{\cbHeight}{4cm}
\newcommand{\xcb}{.5cm}
\newcommand{\ycb}{-.2cm}
\setlength{\ycbTop}{\ycb+\cbHeight}
\setlength{\ycbMid}{\ycb+\cbHeight*\real{.5}}
\newcommand{\trimfigcb}[3]{\includegraphics[width=#2, height=#3, clip, trim=17cm 2.35cm 1.65cm 2.35cm]{#1}}
\newcommand{\figWidth}{5cm}
\newcommand{\trimfig}[2]{\trimh{#1}{#2}{.11}{.11}{.11}{.11}}
\begin{figure}[htb]
\begin{center}
\begin{tikzpicture}[scale=1]
  \useasboundingbox (0,0.25) rectangle (15.5,5.25);  

   \draw(0.0,0.0) node[anchor=south west,xshift=-15pt,yshift=0pt] {\trimfig{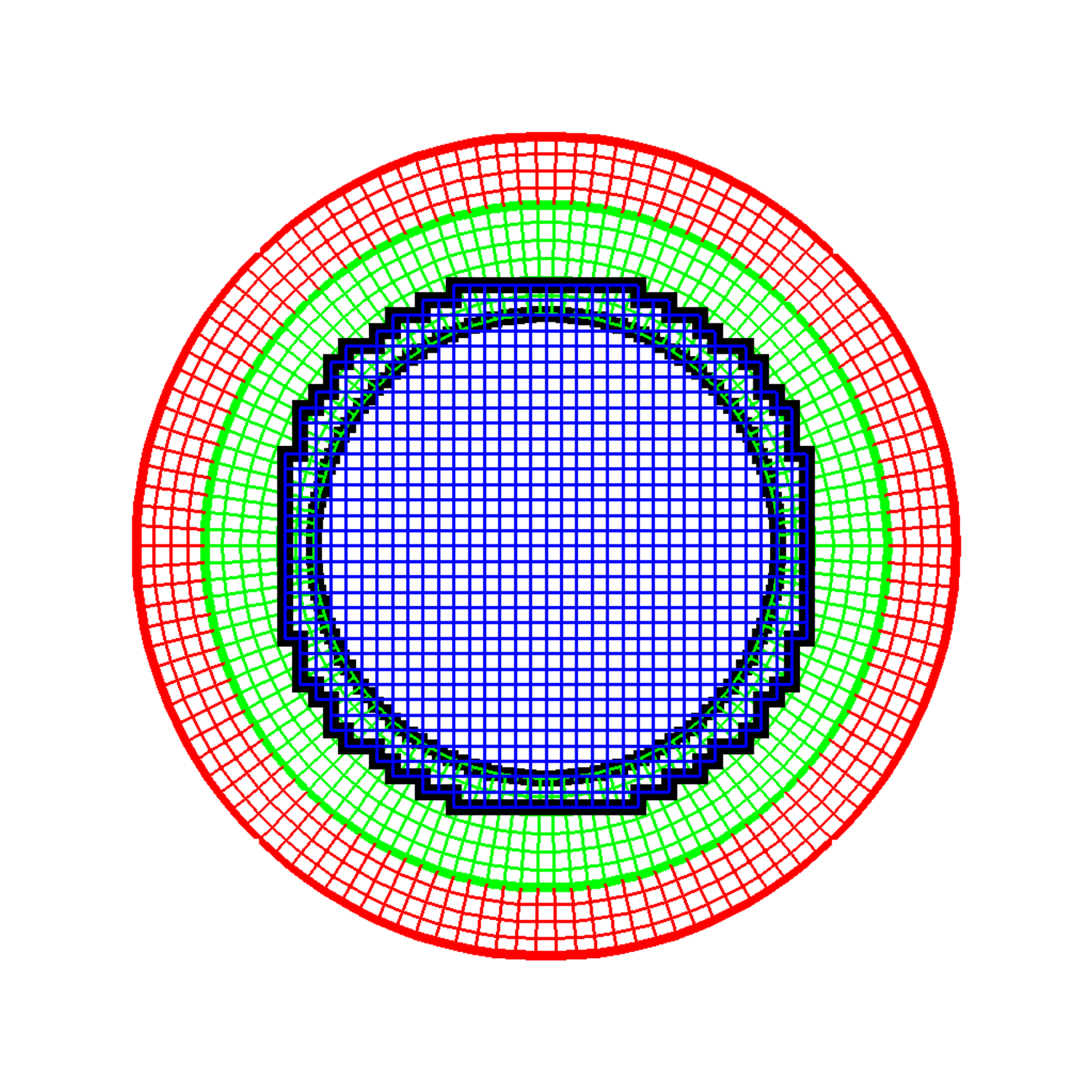}{\figWidth}};
   \drawContour{xshift= 5.cm,yshift=0cm}{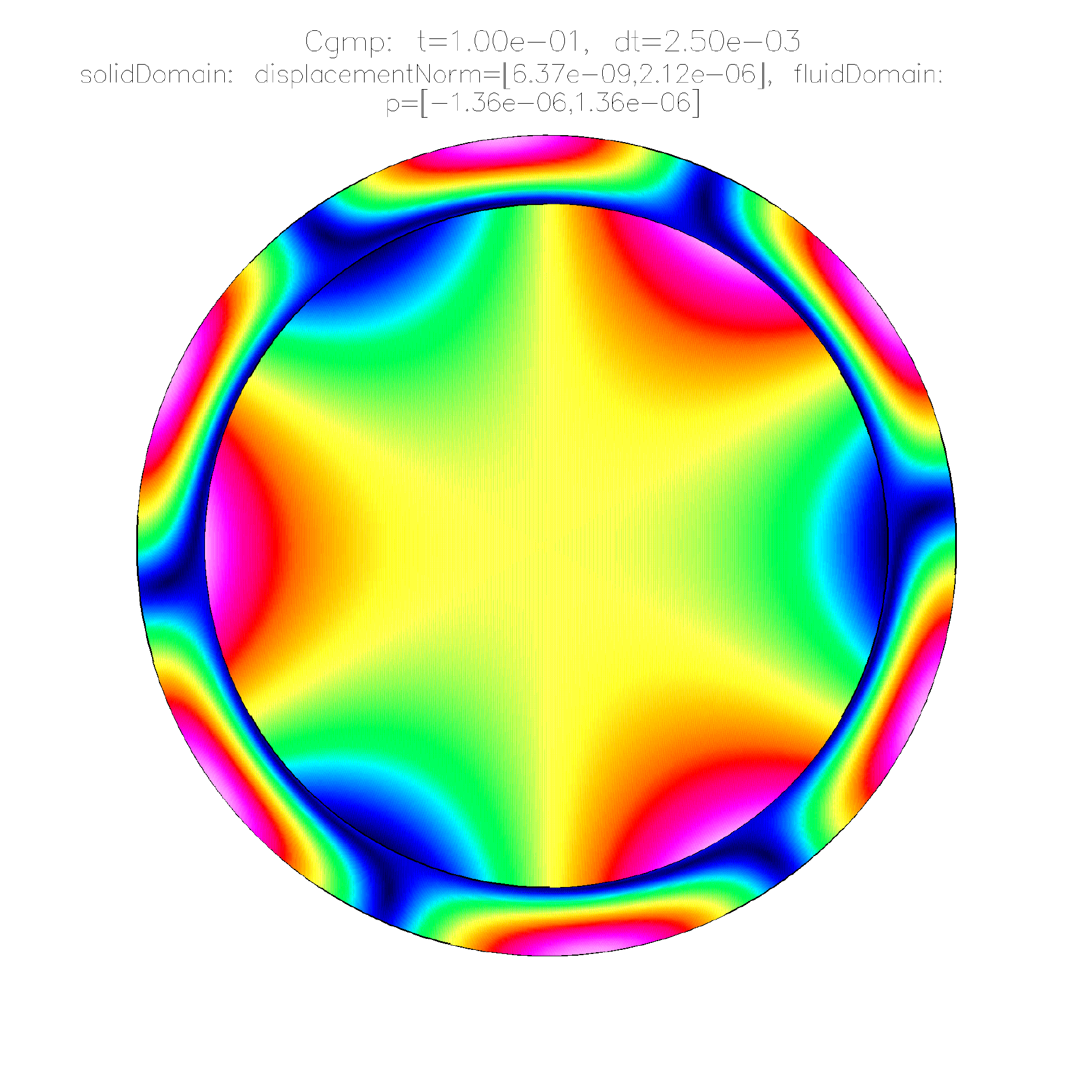}{$|\usv|$}{$p$}{$t=0.1$}{$-1.36e-06$}{$1.36e-06$}{$0$}{$2.12e-06$};
   \drawContour{xshift=10.cm,yshift=0cm}{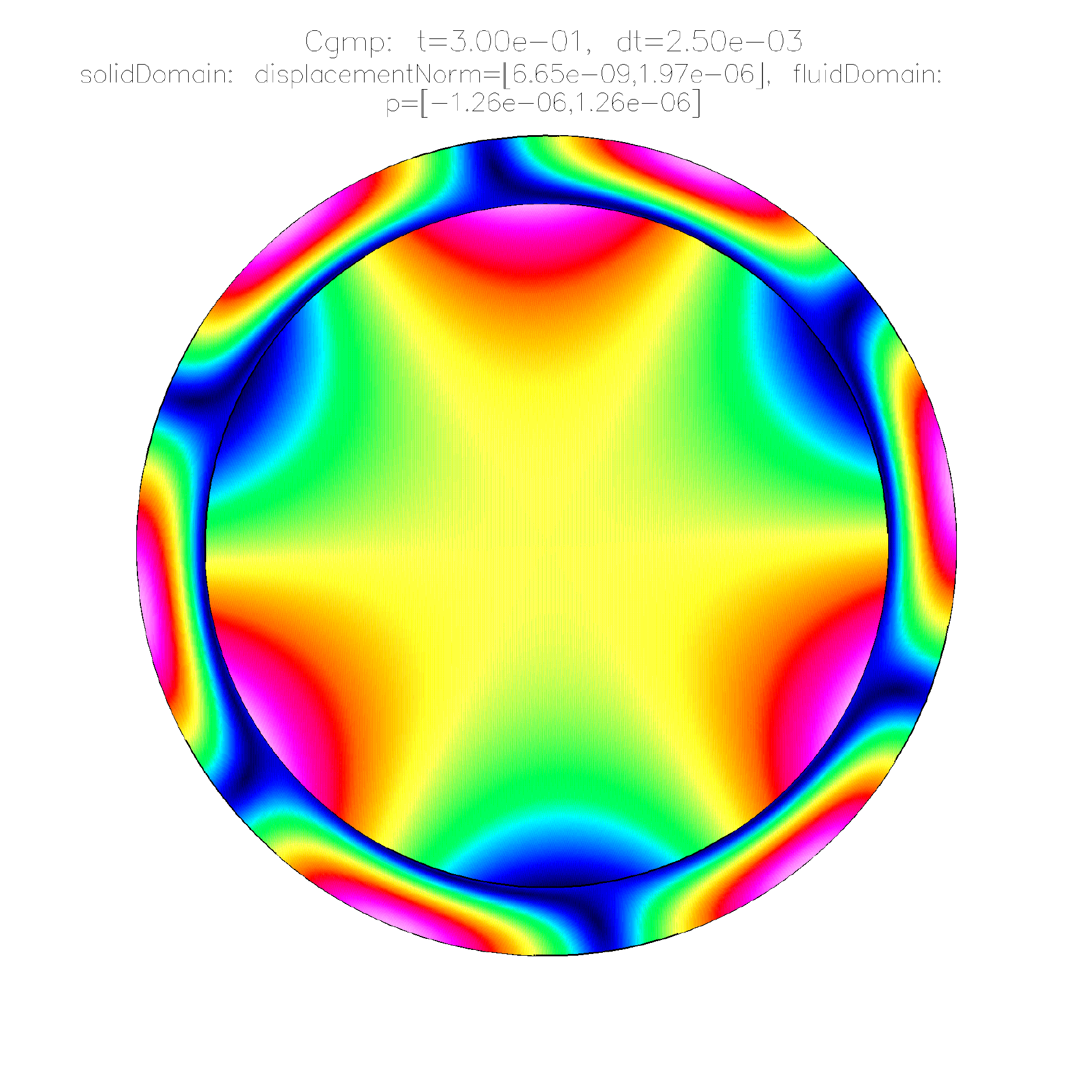}{$|\usv|$}{$p$}{$t=0.3$}{$-1.26e-06$}{$1.26e-06$}{$0$}{$1.96e-06$};
\end{tikzpicture}
\end{center}
  \caption{Left: composite grid $\Gcrtw^{(2)}$ for the radial traveling wave.
    The green fluid grid adjacent to the interface moves over time while the blue fluid background grid remains fixed.
    The red grid for the solid reference domain also remains fixed over time.
    Middle and right: fluid pressure and $\| \usv\|$ at times $t=.1$ and $t=.3$, computed
      on grid $\Gcrtw^{(8)}$ with $\scf=1.$ 
}
  \label{fig:radialTravelingWaveGrid}
\end{figure}
}



{
\def\width{7}
\def\hscale{.75}
\newcommand{\figWidth}{\width cm}
\newcommand{\trimfig}[2]{\trimw{#1}{#2}{.0}{.05}{.0}{.0}}
\begin{figure}[h]
\begin{center}
\begin{tikzpicture}[scale=1]
\useasboundingbox (0,0.25) rectangle (16.25,-15.2);  

\draw(\width,0.0) node[anchor=north west,xshift=-15pt,yshift=8pt] {\trimfig{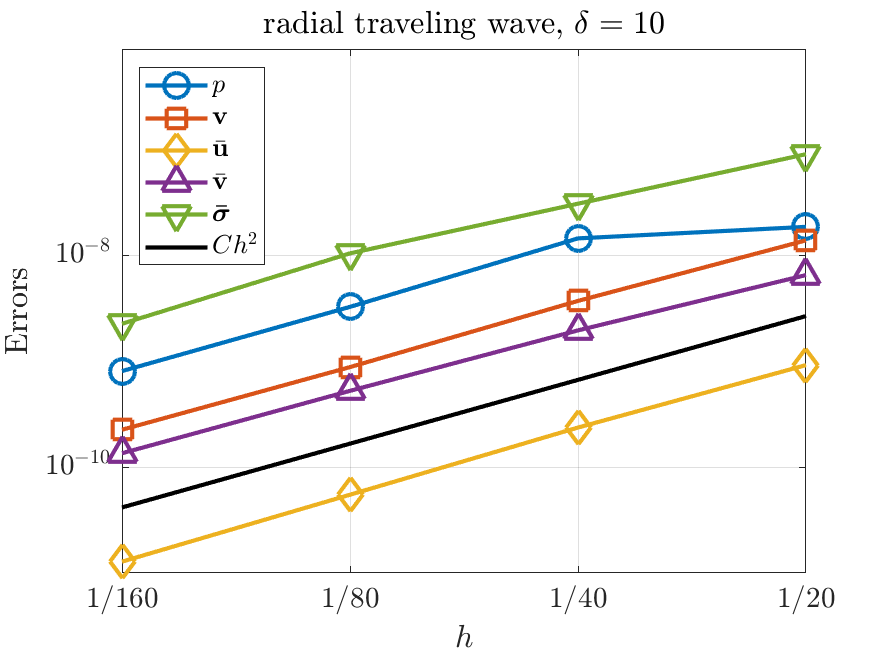}{\figWidth}};
\draw(\width,-\hscale*\width) node[anchor=north west,xshift=-15pt,yshift=6pt] {\trimfig{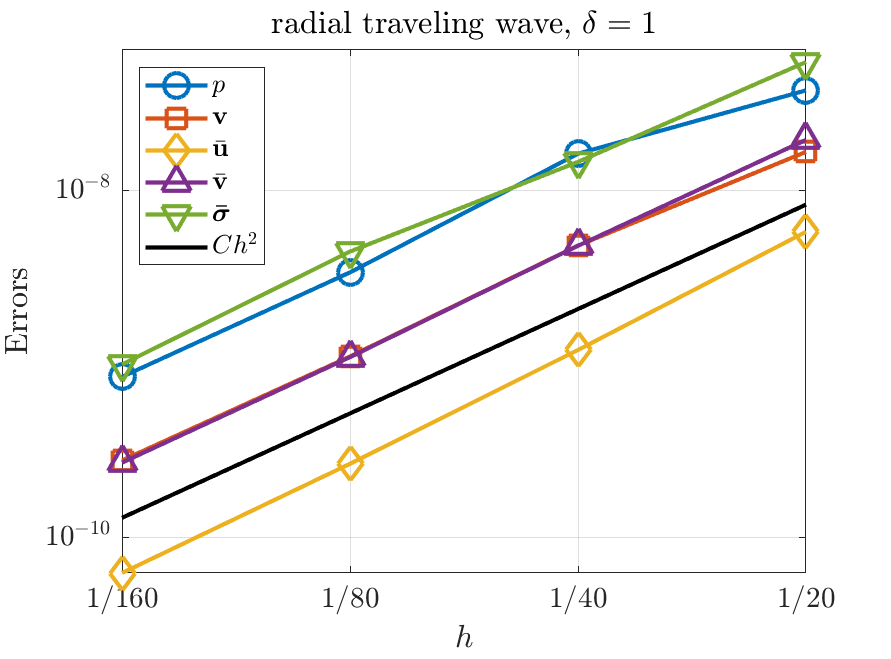}{\figWidth}};
\draw(\width,-2*\hscale*\width) node[anchor=north west,xshift=-15pt,yshift=4pt] {\trimfig{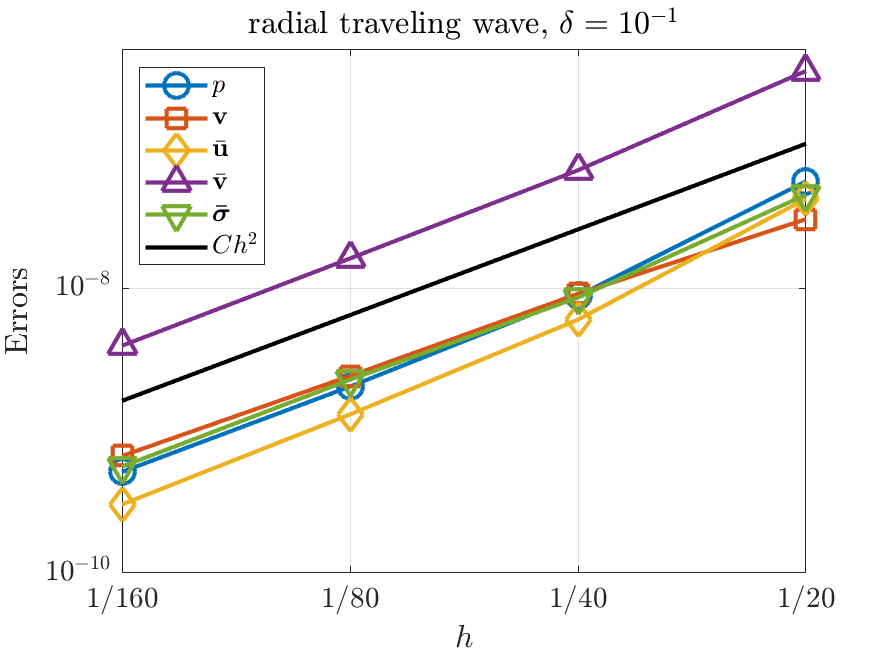}{\figWidth}};

\draw(.5*\width,-2.5*\hscale*\width) node[] 
{
  \footnotesize
  \begin{tabular}{c|ccc}
    \hline
$\delta$  & $10^{-1}$  \\
    \hline
 $\omega$  &  $4.487-0.1529i$  \\
 $d / \us_0$  &  $-0.9447+0.1398i$  \\
 $ p_I / \us_0$  &  $12.41+13.04i$  \\
 $\bar{d}_1 / \us_0$  &  $-1.966-0.1994i$  \\
 $\bar{d}_2 / \us_0$  &  $0.3716+0.03041i$  \\
 $\bar{d}_3 / \us_0$  &  $-0.01249-1.180i$  \\
 $\bar{d}_4 / \us_0$ &  $-0.09547-0.1188i$ \\
\hline
  \end{tabular}
};
\draw(.5*\width,-1.5*\hscale*\width) node[]
{
  \footnotesize
  \begin{tabular}{c|ccc}
    \hline
$\delta$  & $1$  \\
    \hline
 $\omega$  &  $3.491-1.154i$  \\
 $d / \us_0$  &  $-7.889+0.7399i$  \\
 $ p_I / \us_0$  &  $28.98+23.33i$  \\
 $\bar{d}_1 / \us_0$  &  $-5.778-6.619i$  \\
 $\bar{d}_2 / \us_0$  &  $2.166+0.6514i$  \\
 $\bar{d}_3 / \us_0$  &  $-0.4929+4.072i$  \\
 $\bar{d}_4 / \us_0$ &  $5.728+2.071i$ \\
\hline
  \end{tabular}
};

\draw(.5*\width,-0.5*\hscale*\width) node[]
{
  \footnotesize
  \begin{tabular}{c|ccc}
    \hline
$\delta$  & $10$  \\
    \hline
 $\omega$  &  $8.446-0.3587i$  \\
 $d / \us_0$  &  $0.08293+0.09885i$  \\
 $ p_I / \us_0$  &  $-9.442-10.40i$  \\
 $\bar{d}_1 / \us_0$  &  $0.4973-0.4251i$  \\
 $\bar{d}_2 / \us_0$  &  $-3.942-2.816i$  \\
 $\bar{d}_3 / \us_0$  &  $-5.889+6.165i$  \\
 $\bar{d}_4 / \us_0$ &  $0.3925+5.918i$ \\
\hline
  \end{tabular}
};

\end{tikzpicture}
\end{center}
\caption{Grid convergence studies of the max-norm errors in the numerical solution of the radial traveling wave problem at time $t=0.3$
for density ratios
$\delta=10$ (top right), $\delta=1$ (middle right), and $\delta=10^{-1}$ (bottom right).
The frequencies and constants for each case are tabulated to the left of each convergence plot.
For each case, $r_0=1$, $R=1.2$, $\rho=1$, $\nu=0.1$, $\mus= \lambdas= \rhos=\delta$ and $n=3$.
\label{fig:radialTravelingWaveConv}}
\end{figure}
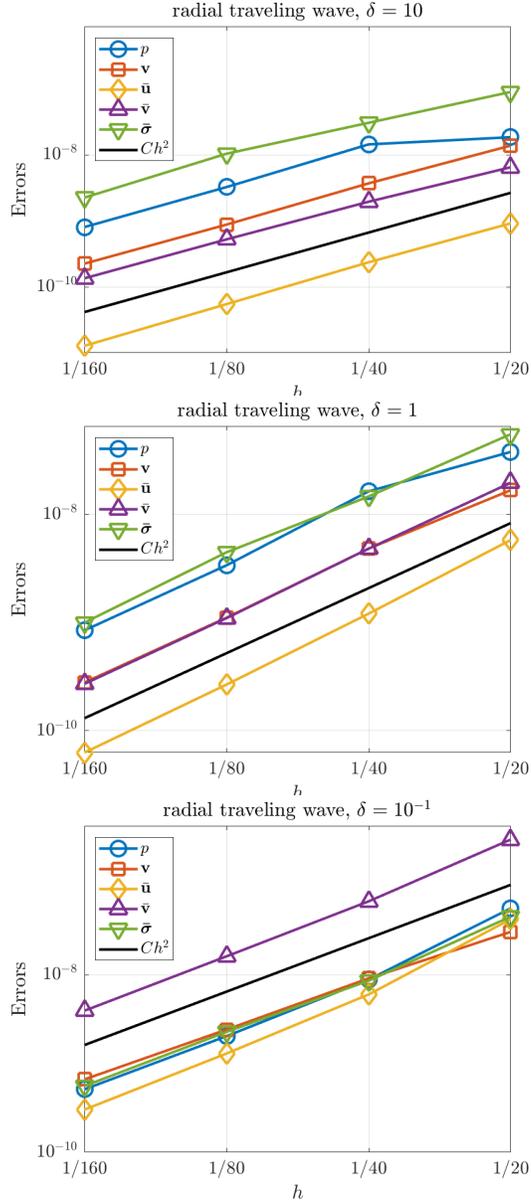
}


\subsection{Elastic annulus in a pressure-driven channel flow} \label{sec:elasticDiskInAChannel}

\newcommand{\ramp}{{\eta}}
In this section, we consider the flow of an incompressible fluid in a channel past an
embedded annular elastic solid which is anchored by a fixed inner boundary.
The fluid with density $\rho=1$ and viscosity $\nu=0.05$ exists in the rectangular region $[-3,3]\times[-1.75,1.75]$.
The solid annulus has center $\xv=(0,0)$ with outer radius $R_d=0.8$ and inner radius $R_i=0.4$.
The solid parameters are
$\rhos=\lambdas=\mus=\delta$, where the density ratio $\delta$ is a parameter
defining light and heavy solids as before.
Initially the fluid is at rest and the displacement of the solid is zero.  A flow in the channel is
driven by a pressure gradient determined by a variable pressure assigned at the inflow boundary on the left.
To define a smooth solution for testing convergence,
the pressure at the left inflow boundary is smoothly ramped from zero to a final value 
of $p_0=1$ at $t=1$.  Specifically, we set $p=p_0\ramp(t)$ at $x=-3$, where $\ramp(t)$ is a ramp function given by
\ba
    \ramp(t) = \begin{cases}
            (35+(-84+(70- 20 t) t) t) t^4,\quad     & \text{for $0\le t \le 1$}, \\
                1,   & \text{for $t>1$}.
               \end{cases}
              \label{eq:ramp}
\ea
The ramp function satisfies $\ramp=\ramp\sp\prime=\ramp\sp{\prime\prime}=\ramp\sp{\prime\prime\prime}=0$ at $t=0$, and it has three continuous derivatives at $t=1$.  We also set the tangential component of the velocity to be zero at the left inflow boundary.
No-slip boundary conditions are applied at the top and bottom walls of the channel, while the pressure and normal derivatives of the
velocity are set to zero
at the right outflow boundary.
The usual matching conditions involving velocity and stress are taken at the fluid-solid interface, and the inner solid boundary is taken as a zero displacement condition.  While an exact solution is not available for this problem, the set-up is very clean so that it makes a good benchmark problem for testing various FSI algorithms as is done below.

{
%
\newcommand{\drawContour}[9]{%
\begin{scope}[#1]
\draw(0.0,0) node[anchor=south west,xshift=-4pt,yshift=+0pt] {\trimfig{fig/#2}{\figWidth}};
  \draw(.5,.5) node[draw,fill=white,anchor=west,xshift=5pt,yshift=1pt] {\scriptsize #3};
  \draw(3.5,2.) node[draw,fill=white,anchor=west,xshift=2pt,yshift=1pt] {\scriptsize #4};
  \draw(.5,5) node[draw,fill=white,anchor=west,xshift=4pt,yshift=-5pt] {\scriptsize #5}; 
\begin{scope}[xshift=2cm,yshift=-9pt]
  \draw (\xcb,\ycb) node[anchor=south west,xshift=0.25cm,yshift=.5cm,rotate=-90] {\trimfigcb{fig/colourBarLines}{\cbWidth}{\cbHeight}};
  \draw (.8,0) node[anchor=north,xshift=+3pt,yshift=+2pt] {\scriptsize $#6$};
  \draw (2.8,0) node[anchor=north,xshift=+0pt,yshift=+2pt] {\scriptsize fluid};
  \draw (4.8,0) node[anchor=north,xshift=+0pt,yshift=+2pt] {\scriptsize $#7$};
  \begin{scope}[yshift=14pt]
    \draw (.8,0) node[anchor=north,xshift=+3pt,yshift=+2pt] {\scriptsize $#8$};
    \draw (2.8,0) node[anchor=north,xshift=+0pt,yshift=+2pt] {\scriptsize solid};
    \draw (4.8,0) node[anchor=north,xshift=+0pt,yshift=+2pt] {\scriptsize $#9$};
 \end{scope}
\end{scope}
\end{scope}
}
\newcommand{\cbWidth}{.2cm}
\newcommand{\cbHeight}{4cm}
\newcommand{\xcb}{.5cm}
\newcommand{\ycb}{-.2cm}
\setlength{\ycbTop}{\ycb+\cbHeight}
\setlength{\ycbMid}{\ycb+\cbHeight*\real{.5}}
\newcommand{\trimfigcb}[3]{\includegraphics[width=#2, height=#3, clip, trim=17cm 2.35cm 1.65cm 2.35cm]{#1}}
\newcommand{\figWidtha}{5cm}
\newcommand{\trimfiga}[2]{\trimhb{#1}{#2}{.11}{.11}{.11}{.11}}
\newcommand{\figWidth}{5cm}
\newcommand{\trimfig}[2]{\trimh{#1}{#2}{.11}{.11}{.275}{.275}}
\begin{figure}[htb]
\begin{center}
\begin{tikzpicture}[scale=1]
  \useasboundingbox (0,0.1) rectangle (14,5.25);  

   \draw(0.0,0.0) node[anchor=south west,xshift=-15pt,yshift=0pt] {\trimfiga{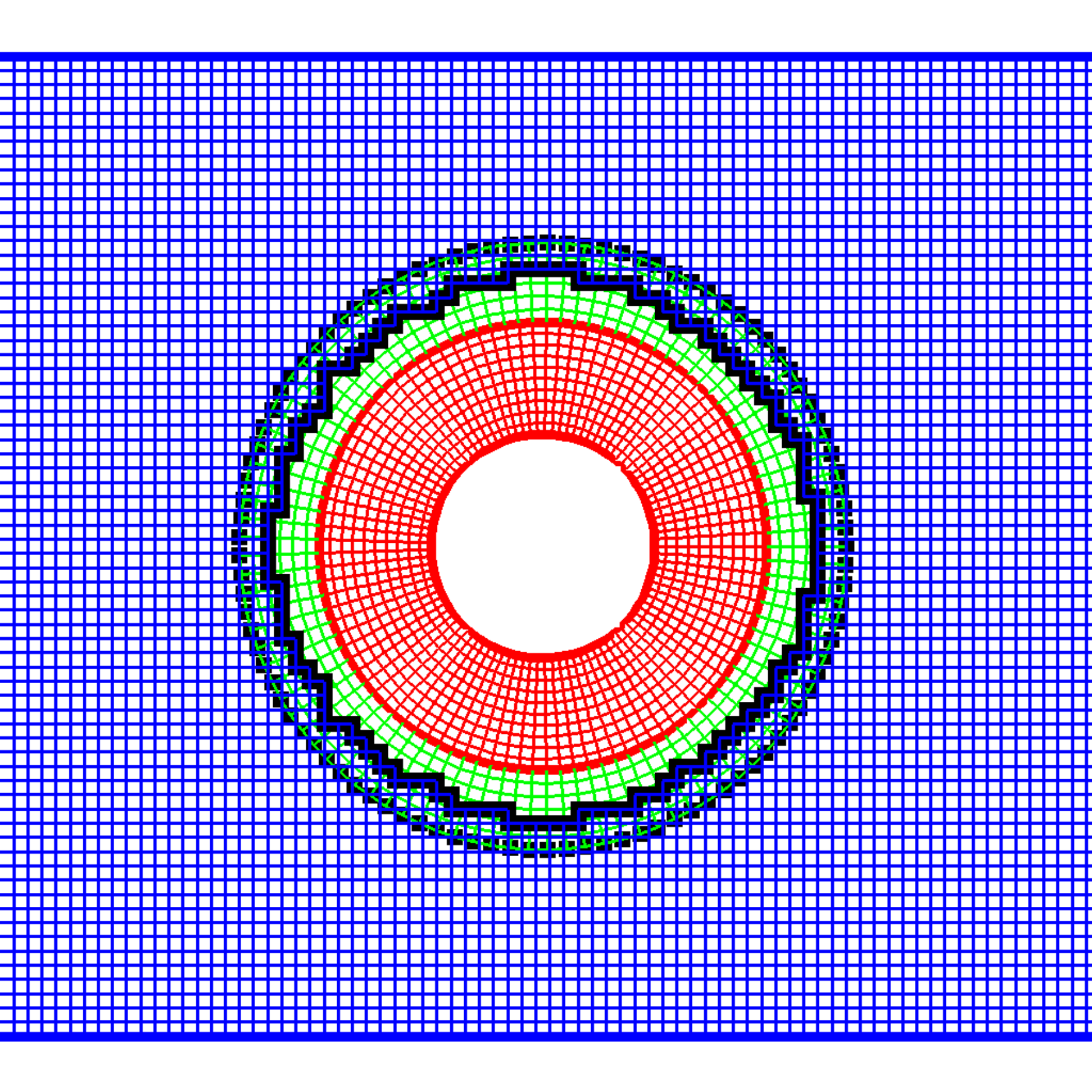}{\figWidtha}};
   \drawContour{xshift= 5.cm,yshift=0cm}{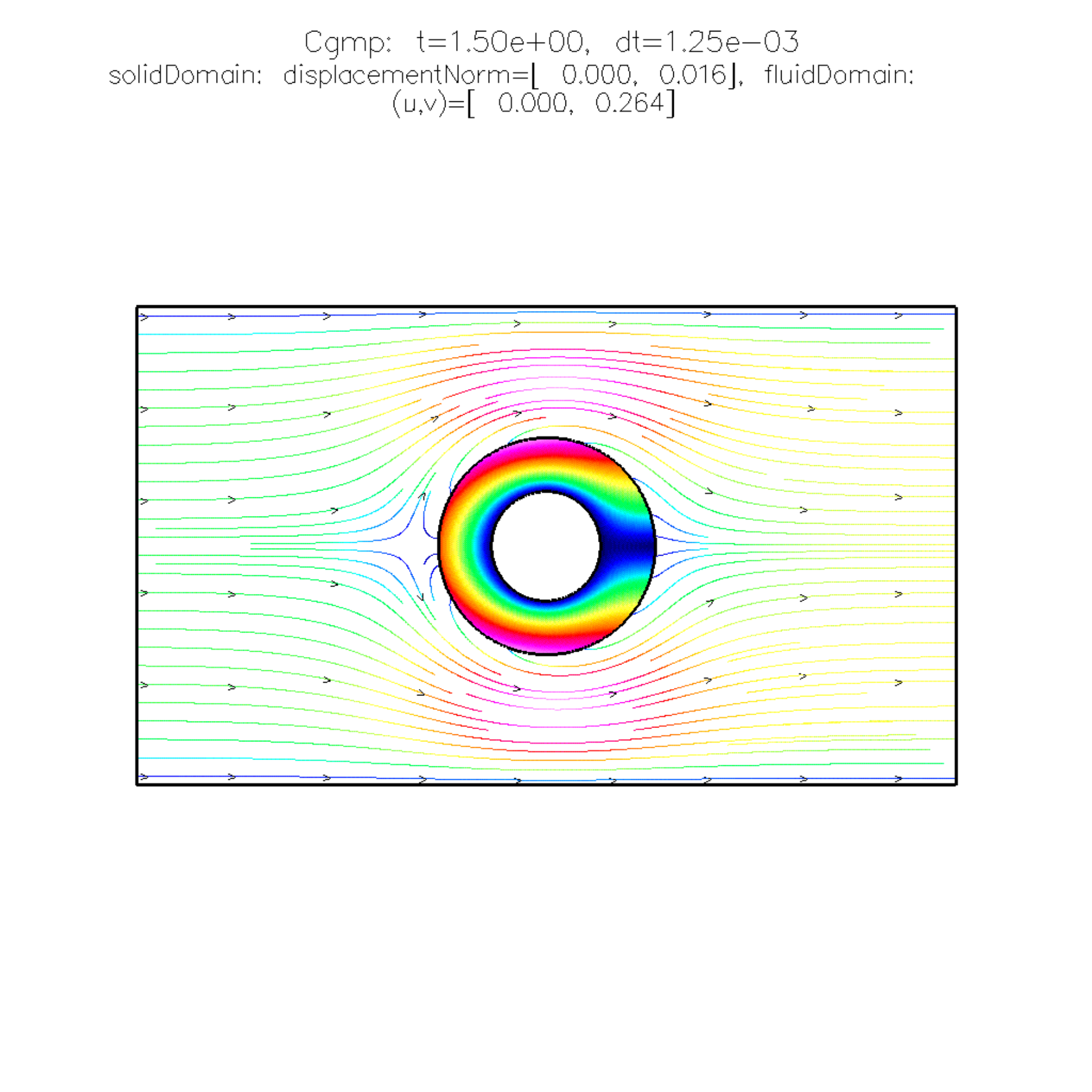}{streamlines}{$|\usv|$}{$t=1$}{$0.0$}{$.264$}{$0.0$}{$.016$};
\end{tikzpicture}
\end{center}
  \caption{Left: composite grid $\Gcd^{(2)}$ for the annulus in a channel.
    The green fluid grid adjacent to the interface moves over time while the blue fluid background grid remains fixed.
    The red grid for the solid reference domain also remains fixed over time.
    Right: streamlines and solid displacement norm $|\usv|$ for $\scf=10$ at $t=1.0$.
}
  \label{fig:diskInAChannel}
\end{figure}
}

The composite grid used for this problem is denoted by $\Gcd^{(j)}$ with target grid spacing equal to $\Delta s=1/(10 j)$.  An
enlarged view of the central portion of the grid is shown in Figure~\ref{fig:diskInAChannel}.
The solid is represented by a single annular (red) grid. The fluid domain
is covered by a background Cartesian (blue) grid and an annular (green) grid
of radial width equal to $0.2$ initially. This latter grid is attached
to the interface and deforms over time.

Figures~\ref{fig:solidDiskInterfaceTime} and~\ref{fig:solidDiskInterface} describe the evolution of the fluid-solid interface
over several snapshots in time. Figure~\ref{fig:solidDiskInterfaceTime} shows the shape of the interface in physical space at different
times by indicating the displacement of a set of marker points on the interface from their reference positions at $t=0$ to their current positions at the times given.
Figure~\ref{fig:solidDiskInterface}, on the other hand, shows the radial and circumferential components of the interface displacement as a function of the polar angle, $\theta$, for the
same set of times.
The angle 
$\theta=0$ corresponds to the trailing edge of the solid and $\theta$ increases in the
counter-clockwise direction.
At time $t=0$ the interface is circular. 
The pressure at the inflow boundary ramps up from $p=0$ at $t=0$ to $p=1$ at $t=1$, and
by $t=2$ the interface has deformed with the left side being pushed to the right.  Following the
initial deformation of the solid to the right by the drag of the fluid, the solid undergoes elastic
oscillations in time about this initially deformed shape.  The oscillations in the solid are damped
over time due to energy transfer to the fluid and the viscous dissipation that occurs there. 
We observe that the outer boundary of the solid achieves its maximum circumferential displacement 
(in magnitude) at approximately $\theta=90^\circ$ and $\theta=270^\circ$, while its maximum radial displacement
occurs at the leading edge $\theta=180^\circ$.

{
\newcommand{\figWidth}{16.0cm}
\newcommand{\trimfig}[2]{\trimw{#1}{#2}{.0}{.00}{.0}{.0}}
\begin{figure}[h]
\begin{center}
\begin{tikzpicture}[scale=1]
\useasboundingbox (0,0.25) rectangle (16.25,-3.2);  

\draw(8.0,0.0) node[anchor=north,xshift=0pt,yshift=0pt] {\trimfig{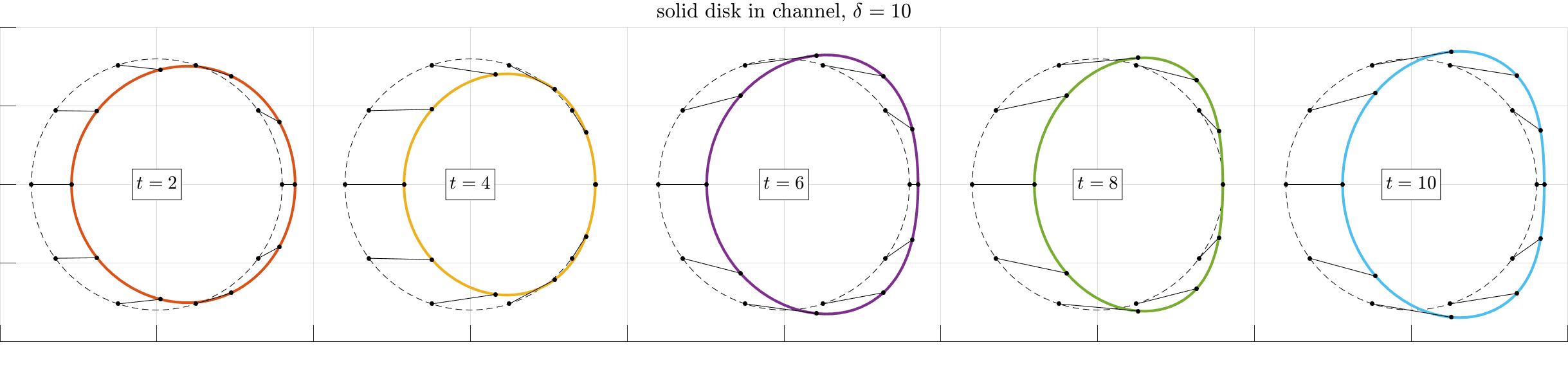}{\figWidth}};

\end{tikzpicture}
\end{center}
\caption{
  Position of the interface over time for the elastic annulus in a channel with density 
  ratio $\scf=10$, computed on grid $\Gcd^{(8)}$. 
  The displacement is magnified by 20 in this plot to accentuate the deformation. 
}
  \label{fig:solidDiskInterfaceTime}
\end{figure}
}

{
\def\width{8}
\def\hscale{.75}
\def\xo{14.5}
\def\yo{-4.25}
\def\R{.9}
\def\RArc{.4}
\def\Angle{45}
\newcommand{\figWidth}{\width cm}
\newcommand{\trimfig}[2]{\trimw{#1}{#2}{.0}{.05}{.0}{.0}}
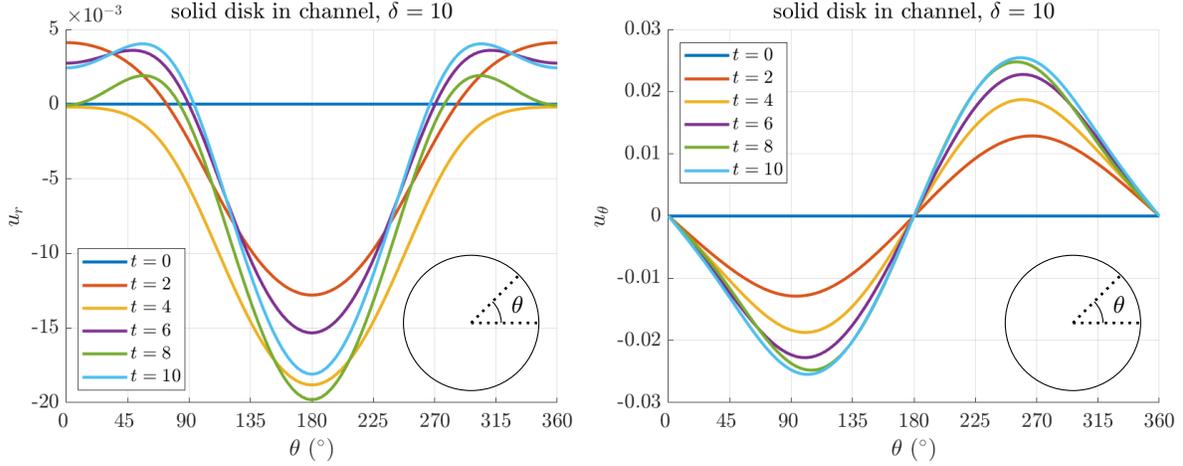
\begin{figure}[h]
\begin{center}
\begin{tikzpicture}[scale=1]
\usetikzlibrary{calc};
\useasboundingbox (-.1,0.1) rectangle (2*\width+.1,-\width*\hscale+.5);  
\draw(0.0,0.0) node[anchor=north west,xshift=-15pt,yshift=8pt] {\trimfig{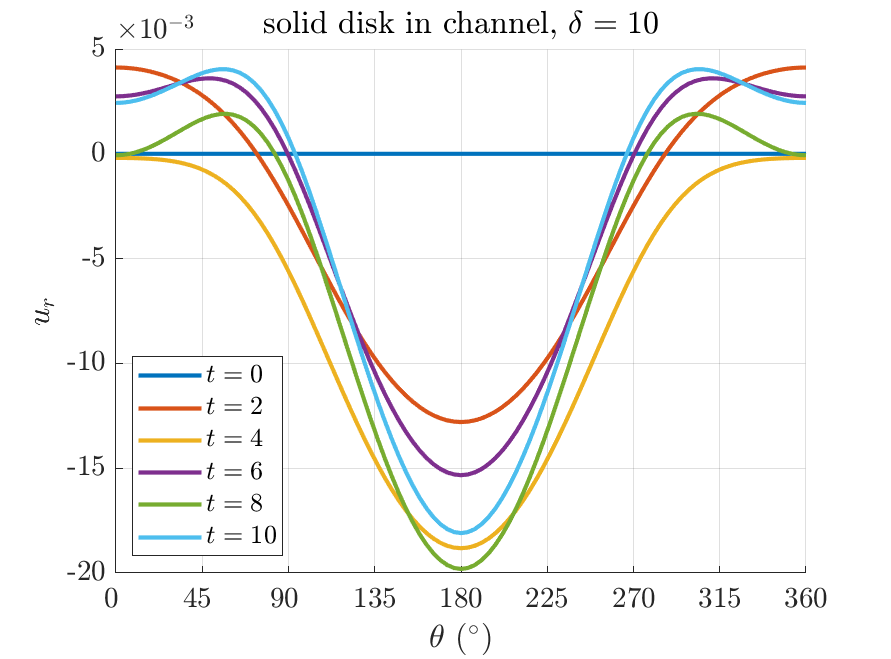}{\figWidth}};
\draw(\width,0.0) node[anchor=north west,xshift=-15pt,yshift=8pt] {\trimfig{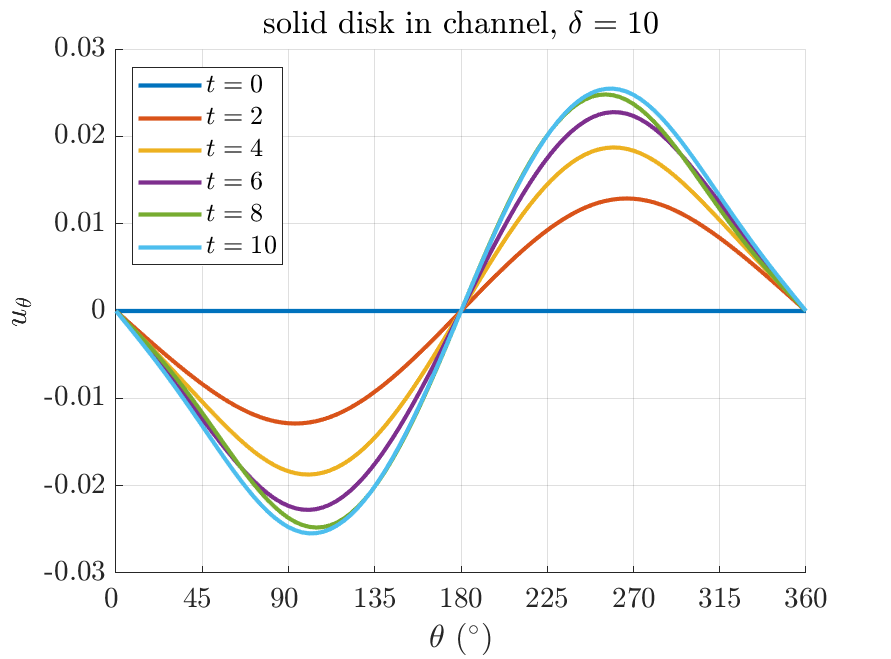}{\figWidth}};

\coordinate (A) at (\xo,\yo);
\coordinate (Rotate) at (\Angle:\R);
\coordinate (NoRotate) at (0:\R);
\draw(A) circle (\R);
\draw[dotted,line width=1pt] (A) -- ($(A)+(NoRotate)$);
\draw[dotted,line width=1pt] (A) -- ($(A)+(Rotate)$);
\draw(\xo+\RArc,\yo) arc (0:\Angle:\RArc);
\draw(\xo+\RArc,\yo) node[anchor=south west] {\footnotesize $\theta$};

\coordinate (A) at (\xo-\width,\yo);
\coordinate (Rotate) at (\Angle:\R);
\coordinate (NoRotate) at (0:\R);
\draw(A) circle (\R);
\draw[dotted,line width=1pt] (A) -- ($(A)+(NoRotate)$);
\draw[dotted,line width=1pt] (A) -- ($(A)+(Rotate)$);
\draw(\xo+\RArc-\width,\yo) arc (0:\Angle:\RArc);
\draw(\xo+\RArc-\width,\yo) node[anchor=south west] {\footnotesize $\theta$};

\end{tikzpicture}
\end{center}
\caption{
  Radial (left) and circumferential (right) components of displacement as a function of the 
  reference interface angle, $\theta,$ for the elastic annulus in a channel with 
  density ratio $\scf=10,$ computed on grid $\Gcd^{(8)}.$
}
  \label{fig:solidDiskInterface}
\end{figure}
}

To estimate the accuracy of the computed solution, a self-convergence test is performed.
Solutions are computed on the composite grids $\Gcd^{(j)}$, $j=2$, 4 and~8, and the errors and
convergence rates are estimated using a Richardson extrapolation procedure
following the approach described in~\cite{pog2008a}. 
Table~\ref{tab:elasticAnnulus} presents the max-norm error results for this self-convergence test.
It is seen that the max-norm errors are converging at close to second-order.

\begin{table}[hbt]\tableFont
\begin{center}
\begin{tabular}{|c|c|c|c|c|c|c|c|c|c|c|} \hline
\multicolumn{11}{|c|}{\strutt Elastic annulus in a channel, $\scf=10$} \\ \hline 
\strutt~~$h_j$~~& $E_j^{(p)}$ & $r$ & $E_j^{(\vv)}$ & $r$ & $E_j^{(\usv)}$ & $r$  & $E_j^{(\vsv)}$  & $r$ & $E_j^{(\sigmasv)}$  & $r$ \\ \hline
 1/20      &  \num{6.3}{-3} &      & \num{7.4}{-3} &      & \num{1.7}{-4} &      & \num{4.1}{-4} &      & \num{1.3}{-2} &       \\ \hline
 1/40      &  \num{1.7}{-3} &  3.7 & \num{1.5}{-3} &  5.0 & \num{4.6}{-5} &  3.7 & \num{9.2}{-5} &  4.5 & \num{3.7}{-3} &  3.6  \\ \hline
 1/80      &  \num{4.7}{-4} &  3.6 & \num{3.0}{-4} &  4.9 & \num{1.2}{-5} &  3.7 & \num{2.0}{-5} &  4.5 & \num{1.0}{-3} &  3.6  \\ \hline
\rateLabel &      1.86      &      &     2.30      &      &    1.90       &      &     2.18      &      &     1.86      &        \\ \hline
\end{tabular}
\caption{Max-norm self-convergence results for the elastic annulus in a channel at $t=1.5$ with density ratio $\scf=10$.}
\end{center}
\label{tab:elasticAnnulus}
\end{table}

Numerical solutions are also computed for this FSI problem using the traditional partitioned (TP) scheme.
Table~\ref{fig:diskInAChannelStability} shows the stability of the computed solution as a function of the
density ratio $\scf=\rhos/\rho$ and the target grid spacing $h$ for the composite grid $\Gcd^{(j)}$.
We observe that the TP scheme is unstable on even the coarsest grid used for $\scf=200$ (and smaller) indicating that
this value is already a \textit{very light} solid case for the TP scheme.  Larger values of $\delta$ are stable on
the coarsest grid, but the results indicate that the TP scheme becomes unstable for a fixed value of $\delta$ as the grid
is refined.  This result is in agreement with the stability results discussed in~\cite{fib2014} and~\cite{fibrmparXiv}.

\newcommand{\bad}{\cellcolor{blue!20}unstable}
\newcommand{\ok}{\cellcolor{green!20}stable}
\begin{table}[H]
\begin{center}
\begin{tabular}{cccccccc}
\hline
$h$    & $\scf=100$ & $\scf=200$ & $\scf=300$ &  $\scf=400$  & $\scf=500$  & $\scf=800$ & $\scf=1000$ \\
\hline
1/20 & \bad&  \bad& \ok & \ok  & \ok & \ok & \ok  \\
1/40 & \bad&  \bad& \bad& \bad & \ok & \ok & \ok  \\
1/80 & \bad&  \bad& \bad& \bad & \bad& \bad& \ok    \\
\hline
\end{tabular}
\end{center}
\caption{ Stability of the traditional-partitioned scheme for the elastic annulus in a channel as a function
of the density ratio $\scf$ and the grid-spacing. Even for very heavy solids the TP scheme will become unstable
when the grid becomes sufficiently fine.}
\label{fig:diskInAChannelStability}
\end{table}






\subsection{Multiple elastic solids in a pressure-driven channel flow} \label{sec:multipleBodies}

To test the robustness of the AMP algorithm, we consider a final example
of an incompressible fluid in a channel flowing past five embedded solids with different 
shapes and densities.
The fluid has density $\rho=1$ and viscosity $\nu=0.02$, and it occupies the rectangular region $[-3.5,5] \times [-2,2]$.
The form of the boundary conditions at the inlet and outlet of the channel, and along its top and bottom walls, are the
same as that described for the pressure-driven flow in Section~\ref{sec:elasticDiskInAChannel}.
For the present problem, the pressure is ramped smoothly from $p=0$ at $t=0$
to $p=p_0=0.1$ at $t = 0.5$  by prescribing $p(-3,y,t)=p_0 \eta (2 t )$ for $-2<y<2$ with the ramp function
$\eta$ given by~\eqref{eq:ramp}.

{
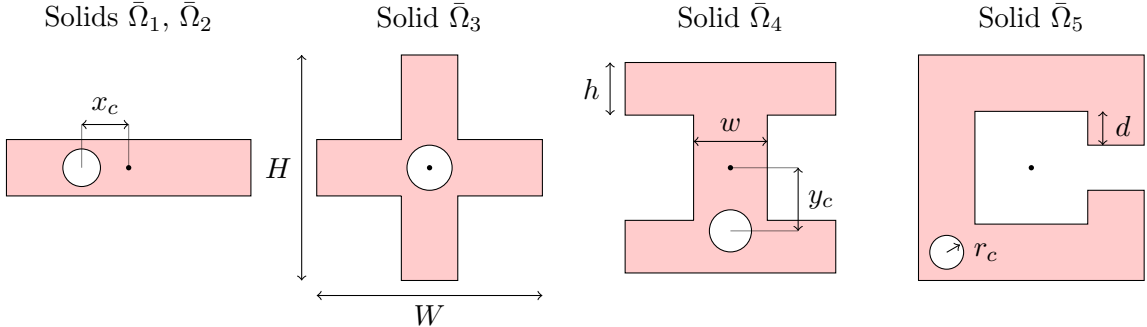
\begin{figure}[h]
\begin{center}
\begin{tikzpicture}
\def\SL{16}
\def\SH{3.5}

\useasboundingbox (0.0,.5) rectangle (\SL,-\SH);

\def\xa{2}
\def\ya{-2}
\def\eh{.3}
\def\ew{1.3}
\def\cr{.1}
\def\xc{.25}
\def\yc{0}
\def\eta{2.5}
\def\scf{.9}
\draw [fill=red!20] (\xa-\eta*.5*\ew,\ya-\eta*.5*\eh) rectangle (\xa+\eta*.5*\ew,\ya+\eta*.5*\eh);
\draw [fill=white] (\xa-\eta*\xc,\ya-\eta*\yc) circle(\eta*\cr);
\fill (\xa,\ya) circle(1pt);

\def\gap{.2}



\draw [<->] (\xa-\eta*\xc,\ya+\eta*.5*\eh+\gap) -- (\xa,\ya+\eta*.5*\eh+\gap);
\draw [] (\xa-\eta*\xc*.5,\ya+\eta*.5*\eh+\gap) node[anchor=south] {$x_{c}$};
\draw [opacity=.5] (\xa-\eta*\xc,\ya+\eta*.5*\eh+\gap) -- (\xa-\eta*\xc,\ya-\eta*\yc);
\draw [opacity=.5] (\xa,\ya+\eta*.5*\eh+\gap) -- (\xa,\ya);
\draw [] (\xa,0) node [] {Solids ${\bar\Omega}_1$, ${\bar\Omega}_2$};

\def\xa{10}
\def\ya{-2}
\def\eh{.25}
\def\ew{1.}
\def\cw{.35}
\def\ch{.5}
\def\xc{0}
\def\yc{.3}
\def\eta{2.8}
\draw [fill=red!20] 
(\xa-.5*\eta*\cw,\ya-.5*\eta*\ch) --
(\xa-.5*\eta*\cw,\ya+.5*\eta*\ch) --
(\xa-.5*\eta*\ew,\ya+.5*\eta*\ch) --
(\xa-.5*\eta*\ew,\ya+.5*\eta*\ch+\eta*\eh) --
(\xa+.5*\eta*\ew,\ya+.5*\eta*\ch+\eta*\eh) --
(\xa+.5*\eta*\ew,\ya+.5*\eta*\ch) --
(\xa+.5*\eta*\cw,\ya+.5*\eta*\ch) --
(\xa+.5*\eta*\cw,\ya-.5*\eta*\ch) --
(\xa+.5*\eta*\ew,\ya-.5*\eta*\ch) --
(\xa+.5*\eta*\ew,\ya-.5*\eta*\ch-\eta*\eh) --
(\xa-.5*\eta*\ew,\ya-.5*\eta*\ch-\eta*\eh) --
(\xa-.5*\eta*\ew,\ya-.5*\eta*\ch) --
(\xa-.5*\eta*\cw,\ya-.5*\eta*\ch);
\draw [fill=white] (\xa-\eta*\xc,\ya-\eta*\yc) circle(\eta*\cr);
\fill (\xa,\ya) circle(1pt);



\draw [<->] (\xa+\eta*.25*\ew+\gap,\ya) -- (\xa+\eta*.25*\ew+\gap,\ya-\eta*\yc);
\draw [opacity=.5] (\xa,\ya) -- (\xa+\eta*.25*\ew+\gap,\ya);
\draw [opacity=.5] (\xa,\ya-\eta*\yc) -- (\xa+\eta*.25*\ew+\gap,\ya-\eta*\yc);
\draw [] (\xa+\eta*.25*\ew+\gap,\ya-.5*\eta*\yc) node[anchor=west] {$y_{c}$};

\draw [<->] (\xa-\eta*.5*\cw,\ya+\eta*.25*\ch) -- (\xa+\eta*.5*\cw,\ya+\eta*.25*\ch);
\draw [] (\xa,\ya+\eta*.25*\ch) node[anchor=south] {$w$};


\draw [<->] (\xa-\eta*.5*\ew-\gap,\ya+\eta*.5*\ch) -- (\xa-\eta*.5*\ew-\gap,\ya+\eta*.5*\ch+\eta*\eh);
\draw [] (\xa-\eta*.5*\ew-\gap,\ya+\eta*.5*\ch+.5*\eta*\eh) node[anchor=east] {$h$};

\draw [] (\xa,0) node [] {Solid ${\bar\Omega}_4$};

\def\xa{6}
\def\ya{-2}
\def\eh{1.}
\def\ew{.25}
\def\cw{.35}
\def\ch{.5}
\def\xc{0}
\def\yc{.3}
\def\eta{3}
\draw [fill=red!20] 
(\xa-.5*\eta*\ew,\ya-.5*\eta*\ew) --
(\xa-.5*\eta*\eh,\ya-.5*\eta*\ew) --
(\xa-.5*\eta*\eh,\ya+.5*\eta*\ew) --
(\xa-.5*\eta*\ew,\ya+.5*\eta*\ew) --
(\xa-.5*\eta*\ew,\ya+.5*\eta*\eh) --
(\xa+.5*\eta*\ew,\ya+.5*\eta*\eh) --
(\xa+.5*\eta*\ew,\ya+.5*\eta*\ew) --
(\xa+.5*\eta*\eh,\ya+.5*\eta*\ew) --
(\xa+.5*\eta*\eh,\ya-.5*\eta*\ew) --
(\xa+.5*\eta*\ew,\ya-.5*\eta*\ew) --
(\xa+.5*\eta*\ew,\ya-.5*\eta*\eh) --
(\xa-.5*\eta*\ew,\ya-.5*\eta*\eh) --
(\xa-.5*\eta*\ew,\ya-.5*\eta*\ew);

\draw [fill=white] (\xa,\ya) circle(\eta*\cr);
\fill (\xa,\ya) circle(1pt);


\draw [<->] (\xa-\eta*.5*\eh,\ya-\eta*.5*\eh-\gap) -- (\xa+\eta*.5*\eh,\ya-\eta*.5*\eh-\gap);
\draw [] (\xa,\ya-\eta*.5*\eh-\gap) node[anchor=north] {$W$};


\draw [<->] (\xa-\eta*.5*\eh-\gap,\ya-\eta*.5*\eh) -- (\xa-\eta*.5*\eh-\gap,\ya+\eta*.5*\eh);
\draw [] (\xa-\eta*.5*\eh-\gap,\ya) node[anchor=east] {$H$};


\draw [] (\xa,0) node [] {Solid ${\bar\Omega}_3$};

\def\xa{14}
\def\ya{-2}
\def\H{1}
\def\W{1}
\def\w{.25}
\def\d{.15}
\def\xc{.3}
\def\yc{.3}
\def\cr{.075}
\def\eta{3}
\draw [fill=red!20]
(\xa-.5*\eta*\W,\ya-.5*\eta*\H) --
(\xa-.5*\eta*\W,\ya+.5*\eta*\H) --
(\xa+.5*\eta*\W,\ya+.5*\eta*\H) --
(\xa+.5*\eta*\W,\ya+.5*\eta*\H-\eta*\w-\eta*\d) --
(\xa+.5*\eta*\W-\eta*\w,\ya+.5*\eta*\H-\eta*\w-\eta*\d) --
(\xa+.5*\eta*\W-\eta*\w,\ya+.5*\eta*\H-\eta*\w) --
(\xa-.5*\eta*\W+\eta*\w,\ya+.5*\eta*\H-\eta*\w) --
(\xa-.5*\eta*\W+\eta*\w,\ya-.5*\eta*\H+\eta*\w) --
(\xa+.5*\eta*\W-\eta*\w,\ya-.5*\eta*\H+\eta*\w) --
(\xa+.5*\eta*\W-\eta*\w,\ya-.5*\eta*\H+\eta*\w+\eta*\d) --
(\xa+.5*\eta*\W,\ya-.5*\eta*\H+\eta*\w+\eta*\d) --
(\xa+.5*\eta*\W,\ya-.5*\eta*\H) --
(\xa-.5*\eta*\W,\ya-.5*\eta*\H);
\draw [fill=white] (\xa-.5*\eta*\W+.5*\eta*\w,\ya-.5*\eta*\H+.5*\eta*\w) circle(\eta*\cr);
\fill (\xa,\ya) circle(1pt);





\draw [<->] (\xa+\eta*.5*\W-\eta*\w+\gap,\ya+\eta*.5*\H-\eta*\w-\eta*\d) -- (\xa+\eta*.5*\W-\eta*\w+\gap,\ya+\eta*.5*\H-\eta*\w);
\draw [] (\xa+\eta*.5*\W-\eta*\w+\gap,\ya+\eta*.5*\H-\eta*\w-.5*\eta*\d) node[xshift=1pt,yshift=1pt,anchor=west] {$d$};



\draw [->] (\xa-.5*\eta*\W+.5*\eta*\w,\ya-.5*\eta*\H+.5*\eta*\w) -- ($(\xa-.5*\eta*\W+.5*\eta*\w,\ya-.5*\eta*\H+.5*\eta*\w)+(30:\eta*\cr*\scf)$);
\draw [] (\xa-.5*\eta*\W+.5*\eta*\w+\eta*\cr,\ya-.5*\eta*\H+.5*\eta*\w) node[yshift=0pt,anchor=west] {$r_c$};

\draw [] (\xa,0) node [] {Solid ${\bar\Omega}_5$};

\end{tikzpicture}
\end{center}
\caption{Geometry of the solids ${\bar\Omega}_k$, $k=1,\ldots,5$, in the channel flow. 
  Each solid shape has an overall height $H$ and width $W$, and may also be described
  by the parameters $h$ and $w$ corresponding to the height and width of rectangular
  segments as shown for ${\bar\Omega}_4$.  The additional parameter $d$ determines
  the height of the top and bottom vertical segments in ${\bar\Omega}_5$.  Each solid includes a circular
  cut-out of radius $r_c$ centered at $(x_c,y_c)$ relative to the geometric center.
  Zero-displacement boundary conditions are applied on the cut-out to anchor each solid
  in the flow.
\label{fig:multiDomainsGeometry}
}
\end{figure}
}

Figure~\ref{fig:multiDomainsGeometry} describes the baseline reference geometries for the elastic
solids ${\bar\Omega}_k$, $k=1,\ldots,5$, in the problem.  These solids are positioned in the channel
as shown in Figure~\ref{fig:multipleBodiesFigGrid}.  The geometric parameters that describe each solid
are given in Table~\ref{tab:multipleBodiesParameters}, along with the position $(x_0,y_0)$ of the
geometric center of each solid in the channel and the angle $\theta_0$ of the solid (in degrees)
relative to the baseline configurations shown in Figure~\ref{fig:multiDomainsGeometry}.
Each solid is anchored to the channel at a small circular hole of radius $r_c$ and center $(x_c,y_c)$ where
the displacement is set to zero.
We set
$\rhos=\lambdas=\mus=\delta$ for each solid with the different values for $\delta$ also listed in
the table. The solid densities, which range from \textit{very light}, $\delta=0.1$, to very \textit{heavy}, $\delta=1000$,
provide a significant test for a partitioned scheme.  Also note that even large values for $\delta$ become difficult for the traditional scheme as the grids are refined as noted from Table~\ref{fig:diskInAChannelStability}.

{
\begin{table}[h]
\centering
\begin{tabular}{ccccccccccc}
Solid & $W$ & $H$ & $w$ & $h$ & $d$ & $(x_c,y_c)$ & $r_c$ & $(x_0,y_0)$ & $\theta_0$ & $\delta$\\
\hline
${\Omegas}_1$ & $1.5$ & $0.25$ & $-$    & $-$    & $-$   & $(0.25,0)$      & $0.05$ & $(-1.5,0.85)$ & $-10^\circ$ & 2\\
${\Omegas}_2$ & $1.5$ & $0.35$ & $-$    & $-$    & $-$   & $(-0.25,0)$     & $0.05$ & $(-1.5,-1)$   & $10^\circ$  & 0.1\\
${\Omegas}_3$ & $1.5$ & $1.5$  & $0.30$ & $0.30$ & $-$   & $(0,0)$         & $0.1$  & $(0,0)$       & $20^\circ$ &  1 \\
${\Omegas}_4$ & $1.0$ & $1.0$  & $0.35$ & $0.25$ & $-$   & $(0,-0.30)$     & $0.05$ & $(1.5,0.75)$  & $45^\circ$  & 10\\
${\Omegas}_5$ & $1.0$ & $1.2$  & $0.25$ & $0.25$ & $0.2$ & $(-0.35,-0.35)$ & $0.05$ & $(1.5,-0.85)$ & $-160^\circ$ & 1000 \\
\hline
\end{tabular}
\caption{Parameters for the solid objects in the channel as shown in Figure~\ref{fig:multipleBodiesFigGrid}.}
  \label{tab:multipleBodiesParameters}
\end{table}
}

\newcommand{\Gmb}{\mathcal{G}_{\rm mb}}
The composite grid, denoted by $\Gmb^{(j)}$ with target grid spacing $\Delta s = 1 / (10 j)$, covers the
fluid domain of the channel and reference domains of the five solids as shown in Figure~\ref{fig:multipleBodiesFigGrid}.
The fluid domain is composed of a (blue) background grid and a (green) deforming 
interface-fitted grid about each of the five solids.
The (static) reference domain of each solid is composed of a (red) background grid, a (purple) interface-fitted 
grid, and a (pink) annular grid about the circular cut-out.  An enlarged view of the composite grid
showing more detail of the overlapping grid structure is also shown in the figure.

{
\newcommand{\figWidth}{7.5cm}
\newcommand{\trimfigZ}[2]{\trimwb{#1}{#2}{.4}{.4}{.45}{.45}}
\newcommand{\trimfig}[2]{\trimwb{#1}{#2}{.0}{.0}{.175}{.15}}
\begin{figure}[htb]
\begin{center}
\begin{tikzpicture}[scale=1]
  \useasboundingbox (0,4.75) rectangle (16,9.25);  

   \draw(0.0,4.0) node[anchor=south west,xshift=-15pt,yshift=0pt] {\trimfig{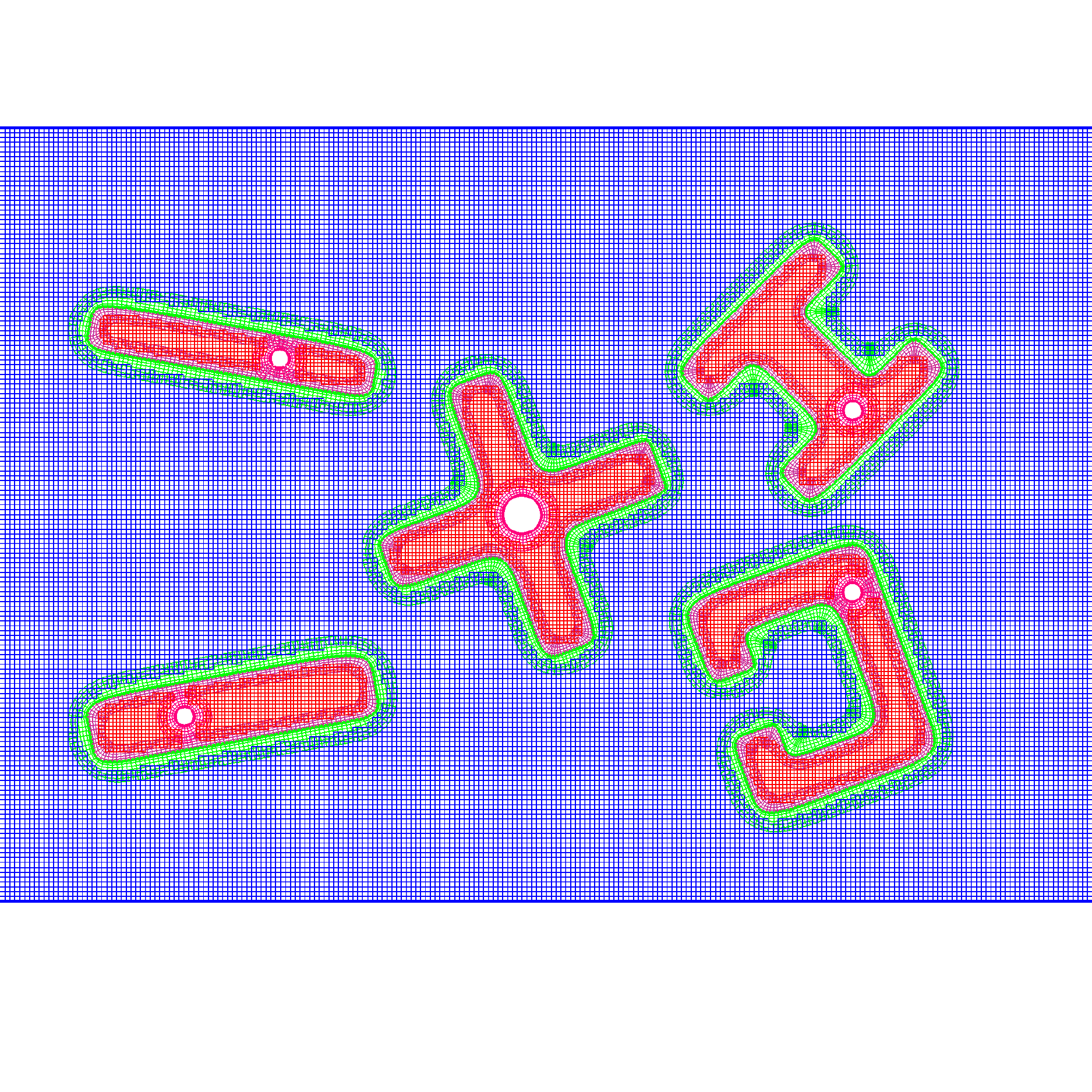}{\figWidth}};
   \draw(8.0,4.0) node[anchor=south west,xshift=-15pt,yshift=0pt] {\trimfig{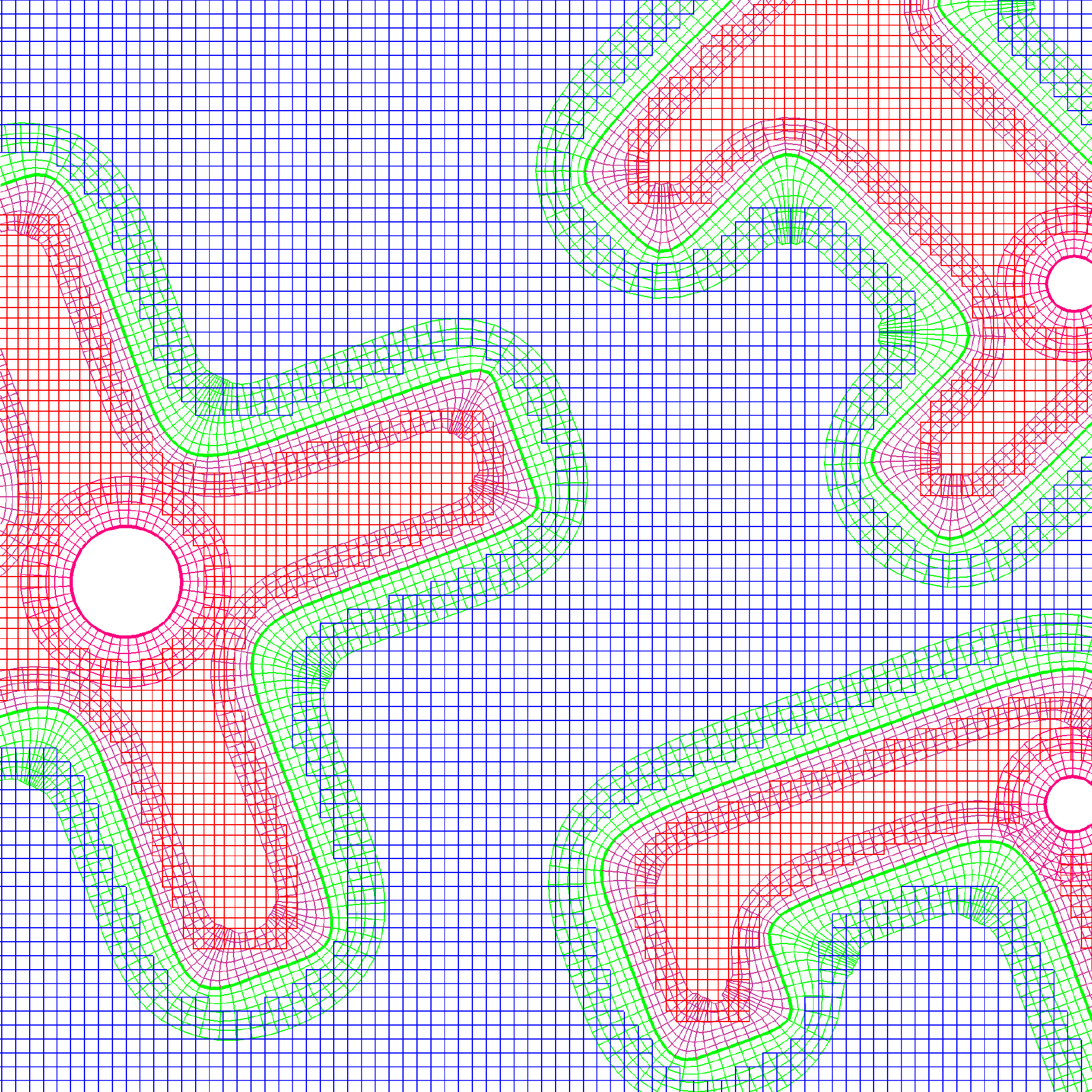}{\figWidth}};





       \draw (1.9,8.5) node[fill=white,draw] {{\footnotesize ${\bar\Omega}_1$}};
       \draw (1.5,6.2) node[fill=white,draw] {{\footnotesize ${\bar\Omega}_2$}};
       \draw (4.1,7.7) node[fill=white,draw] {{\footnotesize ${\bar\Omega}_3$}};
       \draw (6.3,8.5) node[fill=white,draw] {{\footnotesize ${\bar\Omega}_4$}};
       \draw (6.7,6.2) node[fill=white,draw] {{\footnotesize ${\bar\Omega}_5$}};

\end{tikzpicture}
\end{center}
\caption{Left: composite grid $\Gmb^{(4)}$ for the multiple solid domains example.
  Right: zoomed in view of the composite grid $\Gmb^{(4)}$. Only the green fluid grids deform in time, all other grids are static.
  \label{fig:multipleBodiesFigGrid}}
\end{figure}
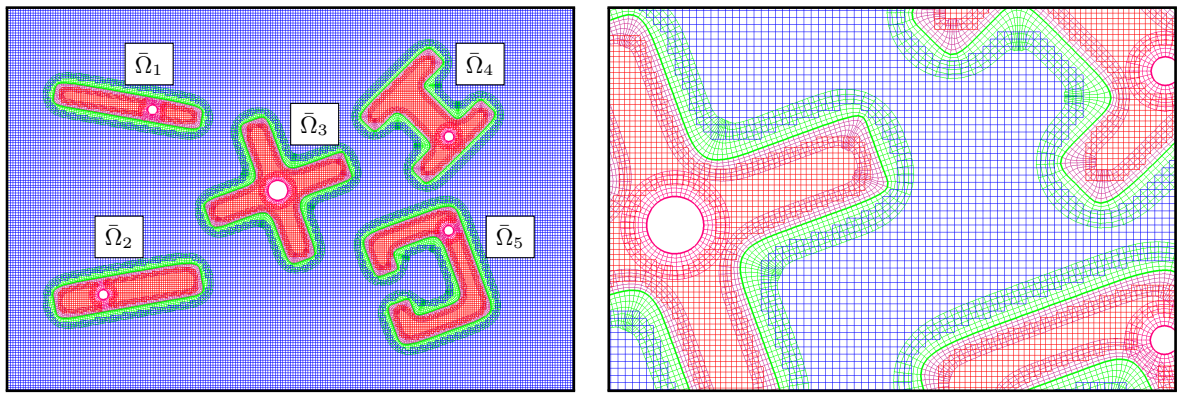
}

\bogus{
\newcommand{\figWidth}{7.8cm}
\newcommand{\trimfigZ}[2]{\trimwb{#1}{#2}{.4}{.4}{.45}{.45}}
\newcommand{\trimfig}[2]{\trimwb{#1}{#2}{.15}{.25}{.35}{.35}}
\begin{figure}[htb]
\begin{center}
\begin{tikzpicture}[scale=1]
  \useasboundingbox (0,4.75) rectangle (16,8.25);  

   \draw(0.0,4.0) node[anchor=south west,xshift=-15pt,yshift=0pt] {\trimfig{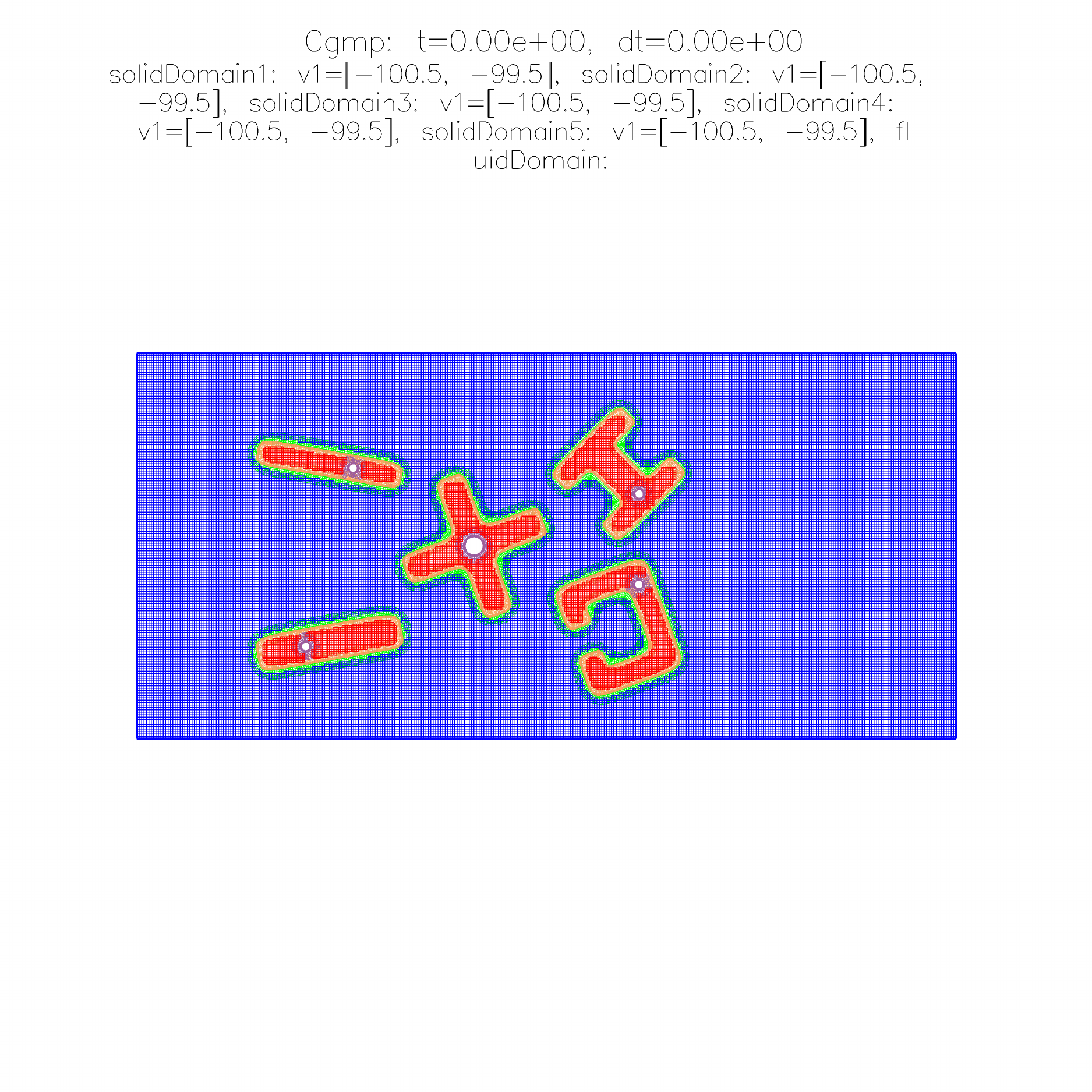}{\figWidth}};
   \draw(8.0,4.0) node[anchor=south west,xshift=-15pt,yshift=0pt] {\trimfigZ{fig/multiObjGrid0}{\figWidth}};





       \draw (2,7.7) node[fill=white,draw] {{\small ${\bar\Omega}_1$}};
       \draw (2,5.6) node[fill=white,draw] {{\small ${\bar\Omega}_2$}};
       \draw (4.0,7) node[fill=white,draw] {{\small ${\bar\Omega}_3$}};
       \draw (6,7.5) node[fill=white,draw] {{\small ${\bar\Omega}_4$}};
   \draw (6.4,5.5) node[fill=white,draw] {{\small ${\bar\Omega}_5$}};

\end{tikzpicture}
\end{center}
\caption{Left: composite grid $\Gmb^{(4)}$ for the multiple solid domains example.
  Right: zoomed in view of composite grid $\Gmb^{(4)}$.
  \label{fig:multipleBodiesFigGrid}}
\end{figure}
}

\bogus{
\newcommand{\figWidth}{7.8cm}
\newcommand{\trimfigZ}[2]{\trimwb{#1}{#2}{.4}{.4}{.45}{.45}}
\newcommand{\trimfig}[2]{\trimwb{#1}{#2}{.15}{.25}{.35}{.35}}
\begin{figure}[htb]
\begin{center}
\begin{tikzpicture}[scale=1]
  \useasboundingbox (0,.5) rectangle (16,8);  

   \draw(8.0,4.0) node[anchor=south west,xshift=-15pt,yshift=0pt] {\trimfigZ{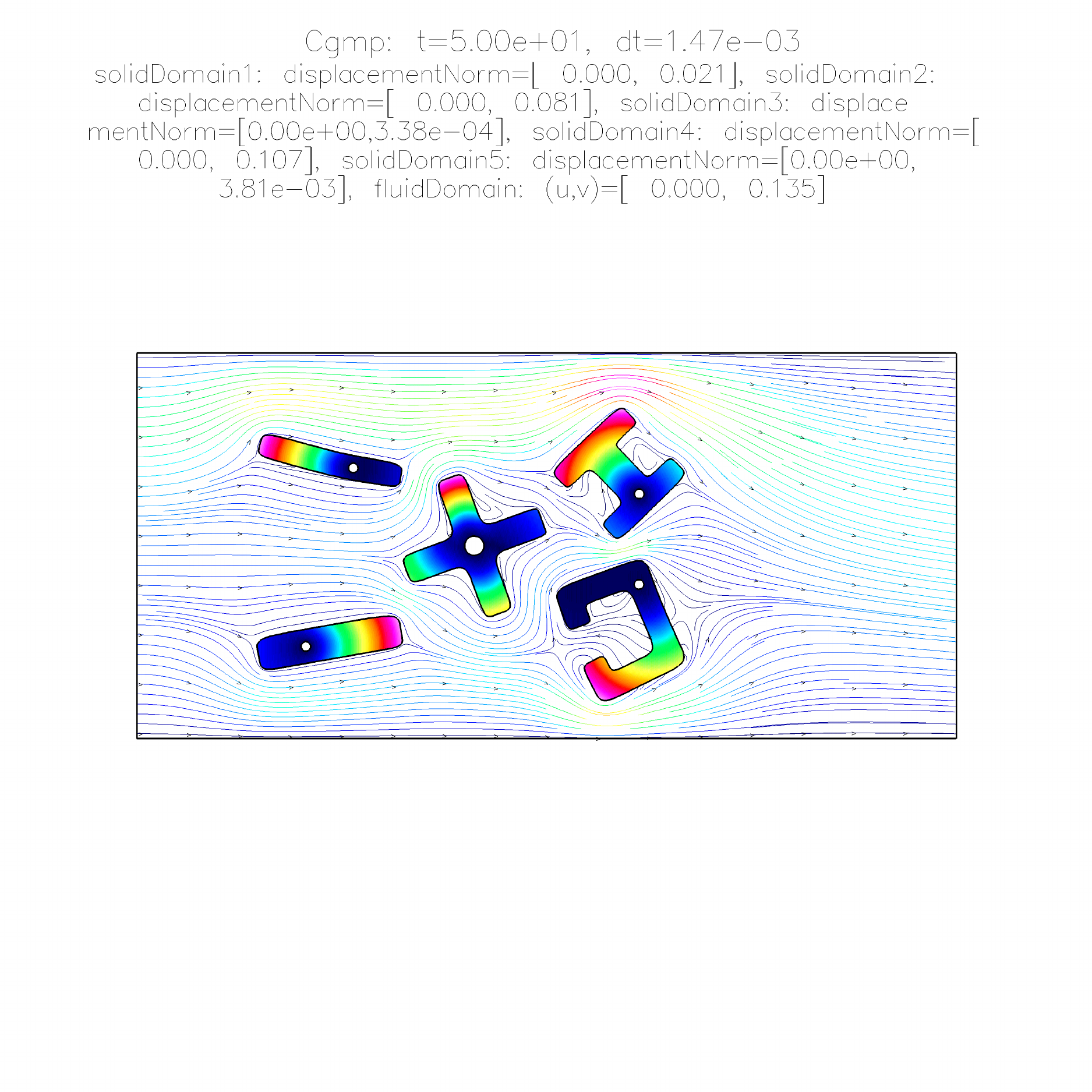}{\figWidth}};
   \draw(8.0,0.0) node[anchor=south west,xshift=-15pt,yshift=-5pt] {\trimfigZ{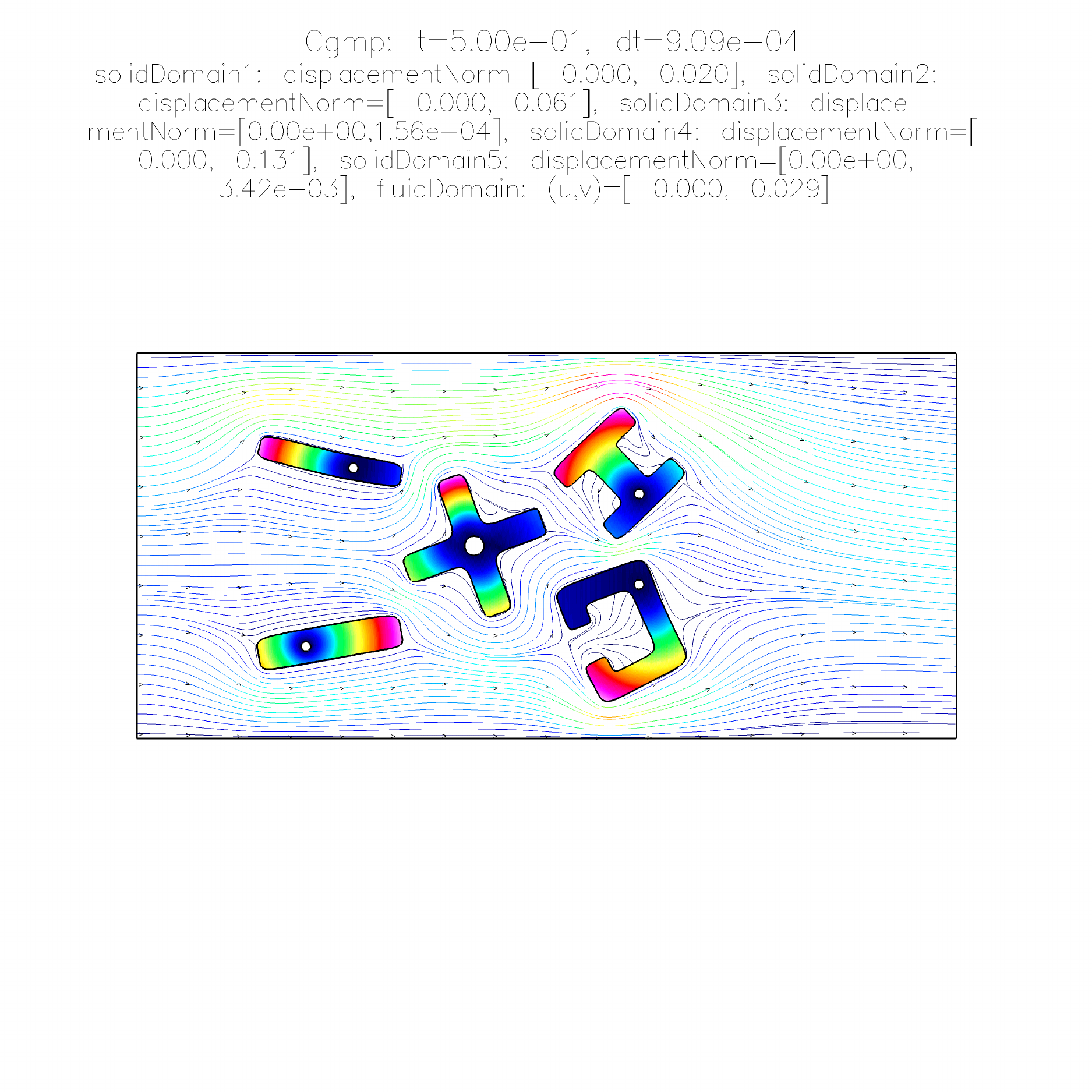}{\figWidth}};

   \draw(0.0,4.0) node[anchor=south west,xshift=-15pt,yshift=0pt] {\trimfig{fig/fiveBodiesG4_Pp1}{\figWidth}};
   \draw(0.0,0.0) node[anchor=south west,xshift=-15pt,yshift=-5pt] {\trimfig{fig/fiveBodiesG8_Pp1}{\figWidth}};

   \draw(0.0,4) node[yshift=-4pt,xshift=2pt,draw,fill=white,anchor=north west] {\small $t=50$};
   \draw(0.0,8) node[yshift=1pt,xshift=2pt,draw,fill=white,anchor=north west] {\small $t=50$};

   \draw(8,8) node[yshift=1pt,xshift=-6pt,draw,fill=white,anchor=north east] {\small $\Gmb^{(4)}$};
   \draw(8,4) node[yshift=-4pt,xshift=-6pt,draw,fill=white,anchor=north east] {\small $\Gmb^{(8)}$};



\end{tikzpicture}
\end{center}
\caption{Plot of streamlines in the fluid and displacement norm in the solid at $t=50$ 
  computed on $\Gmb^{(4)}$ (top) and $\Gmb^{(8)}$ (bottom).
  \label{fig:multipleBodiesFigFlow}}
\end{figure}
}

\input texFiles/fiveBodiesFlowFig

Figure~\ref{fig:multipleBodiesFigFlow} shows computed results at two times, $t=8$ and $t=10$, with a comparison
between two different grid resolutions, a `coarse' grid $\Gmb^{(4)}$ (left column) and a `fine' grid $\Gmb^{(8)}$ (right column).
The plots in the figure  show (instantaneous) streamlines in the fluid domain
and shaded contours of the magnitude of displacement in each solid domain. 
The shaded contours are colored on different scales for each solid as noted in the figure caption.
The streamlines are colored by the magnitude of the flow speed. The solids are seen to divert the flow with the
fastest flow found above the top of the I-shaped solid $\Omegas_4$. 
The flow is time-dependent and the
solids undergo oscillations over time as elastic waves propagate through the bodies. 
The lighter solids are seen to deform more than the heavier ones. The rectangular beam $\Omegas_2$ at the bottom
near the inlet of the channel, for example, is very light and is easily deformed by the fluid.  As a result, it also undergoes
a fairly large compression and expansion. Three arms of the plus-shaped
solid $\Omegas_3$,
are seen to be significantly pushed by the fluid while the top of the I-shaped solid is deformed to the right.
Note that since the deformation of the bodies is time-dependent, the matching fluid velocity can be non-zero on the fluid-solid interfaces and
the fluid streamlines appear to enter or exit the bodies
at some locations. The results from the coarse and fine grid are very similar indicating a
measure of grid convergence. Importantly, the AMP scheme is shown to be robust and stable for solid bodies
with a wide range of densities.

\section{Conclusions} \label{sec:conclusions}

A stable added-mass partitioned (AMP) algorithm was developed for fluid-structure interaction
problems involving viscous incompressible flow and compressible elastic solids.  The new algorithm
is stable, without sub-time-step iterations, for both heavy and very light solids and effectively
suppresses both added-mass and added-damping effects.
Key elements of the new AMP scheme are a Robin interface condition for
the pressure and an impedance-weighted interface projection based on a new form for the fluid
impedance.  The scheme was implemented in two dimensions using deforming composite grids to treat
complex geometry with multiple deforming bodies. The fluid is advanced using a fractional-step IMEX
scheme with the viscous terms treated implicitly, while the hyperbolic system for the solid is
integrated using an explicit upwind scheme.  Both domain solvers are individually second-order accurate in both time
and space.  A number of new benchmark problems were developed
including a \textit{radial-piston} problem where exact solutions for radial and azimuthal motions
were found. Traveling-wave exact solutions for a solid disk surrounded by an annulus of fluid were
also derived.  These exact solutions were used to verify the stability and second-order accuracy of
the scheme.  A clean benchmark problem involving a pressure-driven flow in a channel past a deformable solid annulus was
developed and simulated; a self-convergence grid-refinement study was used to demonstrate
second-order accuracy.  A final example of flow past five solid bodies with densities ranging from $10^{-1}$
to $10^3$ illustrated the flexibility and robustness of the new scheme.

\appendix

\section{Radial piston exact solution} \label{sec:radialElasticPistonExactSolution}

The geometry of a FSI problem involving radial motion of a circularly symmetric elastic piston is shown in Figure~\ref{fig:radialElasticPistonGeometry}.  An incompressible fluid occupies the annular region between the fluid-solid interface at $r=r_I(t)$ and an outer boundary at $r=R$.  The fluid surrounds an elastic solid disk for $0<r<r_I(t)$.  It is assumed that the solid deforms in the radial direction only (as shown in the figure), and that the corresponding fluid motion is also confined to the radial direction.  Under these assumptions, the fluid continuity equation \eqref{eq:fluidContinuity} in cylindrical coordinates reduces to 
\begin{align}
\frac{1}{r} \frac{\partial }{\partial r} (r v_r) = 0, \qquad r_I(t)<r<R,
\end{align}
so that the radial component of the fluid velocity is given by
\begin{align}
v_r(r,t) = \frac{R}{r} V(t),
\label{eq:radialSolidFluidVelocity}
\end{align}
where $V(t)$ is the radial velocity of the fluid at $r=R$.  The radial component of the momentum equation \eqref{eq:fluidMomentum} reduces to
\begin{align}
\frac{\partial p}{ \partial r} = \frac{\rho R\sp2}{r^3} V(t)^2 - \frac{\rho R}{r} \dot{V}(t), \qquad r_I(t)<r<R.
\end{align}
Integration with respect to $r$ gives
\begin{align}
p(r,t) = P(t)+{\rho\over2}\left(1-{R\sp2\over r\sp2}\right)V(t)\sp2+\rho R\log\left({R\over r}\right)\dot{V}(t),
\label{eq:radialSolidFluidPressure}
\end{align}
where $P(t)$ is the fluid pressure at $r=R$.  The evolution of the fluid velocity and pressure is completely determined by the interface position, $r_I(t)$, and the velocity and pressure at the outer radius, $V(t)$ and $P(t)$.  These functions are found by matching the evolution of the fluid with that of the solid, which is considered next.

In the solid reference domain, the radial component of the momentum equation reduces to
\begin{align}
\frac{\partial^2 \us_r}{\partial t^2} = 
\cp^2 \frac{\partial}{\partial \rs} 
\left(  \frac{1}{\rs} \frac{\partial}{\partial \rs} \rs \us_r\right), \qquad 0<\rs<r_0,
\label{eq:radialSolidMomentum}
\end{align}
where $\us_r(\rs,t)$ is the radial component of the solid displacement and $r_0=r_I(0)$ is the initial radius of the interface.  Bounded solutions of~\eqref{eq:radialSolidMomentum} as $\rs\rightarrow0$ are sought in the form
\begin{align}
\us_r(\rs,t) = \beta J_1(\omega\rs/\cp) \sin(\omega t),
\label{eq:radialSolidDisplacement}
\end{align}
where $J_1$ is the Bessel function of the first kind of order one, $\omega$ is the frequency of the oscillating displacement, and $\beta$ determines the magnitude of the oscillation.  In terms of the radial displacement in~\eqref{eq:radialSolidDisplacement}, the interface position in physical space is given by
\begin{equation}
r_I(t)=r_0+\us_r(r_0,t)=r_0+b\sin(\omega t),\qquad b=\beta J_1(\omega r_0/\cp).
\label{eq:radialInterfacePosition}
\end{equation}
Matching the radial velocity of the fluid given in~\eqref{eq:radialSolidFluidVelocity} to that of the solid determined from~\eqref{eq:radialSolidDisplacement} on the interface gives
\begin{equation}
{R\over r_I(t)}V(t)=\omega\beta J_1(\omega r_0/\cp) \cos(\omega t)
\label{eq:radialOuterFluidVelocity}
\end{equation}
which specifies $V(t)$.  The matching condition on the radial component of stress is
\begin{equation}
-p(r,t)={\lambdas\over\rs}{\partial\over\partial\rs}(\rs\us_r)+2\mus{\partial\us_r\over\partial\rs},\qquad r=r_I(t),\quad \rs=r_0.
\label{eq:radialMatchingStress}
\end{equation}
Using~\eqref{eq:radialSolidFluidPressure} and~\eqref{eq:radialSolidDisplacement} in the stress condition~\eqref{eq:radialMatchingStress} gives
\begin{equation}
\begin{array}{l}
\displaystyle{
-P(t)={\rho\over2}\left(1-{R\sp2\over r_I(t)\sp2}\right)V(t)\sp2+\rho R\log\left({R\over r_I(t)}\right)\dot{V}(t)
}\medskip\\
\qquad\qquad\quad\displaystyle{
+\beta\left[(\lambdas+2\mus){\omega\over\cp}J_1\sp\prime(\omega r_0/\cp)+{\lambdas\over r_0}J_1(\omega r_0/\cp)\right]\sin(\omega t)
}
\end{array}
\label{eq:radialOuterFluidPressure}
\end{equation}

An exact solution of the radial elastic piston problem is available based on the choice for the solid displacement in~\eqref{eq:radialSolidDisplacement}, which leads to the radial position of the fluid-solid interface given in~\eqref{eq:radialInterfacePosition}.  Matching the radial velocity and stress at the interface specifies $V(t)$ and $P(t)$ from~\eqref{eq:radialOuterFluidVelocity} and~\eqref{eq:radialOuterFluidPressure}, respectively, and these functions determine the fluid velocity and pressure given in~\eqref{eq:radialSolidFluidVelocity} and~\eqref{eq:radialSolidFluidPressure}.

\section{Exact solution for the rotating elastic disk in a fluid annulus} \label{sec:rotatingElasticDiskExactSolution}

An exact solution can be constructed for the case of circumferential motion of a circularly symmetric elastic disk.  It is assumed that the radial components of the fluid velocity and the solid displacement are zero, and thus the interface is located at $r=r_0$ for all time.  It is also assumed that the circumferential components of the fluid velocity, $v_\theta$, and the solid displacement, $\us_\theta$, depend on $r$ and $t$ alone, so that the equations governing these quantities reduce to
\begin{equation}
\frac{\partial v_\theta}{\partial t} = \nu \frac{\partial}{\partial r}
\left( \frac{1}{r} \frac{\partial}{\partial r} r v_\theta\right), \qquad r_0<r<R,
\label{eq:fluidRotating}
\end{equation}
and
\begin{equation}
\frac{\partial^2 \us_\theta}{\partial t^2} = \cs^2 \frac{\partial}{\partial r}
\left( \frac{1}{r} \frac{\partial}{\partial r} r \us_\theta\right), \qquad 0<r<r_0,
\label{eq:solidRotating}
\end{equation}
respectively.  Since there is no displacement in the radial direction, the reference coordinate, $\rs$, is equivalent to the physical coordinate, $r$.  Setting
\begin{equation}
v_\theta(r,t) = \hat{v}_\theta(r) e^{i \omega t}, \qquad \us_\theta(r,t) = \hat{\us}_\theta(r) e^{i \omega t},
\label{eq:radialTransSolutions}
\end{equation}
and substituting these into~\eqref{eq:fluidRotating} and~\eqref{eq:solidRotating}
it is found that the coefficient functions, $\hat{v}_\theta(r)$ and $\hat{\us}_\theta(r)$, both satisfy Bessel equations of order one.  Assuming a no-slip boundary condition on the fluid velocity at $r=R$ and a bounded solution for the solid displacement as $r\rightarrow0$, gives solutions of the form
\begin{equation}
\hat{v}_\theta(r)=b\bigl[J_1(\lambda r)Y_1(\lambda R)-J_1(\lambda R)Y_1(\lambda r)\bigr],\qquad \hat{\us}_\theta(r)=\bar bJ_1(k_sr),
\label{eq:radialTransCoefficients}
\end{equation}
where $J_1$ and $Y_1$ are Bessel functions of the first and second kind, respectively.  
Here, $\lambda^2= i \omega / \nu,$ $k_s = \omega / \cs,$ 
and $(b,\bar b)$ are constants.  
Matching conditions on velocity and stress imply
\begin{equation}
\hat{v}_\theta(r) = i\omega\hat{\us}_\theta(r),  \qquad \mu{d\over dr}\left({\hat{v}_\theta(r)\over r}\right) = \mus{d\over dr}\left({\hat{\us}_\theta(r)\over r}\right),\qquad \hbox{at $r=r_0$.}
\label{eq:radialTransMatching}
\end{equation}
Substituting~\eqref{eq:radialTransCoefficients} into the matching conditions in~\eqref{eq:radialTransMatching}, and using the identity
\[
{d\over d\xi}\left({Z_1(\xi)\over\xi}\right)=-{Z_2(\xi)\over\xi},\qquad\hbox{$Z_n=J_n$ or $Y_n$,}
\]
leads to the system
\begin{equation}
\left[\begin{array}{cc}
J_1(\lambda r_0)Y_1(\lambda R)-J_1(\lambda R)Y_1(\lambda r_0) & -i \omega J_1(k_sr_0) \smallskip\\
-\mu \lambda\bigl[J_2(\lambda r_0)Y_1(\lambda R)-J_1(\lambda R)Y_2(\lambda r_0)\bigr] & \mus k_s J_2(k_sr_0)
\end{array}\right]
\left[\begin{array}{c}
b \smallskip\\ \bar{b}
\end{array}\right]
= 0 .
\label{eq:theBs}
\end{equation}
Nontrivial solutions for $(b,\bar b)$ exist if the following determinant condition is met:
\begin{equation}
\begin{array}{l}
\displaystyle{
\mathcal{D}_2(\omega)=\mus k_s J_2(k_sr_0)\bigl[J_1(\lambda r_0)Y_1(\lambda R)-J_1(\lambda R)Y_1(\lambda r_0)\bigr]
}\medskip\\
\qquad\qquad\quad\displaystyle{
-i \omega\mu \lambda J_1(k_sr_0)\bigl[J_2(\lambda r_0)Y_1(\lambda R)-J_1(\lambda R)Y_2(\lambda r_0)\bigr]
}.
\end{array}
\label{eq:radialTransConstraint}
\end{equation}

Values of $\omega$ with $\Re(\omega) > 0$ satisfying $\mathcal{D}_2(\omega)=0$ imply solutions of the rotating elastic disk problem that decay in time.  For such a value of $\omega$, a unique choice for $(b,\bar b)$, satisfying~\eqref{eq:theBs}, can be expressed as
\begin{align}
\bar{b} = \frac{\us_0}{J_1(k_s r_0)}, \qquad
b = \frac{i \omega \us_0}
{J_1(\lambda r_0) Y_1 (\lambda r_0) - J_1 (\lambda R) Y_1(\lambda r_0)},
\end{align}
where $\us_0=\us_\theta(r_0,0)$ is the initial displacement of the interface.
%
Solutions for the circumferential components of the fluid velocity and the solid displacement are taken from the real parts of $v_\theta(r,t)$ and $\us_\theta(r,t)$ given in~\eqref{eq:radialTransSolutions}, while the fluid pressure is given by
\[
p(r,t) = \rho \int_{r_0}^r \frac{v_\theta(s,t)^2}{s} \, ds,\qquad r_0<r<R,
\]
assuming $p=0$ on the interface $r=r_I(0)=r_0$.



\section{Radial traveling wave solution} \label{sec:radialTravelingWaveSolution}

We consider traveling wave solutions of an FSI problem involving a Stokes fluid in the circular region
$0 < r < r_0$ initially surrounded by a linearly elastic solid in the annular region $r_0 < r < R$.
The solution is linearized about small interface deformations, therefore 
it is assumed that the interface is located at $r=r_0$ for all time. 

In the fluid, solutions are sought of the form
\begin{align*}
v_r(r,\theta,t) = \hat{v}_r(r) e^{i (n \theta - \omega t)}, \qquad
v_\theta(r,\theta,t) = \hat{v}_\theta(r) e^{i (n \theta - \omega t)} , \qquad
p(r,\theta,t) = \hat{p}(r) e^{i (n \theta - \omega t)} ,
\end{align*}
where $\omega$ is a frequency, $n$ is a positive integer, and $\hat{v}_r(r)$, $\hat{v}_\theta(r)$ and $\hat{p}(r)$
are coefficient functions.
Upon substitution, the equations governing the Stokes fluid become
\begin{align}
\frac{1}{r} (r \hat{v}_r)' + \frac{i n}{r} \hat{v}_\theta &= 0  ,
\label{eq:smallDistContinuity} \\
-i \omega \rho \hat{v}_r + \hat{p}' &= \mu 
\left(\frac{1}{r} (r \hat{v}_r')' 
- \frac{n^2}{r^2} \hat{v}_r 
- \frac{\hat{v}_r}{r^2} 
- \frac{2 i n}{r^2} \hat{v}_\theta
\right), \label{eq:smallDistRadialMomentum}\\
-i \omega \rho \hat{v}_\theta + \frac{i n}{r} \hat{p}' &= \mu
\left(\frac{1}{r} (r \hat{v}_\theta')'
- \frac{n^2}{r^2} \hat{v}_\theta
+ \frac{2 i n}{r^2} \hat{v}_r
- \frac{\hat{v}_\theta}{r^2} 
\right).\label{eq:smallDistCircumMomentum}
\end{align}
An elliptic equation for the pressure can be obtained from~\eqref{eq:smallDistContinuity}--\eqref{eq:smallDistCircumMomentum}, and it takes the form
\begin{align}
\hat{p}'' 
+ \frac{1}{r} \hat{p}'
- \frac{n^2}{r^2} \hat{p} &= 0. 
\label{eq:phat}
\end{align}
Bounded solutions of the Cauchy-Euler equation in~\eqref{eq:phat} at $r=0$ are
\begin{align}
\hat{p}(r) = p_I \left(\frac{r}{r_0} \right)^n ,
\end{align}
where $p_I$ is a constant corresponding to the pressure at the (linearized) fluid-solid interface.
Using~\eqref{eq:smallDistContinuity} to eliminate $\hat{v}_\theta$
from~\eqref{eq:smallDistRadialMomentum} gives
\begin{align}
r^2 \hat{v}_r'' + 3 r \hat{v}_r' + (\lambda^2 r^2 - n^2+1) \hat{v}_r
&= p_I \frac{n r^{n+1}}{\mu r_0^n},
\label{eq:besselvr}
\end{align}
where $\lambda^2 = i \omega / \nu$.  Bounded solutions at $r=0$ have the form
\begin{align}
\hat{v}_r(r) &= d \frac{J_n(\lambda r)}{r}
+ p_I \frac{n r^{n-1}}{\mu \lambda^2 r_0^n},
\end{align}
where $d$ is a constant and $J_n(z)$ is the Bessel function of the first kind of order $n$.
The corresponding solution for $\hat{v}_\theta$, obtained from~\eqref{eq:smallDistContinuity}, is given by
\begin{align}
\hat{v}_{\theta}(r) &=
d \frac{i\lambda}{n} J_n'(\lambda r)
{+} p_I \frac{i n r^{n-1}}{\mu \lambda^2 r_0^n}.
\end{align}

In the solid, we consider the solid displacement in the form
\begin{align}
\usv = \grad \phi + \grad \curl \Hv,
\label{eq:displacementDecomp}
\end{align}
where $\phi$ is a potential and $\Hv$ is a vector field.
Let $\ev_r = [\cos \theta, \sin \theta,0]^T,$ $\ev_\theta = [-\sin \theta, \cos \theta,0]^T,$
and $\ev_z = [0,0,1]^T$ represent the polar basis.  Since we are interested in displacements in the two-dimensional $(r,\theta)$ plane,
we take $\Hv = H \ev_z$. With these assumptions, the equations governing linear elasticity reduce to a longitudinal and shear wave equations for
$\phi(r,\theta,t)$ and $H(r,\theta,t)$, respectively, having the form
\begin{align}
{\partial\sp2\phi\over\partial t\sp2} = \cp^2 \left(\frac{1}{r} \frac{\partial }{\partial r} \left( r \frac{\partial \phi}{\partial r} \right) 
  + \frac{1}{r^2} \frac{\partial^2 \phi }{\partial \theta^2} \right), \qquad
{\partial\sp2 H\over\partial t\sp2} = \cs^2 \left(\frac{1}{r} \frac{\partial }{\partial r} \left( r \frac{\partial H}{\partial r} \right) 
  + \frac{1}{r^2} \frac{\partial^2 H }{\partial \theta^2} \right) .
\label{eq:solidequations}
\end{align}
As before, we seek solutions of the form
\begin{align}
\phi(r,\theta,t) = \hat{\phi}(r) e^{i (n \theta - \omega t)}, \qquad
H(r,\theta,t) = \hat{H}(r) e^{i (n \theta - \omega t)}.
\end{align}
Upon substitution into~\eqref{eq:solidequations}, we find
\begin{align}
r^2 \hat{\phi}'' + r \hat{\phi} + \left((k_p r)^2 - n^2 \right) \hat{\phi}
= 0, \qquad
r^2 \hat{H}'' + r \hat{H} + \left((k_s r)^2 - n^2 \right) \hat{H}
= 0,
\label{eq:solidcoeffeqns}
\end{align}
where $k_p=\omega/\cp$ and $k_s=\omega/\cs$.  General solutions of the two equations in~\eqref{eq:solidcoeffeqns} involve Bessel functions of the first and second kinds, and are given by
\newcommand{\Bn}{\db_1 J_n(k_p r) + \db_2 Y_n(k_p r)}
\newcommand{\Bnp}{\db_1 J_n'(k_p r) + \db_2 Y_n'(k_p r)}
\newcommand{\Bnpp}{\db_1 J_n''(k_p r) + \db_2 Y_n''(k_p r)}
\newcommand{\Cn}{\db_3 J_n(k_s r) + \db_4 Y_n(k_s r)}
\newcommand{\Cnp}{\db_3 J_n'(k_s r) + \db_4 Y_n'(k_s r)}
\newcommand{\Cnpp}{\db_3 J_n''(k_s r) + \db_4 Y_n''(k_s r)}
\begin{align}
\hat{\phi}(r) = \Bn , \qquad
\hat{H}(r) = \Cn ,
\end{align}
where $(\db_1,\db_2,\db_3,\db_4)$ are constants.
%
%
Using~\eqref{eq:displacementDecomp}, the coefficient functions for the solid displacements in the radial and circumferential directions are given by
\begin{align}
\hat{\us}_r(r) = \hat{\phi}\sp\prime + {in\over r}\hat{H},\qquad \hat{\us}_\theta(r) = {in\over r}\hat{\phi} - \hat{H}\sp\prime,
\end{align}
respectively.  For later purposes, we also require coefficient functions for the stress components $(\sigmas_{rr},\sigmas_{r\theta})$ in terms of $\hat\phi$ and $\hat H$, which are given by
\begin{align}
\hat{\sigmas}_{rr}(r)&=
(2 \mus + \lambdas) \hat{\phi}''
+ \frac{\lambdas}{r} \left(\hat{\phi}' - \frac{n^2}{r} \hat{\phi} \right)
+ \frac{2 i \mus n}{r} \left(\hat{H}' - \frac{1}{r} \hat{H} \right),\\
 \hat{\sigmas}_{r\theta}(r)&=
-\mus \hat{H}''
+ \frac{\mus}{r} \left(\hat{H}' - \frac{n^2}{r} \hat{H} \right)
+ \frac{2 i \mus n}{r} \left(\hat{\phi}' - \frac{1}{r} \hat{\phi} \right)
 .
\end{align}

The solution for the fluid involves the two constants $(p_I,d)$, while the solution for the solid involves the four constants $(\db_1,\db_2,\db_3,\db_4)$.  Constraints for these six constants can be obtained by imposing boundary conditions on the solid at $r=R$ and the matching conditions at the fluid-solid interface, $r=r_0$.  At the outer boundary of the solid, we assume traction-free boundary conditions of the form
\begin{align}
\hat{\sigmas}_{rr}(R) = \hat{\sigmas}_{r\theta}(R) = 0,
\label{eq:radialWaveFreeSurface}
\end{align}
while matching conditions on the velocity and stress at the interface give
\begin{align}
\hat{v}_r(r_0) = -i \omega \hat{\us}_r(r_0),
\quad
\hat{v}_\theta(r_0) = -i \omega \hat{\us}_\theta(r_0),
\quad
\hat{\sigma}_{rr}(r_0) = \hat{\sigmas}_{rr}(r_0),
\quad
\hat{\sigma}_{r\theta}(r_0) = \hat{\sigmas}_{r\theta}(r_0).
\label{eq:radialWaveInterface}
\end{align}
The two constraints in~\eqref{eq:radialWaveFreeSurface} with the four constraints in~\eqref{eq:radialWaveInterface} lead to a homogeneous system of equations of the form
\begin{align}
M \dv = 0, \label{eq:travelingWaveSystem}
\end{align}
for the vector of constants $\dv = [p_I,d,\db_1,\db_2,\db_3,\db_4]^T$.  Nontrivial solutions for $\dv$ exist if $\det(M)=0$, which leads to a dispersion relation between the frequency $\omega$ and the circumferential wave number~$n$ of the form
\begin{align}
{\cal F}(n,\omega) = 0.
\label{eq:dispersion}
\end{align}
For a given $(n,\omega)$-pair satisfying~\eqref{eq:dispersion}, a nontrivial $\dv$ can be found subject to the normalization
\begin{align}
\|\usv(r_0,0,0)\|\sp2 = \us_r(r_0,0,0)^2 + \us_\theta(r_0,0,0)^2 = \us_0\sp2,
\end{align}
where $\us_0$ is an initial displacement of the interface.  Solutions of~\eqref{eq:dispersion} for $(n,\omega)$ and the corresponding constants in $\dv$ are given in Section~\ref{sec:travelingWave} for the case $\us_0=10\sp{-7}$.

\bibliographystyle{elsart-num}
\bibliography{journal-ISI,jwb,henshaw,henshawPapers,fsi}

\end{document}